\documentclass[11pt,a4paper]{article}
\usepackage[margin=2.5cm]{geometry}
\usepackage{amsmath,amssymb,amsthm}
\usepackage[title]{appendix}
\usepackage{booktabs}
\usepackage{graphicx}
\usepackage{enumitem}
\usepackage{xcolor}
\usepackage{framed}
\usepackage[normalem]{ulem}
\usepackage{mathtools}
\usepackage{longtable}
\usepackage[colorlinks=true,linkcolor=blue!50!black,citecolor=blue!50!black,urlcolor=blue!50!black]{hyperref}
\usepackage{orcidlink} 

\usepackage{thmtools}
\usepackage{thm-restate}
\usepackage{hyperref}

\definecolor{olsthoornblue}{HTML}{2166AC}
\definecolor{velocityorange}{HTML}{D95F0E}
\definecolor{barred}{HTML}{CF222E}
\definecolor{barlightred}{HTML}{E99BA0}
\definecolor{cycleblue}{HTML}{0969DA}
\definecolor{deathpurple}{HTML}{8250DF}

\usepackage{tikz}
\tikzset{
  landscape axis/.style={
    ->,
    thin,
    gray!80
  },
  landscape guide/.style={
    dashed,
    thin,
    gray!55
  },
  landscape curve/.style={
    very thick,
    line join=round
  },
  landscape label/.style={
    font=\small,
    fill=white,
    inner sep=2pt
  }
}
\usetikzlibrary{patterns}
\usetikzlibrary{positioning}

\makeatletter
\@ifundefined{IfFormatAtLeastTF}{}{%
  \IfFormatAtLeastTF{2026-06-01}{%
    \renewcommand{\thmt@autorefsetup}{}%
  }{}%
}
\makeatother

\theoremstyle{plain}
\newtheorem{proposition}{Proposition}
\newtheorem{theorem}[proposition]{Theorem}
\newtheorem{lemma}[proposition]{Lemma}
\newtheorem{corollary}[proposition]{Corollary}

\theoremstyle{definition}
\newtheorem{definition}[proposition]{Definition}

\theoremstyle{remark}
\newtheorem{remark}[proposition]{Remark}
\newtheorem{example}[proposition]{Example}

\theoremstyle{definition}
\newtheorem{assumptionsinner}[proposition]{Assumption}

\newenvironment{assumptions}[1][]
  {%
    \begin{framed}%
    \if\relax\detokenize{#1}\relax
      \begin{assumptionsinner}%
    \else
      \begin{assumptionsinner}[#1]%
    \fi
  }
  {%
    \end{assumptionsinner}%
    \end{framed}%
  }

\title{Persistence landscape velocity:\\[3pt]
Closed formulas and a topological probe of quantum criticality}

\author{%
  Francisco Belchí
  \thanks{Corresponding author:
  \href{mailto:francisco.belchi@salle.url.edu}
    {\texttt{francisco.belchi@salle.url.edu}}.
  Francisco Belchí acknowledges support from the
  Departament de Recerca i Universitats de la Generalitat
  de Catalunya and Universitat Ramon Llull.}
    \\[4pt]
  {\small Universitat Ramon Llull, La Salle Campus Barcelona}\\
  {\small Human-Environment Research (HER), Spain}\\
  {\small ORCID:
  \href{https://orcid.org/0000-0001-5863-3343}
    {\texttt{0000-0001-5863-3343}}}
  \and
  Nahuel Statuto
  \thanks{Nahuel Statuto was supported by the Generalitat de Catalunya (AGAUR, 2025 SGR 01604).}
    \\[4pt]
  {\small Departament de Matemàtiques i Informàtica}\\
  {\small Universitat de Barcelona, Barcelona, Spain}\\
  {\small ORCID:
  \href{https://orcid.org/0000-0002-0062-760X}
    {\texttt{0000-0002-0062-760X}}}
}

\begin{document}
\maketitle

\begin{abstract}
We introduce the landscape velocity~\(\xi\), a
scalar measure of variation of a parameter-dependent
persistence landscape.
For families of finite barcodes whose bars are all
born at the fixed filtration value \(0\), we show
that \(\xi\) is the \(L^2\) speed of the
landscape curve where the sorted death times are differentiable and derive an exact closed formula
that requires no persistence landscape at all: the death times of the
barcodes, as functions of the parameter, suffice.
For the velocity of the empirical mean of several
landscapes we also prove a closed formula which,
again, strikingly requires neither the
empirical mean landscape nor any of the
individual landscapes, only the individual barcodes.

As an application, we revisit Olsthoorn's persistent
homology analysis of the transverse-field Ising chain.
Whereas his approach uses different homology degrees and
summaries depending on the
boundary conditions, ours uses \(H_0\) only and a single
rule: the position of the maximum of \(\xi\).
This rule locates the quantum critical point for both
periodic and open boundary conditions, and for both
mutual-information and correlator-based dissimilarities.
The position of the maximum converges to the critical
point faster than all the estimators we can extract from Olsthoorn's work. 
In the 4 scenarios considered (different boundary conditions and different dissimilarity matrices), 
the peak of the normalized velocity becomes increasingly pronounced with system size,
with observed logarithmic growth.
This logarithmic growth is proved for periodic boundaries with correlator dissimilarity under stated hypotheses, identifying landscape velocity as a singular topological
response to quantum criticality.
\end{abstract}
\noindent
\textbf{Keywords:} persistence landscape, landscape velocity, topological data analysis, quantum criticality, transverse-field Ising chain, closed formula, common-birth barcodes, topological response.

\tableofcontents

\section*{Introduction}

Topological data analysis reveals structure in data by asking what persists
across scales. As a filtration parameter \(t\) increases, connected
components merge, holes appear and disappear, and persistent homology
records their lifetimes in a barcode. Barcodes make this structure visible;
persistence landscapes make it accessible to linear methods. By
representing each barcode as a point in a Hilbert space, landscapes provide
a setting in which averaging, statistical inference, and learning become
natural. But this passage comes with an apparent computational cost: to
work with a landscape, one must first construct it.

Now suppose that the data themselves change with a parameter \(g\). There
are then two parameters with distinct roles: \(t\) probes the data across
scales, while \(g\) changes the data being probed. Each value of \(g\)
produces a barcode and hence a landscape; together, these landscapes trace
a curve in Hilbert space. The question is no longer just what persists,
but how rapidly the persistent structure changes. Can we measure this
rate of change without constructing the landscapes? Can such a measure
reveal changes that existing summaries leave less clearly resolved?

We introduce the \emph{landscape velocity}, \(\xi(g)\), a scalar measure of
how rapidly a persistence landscape changes with \(g\). Our main
theoretical result removes the apparent computational cost: for families
of finite barcodes whose bars are all born at filtration value \(0\), we
derive an exact formula involving only the death times and their
derivatives with respect to \(g\). No landscape need be constructed. A
second formula extends this principle to the velocity of an empirical mean
landscape, using only the individual barcodes. As an application, we study
quantum phase transitions in the transverse-field Ising chain, where peaks
in landscape velocity locate the transition by identifying where the
topology reorganizes most rapidly.


\paragraph{Why measure the motion of landscapes rather than barcodes?}
Because we
want to average before measuring motion. A collection of barcodes need
not have a unique Fr\'echet mean~\cite{mileyko-mukherjee-harer,turner-frechet}.
A collection of landscapes, by contrast, has a unique pointwise average
in the Hilbert space \(\mathcal H\) of Section~\ref{subsec:landscapes}
(see Figure~\ref{fig:frechet-vs-landscape}). This makes landscapes a
natural setting for studying parameter-dependent data with repeated
observations or disorder.

Averaging introduces a subtlety: the mean of several landscapes need
not itself be the landscape of any barcode. Our formula for a single
barcode (Theorem~\ref{prop:closed}) therefore does not automatically
extend to the mean. Nevertheless, Theorem~\ref{prop:mean-speed} gives
a second closed formula that computes the velocity of the empirical
mean directly from the individual death times and their parameter
derivatives. Neither the individual landscapes nor their mean need be
constructed. The linear setting makes averaging possible; the closed
formulas make landscape construction unnecessary.

Landscapes also support laws of large numbers, central limit theorems,
confidence bands, and hypothesis tests~\cite{bubenik,chazal-stochastic,bubenik-dlotko,bubenik-properties}.
Statistical inference based on landscape derivatives additionally requires
suitable regularity and moment assumptions. Averaging will be essential
in our forthcoming companion paper on quantum spin systems, where
disorder averages are unavoidable.

\paragraph{Relation to the companion physics paper.}
The present paper develops the mathematics of this new
topological summary, with a topological data analysis audience
in mind. A companion paper, in preparation, is written for a
quantum physics audience. Since there is little overlap between
the two audiences, we have split the technical details into two
separate papers. Here we include only a small application, which
already exhibits advantages over other persistent homology
summaries; the full physical results will appear in the physics
paper.

Section~\ref{sec:application} is this application.
We repeat Olsthoorn's analysis of the ground states of the
transverse-field Ising chain \cite{olsthoorn-entanglement}, with the
mutual-information and correlator distances, and replace his
summaries by \(\xi\).
Olsthoorn adapts two choices to the boundary conditions: the
homology degree and the summary.
With periodic boundaries he uses the position of the maximum of the
\(H_1\) lifetime, and with open boundaries the position of the maximum of the
range of \(H_0\) death times (we also test, in his favour, the position of the
corner that follows this maximum, which is closer to the critical point).
We always use the same rule: \(H_0\) of the Vietoris--Rips
filtration, the landscape velocity \(\xi\) and the position of its
maximum.
This rule locates the transition in both cases, it needs only
\(H_0\), which is much cheaper to compute than \(H_1\), and the
position of the maximum of \(\xi\) converges to the critical point
faster.
With open boundaries neither of Olsthoorn's two candidates converges to the critical point: the maximum of the range of
death times stays near \(g\approx0.75\), and the corner moves away from the critical point as the chain grows.
The advantages are listed in Section~\ref{subsec:app-advantages}.
In addition, under 
Assumption~\ref{hyp:app-mst}, we
mathematically prove that the maximizers of the normalized landscape
velocity converge to \(g_c=1\) for
the periodic chain with correlator dissimilarity. This conditional
result, proved in Section~\ref{subsec:app-log-proof}, complements the
numerical evidence for both boundary conditions and both
dissimilarities.


The exposition unfolds as follows.

Section~\ref{sec:landscapes} builds the foundation: \emph{(1)} We recall persistent homology and barcodes, and motivate landscapes as their statistical embedding;  \emph{(2)} we introduce the landscape velocity and its geometric interpretation as an \(L^2\) speed in the Hilbert space of landscapes, motivated by response functions in quantum physics; \emph{(3)} we review persistent homology summaries and topological applications to quantum systems, including the transverse-field Ising model central to our application.

Section~\ref{sec:common-birth} develops the common-birth
setting and proves the closed formula and its
interpretation as the \(L^2\) speed of the landscape curve.
For Vietoris--Rips \(H_0\) filtrations, the closed formula for $\xi$ 
can be written in terms of minimum-spanning-tree edge weights.
We also show that, when those dissimilarities are locally
Lipschitz in the parameter, the sorted death times are
differentiable almost everywhere.
Thus the differentiability hypothesis of the closed
formula follows almost everywhere from regularity of
the input.
Finally, we discuss what can happen at crossings.
Section~\ref{sec:application} applies \(\xi\) to the
transverse-field Ising chain and compares it with the summaries of
Olsthoorn \cite{olsthoorn-entanglement}. 
Section~\ref{subsec:disorder} distinguishes
three scalar disorder summaries, evaluates their derivative
field overlaps, and proves a closed formula for the velocity of the
empirical mean landscape.

\section{Persistence landscapes and landscape velocity}\label{sec:landscapes}

This section has three parts. We first recall persistent homology
and barcodes (Section~\ref{subsec:ph-barcodes}), then persistence
landscapes and the reasons for working with them rather than with
barcodes (Section~\ref{subsec:landscapes}), and finally we introduce
the new summary studied in this paper, the landscape velocity
(Section~\ref{subsec:velocity}).

\subsection{Persistent homology and barcodes}
\label{subsec:ph-barcodes}

Persistent homology~\cite{elz,zomorodian-carlsson} records how the
homology of a nested family of spaces changes as a scale parameter
increases; see~\cite{edelsbrunner-harer,ghrist-barcodes,carlsson-topology-data,chazal-michel}
for introductions and surveys. We use the following setting.

\begin{definition}[Filtration and barcode]
\label{def:filtration}
A \emph{filtration} of simplicial complexes is a family
\(\mathcal K=\{K(t)\}_{t\in T}\) of simplicial complexes,
indexed by a subset \(T\subseteq\mathbb R\), such that
\(K(s)\subseteq K(t)\) whenever \(s\leq t\).
It is \emph{right-continuous} if
\(K(t)=\bigcap_{s>t}K(s)\) for every \(t\in T\).
Applying simplicial homology \(H_p(-;\mathbb F)\) in degree
\(p\geq0\), with coefficients in a field \(\mathbb F\), yields a
persistence module, namely the vector spaces
\(H_p(K(t);\mathbb F)\) together with the linear maps induced by
the inclusions. When the module is pointwise finite-dimensional
(for instance, when all \(K(t)\) are finite) it decomposes into
interval modules, and the multiset of these intervals is the
\emph{persistence barcode} of \(\mathcal K\) in degree \(p\)
\cite{zomorodian-carlsson,crawley-boevey,edelsbrunner-harer}.
The left endpoint of a bar is its \emph{birth time} and the
right endpoint its \emph{death time}.
\end{definition}

Barcodes are stable: small perturbations of the filtration, for
instance of the input data, produce small changes of the barcode in
the bottleneck distance~\cite{cohen-steiner-stability,chazal-desilva-oudot}.
In this paper we are interested in families of barcodes indexed by
an external parameter~\(g\); such one-parameter families of
persistence diagrams are known as
vineyards~\cite{vineyards}, and their statistics have been studied
in~\cite{munch-frechet}.

Let \(I\subseteq\mathbb{R}\) be an open interval of parameter
values. For each \(g\in I\), consider a persistence barcode,
and denote its finite part by
\[
\mathcal{B}_g^{\mathrm{fin}}
=
\bigl\{[b_m(g),d_m(g))\bigr\}_{m=1}^{M(g)},
\qquad
-\infty<b_m(g)<d_m(g)<\infty.
\]
The bars are counted with multiplicity, meaning that each finite bar contributes one entry to this list, so if
several bars have the same birth and death time, that pair [birth, death) is repeated
once for each such bar.
No restriction on the homological degree or on the
construction of the underlying filtration is imposed
at this stage.

Here \(t\) denotes the filtration parameter, whereas
\(g\) parametrizes the family of barcodes.
Throughout, we construct landscapes from the finite
bars only; any bars with infinite endpoints are omitted.

\subsection{Persistence landscapes}
\label{subsec:landscapes}

\begin{definition}[Persistence landscape]
\label{def:landscape}
Following Bubenik~\cite{bubenik}, for a finite bar
\([b,d)\), define its tent function by
\[
\Lambda_{(b,d)}(t)
=
\max\{0,\min\{t-b,d-t\}\},
\qquad t\in\mathbb{R}.
\]
The persistence landscape associated with
\(\mathcal{B}_g^{\mathrm{fin}}\) is the sequence
\[
\lambda(g)
=
\bigl(\lambda_k(g;\cdot)\bigr)_{k\geq1},
\]
where
\[
\lambda_k(g;t)
=
k\text{-th largest value of }
\bigl\{
\Lambda_{(b_m(g),d_m(g))}(t)
\bigr\}_{m=1}^{M(g)}.
\]
We set \(\lambda_k(g;t)=0\) whenever \(k>M(g)\).
\end{definition}

\begin{figure}[htbp]
\centering
\begin{tikzpicture}[x=0.85cm,y=0.85cm]

  \tikzset{
    texture blue/.style={solid},
    texture red/.style={dash pattern=on 4pt off 2pt},
    texture teal/.style={dash pattern=on 3pt off 1.5pt on 1pt off 1.5pt},
    texture violet/.style={solid},
    texture orange/.style={dash pattern=on 4pt off 2pt},
    texture green/.style={dash pattern=on 3pt off 1.5pt on 1pt off 1.5pt},
  }

  %

  \foreach \x in {0.5,2,3.5,6.5,8,9}{
    \draw[landscape guide] (\x,6) -- (\x,-3);
  }

  \draw[landscape guide] (4.25,2.25) -- (4.25,8.25);
  \draw[landscape guide] (5,1.5) -- (5,7.5);

  \begin{scope}[shift={(0,6)}]

    \draw[landscape axis]
      (-0.4,0) -- (10,0) node[right] {$t$};

    \draw[landscape axis]
      (0,-0.15) -- (0,4.3) node[above] {$\lambda_k$};

    \draw[landscape curve,violet,texture violet]
      (0.5,0)
      -- (3.5,3)
      -- (4.25,2.25)
      -- (5.5,3.5)
      -- (9,0);

    \draw[landscape curve,orange,texture orange]
      (2,0)
      -- (4.25,2.25)
      -- (5,1.5)
      -- (5.75,2.25)
      -- (8,0);

    \draw[landscape curve,green!50!black,texture green]
      (3.5,0)
      -- (5,1.5)
      -- (6.5,0);

    \node[violet,above=5pt] at (5.5,3.4)
      {$\lambda_1$};

    \node[orange,above right=4pt] at (5.2,2.25)
      {$\lambda_2$};

    \node[
      text=green!50!black,
      inner sep=2pt
    ] at (5,1)
      {$\lambda_3$};

  \end{scope}

  \begin{scope}

    \draw[landscape axis]
      (-0.4,0) -- (10,0) node[right] {$t$};

    \draw[landscape axis]
      (0,-0.15) -- (0,4.3) node[above] {$\Lambda$};

    \draw[landscape curve,blue,texture blue]
      (0.5,0) -- (3.5,3) -- (6.5,0);

    \draw[landscape curve,red,texture red]
      (2,0) -- (5.5,3.5) -- (9,0);

    \draw[landscape curve,teal,texture teal]
      (3.5,0) -- (5.75,2.25) -- (8,0);

    \node[
      text=blue,
      above=0pt,
      inner sep=2pt
    ] at (2.75,2.9)
      {$\Lambda_{(b_1,d_1)}$};

    \node[red,above=0pt] at (5.75,3.5)
      {$\Lambda_{(b_2,d_2)}$};

    \node[
      text=teal,
      above=0pt,
      inner sep=2pt
    ] at (5.7,2.25)
      {$\Lambda_{(b_3,d_3)}$};

  \end{scope}

  \begin{scope}[shift={(0,-3)}]

    \draw[landscape axis]
      (-0.4,0) -- (10,0) node[right] {$t$};

    \draw[landscape curve,blue,texture blue]
      (0.5,1.65) -- (6.5,1.65);

    \draw[landscape curve,red,texture red]
      (2,1.05) -- (9,1.05);

    \draw[landscape curve,teal,texture teal]
      (3.5,0.45) -- (8,0.45);

    \node[
      anchor=east,
      text=blue,
      inner sep=2pt
    ] at (0.45,1.65)
      {$[b_1,d_1)$};

    \node[
      anchor=east,
      text=red,
      inner sep=2pt
    ] at (1.95,1.05)
      {$[b_2,d_2)$};

    \node[
      anchor=east,
      text=teal,
      inner sep=2pt
    ] at (3.45,0.45)
      {$[b_3,d_3)$};

    \node[
      anchor=north,
      font=\small
    ] at (0,-0.65) {$0$};

    \foreach \x/\lab/\col in {
      0.5/{b_1}/blue,
      2/{b_2}/red,
      3.5/{b_3}/teal,
      6.5/{d_1}/blue,
      8/{d_3}/teal,
      9/{d_2}/red}{
      \draw (\x,0) -- (\x,-0.08);
      \node[
        below=5pt,
        font=\small,
        text=\col
      ] at (\x,0) {$\lab$};
    }

  \end{scope}

\end{tikzpicture}
\caption{
Illustration of Definition \ref{def:landscape}:
from a barcode to its persistence landscape.
The 3 bars (bottom) determine individual tent functions
(middle), whose pointwise ordered values form the
landscape layers (top).
Each bar \([b_i,d_i)\) determines
\(\Lambda_{(b_i,d_i)}(t)
=\max\{0,\min\{t-b_i,d_i-t\}\}\).
Geometrically, the two sloping sides form an isosceles
right triangle with base vertices \((b_i,0)\) and
\((d_i,0)\), and right-angle apex
\(\bigl(\frac{b_i+d_i}{2},\frac{d_i-b_i}{2}\bigr)\).
At each \(t\), {\color{violet}\(\lambda_1(t)\)}, {\color{orange}\(\lambda_2(t)\)}, and
{\color{green!50!black}\(\lambda_3(t)\)} are the largest, second-largest, and
third-largest tent values, respectively.
Since the tents are not nested, a layer can switch
between tents at their intersections.
For example, {\color{violet}\(\lambda_1\)} follows
portions of {\color{blue}\(\Lambda_{(b_1,d_1)}\)}
and {\color{red}\(\Lambda_{(b_2,d_2)}\)}.
Each data series also has a distinct line texture:
the first interval and first landscape layer are solid,
the second are dashed, and the third are dash-dotted.
All panels use the same horizontal scale.
The continuous dashed guides align birth and death
times across the panels; the two additional dashed
guides connect the tent intersections in the middle
panel to the corresponding layer intersections in the
top panel.
}
\label{fig:general-barcode-to-landscape}
\end{figure}
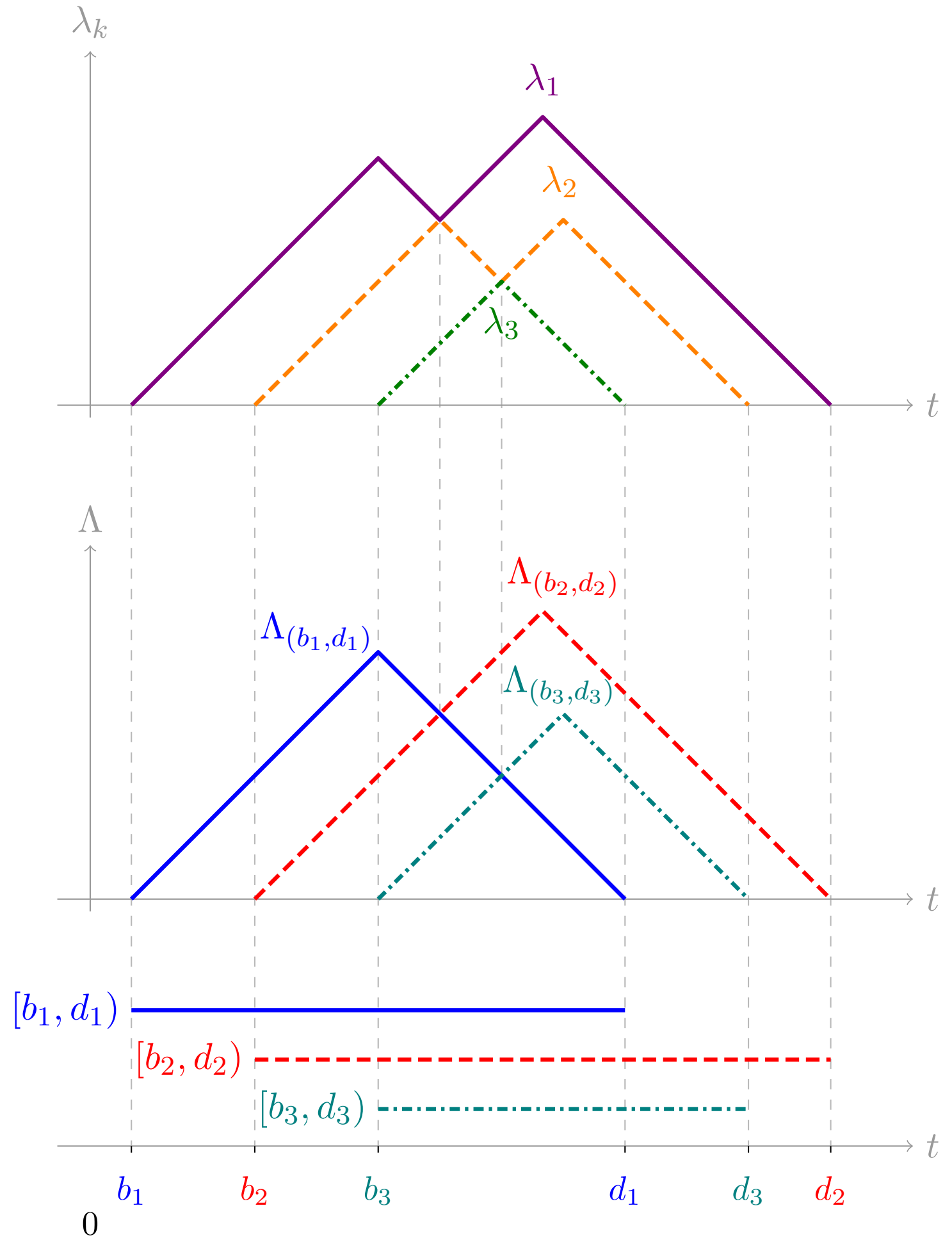

\noindent Figure~\ref{fig:general-barcode-to-landscape}
illustrates the construction of a persistence landscape
from a finite barcode.

\noindent We work in the real Hilbert space
\[
\mathcal H=\ell^2\bigl(\mathbb N;L^2(\mathbb R;\mathbb R)\bigr)
=\left\{u=(u_k)_{k\geq1}:u_k\in L^2(\mathbb R;\mathbb R),\quad
\sum_{k=1}^{\infty}\int_{\mathbb R}|u_k(t)|^2\,dt<\infty\right\}.
\]
Recall that \(L^2(\mathbb R;\mathbb R)\) denotes the space of
(equivalence classes, modulo equality almost everywhere, of)
Lebesgue-measurable functions \(u\colon\mathbb R\to\mathbb R\)
that are square-integrable, that is,
\(\int_{\mathbb R}|u(t)|^2\,dt<\infty\).
Thus an element of \(\mathcal H\) is a sequence of such functions
whose squared \(L^2\) norms are summable.
Its inner product and norm are
\begin{equation}\label{eq:hilbert}
\langle u,w\rangle_{\mathcal H}
=\sum_{k=1}^{\infty}\int_{\mathbb R}u_k(t)w_k(t)\,dt,
\qquad
\|u\|_{\mathcal H}^2
=\sum_{k=1}^{\infty}\int_{\mathbb R}|u_k(t)|^2\,dt.
\end{equation}
Every landscape of a finite barcode belongs to \(\mathcal H\),
since it has finitely many nonzero, compactly supported layers.
Thus distances and inner products compare corresponding layers,
not the sum of all layers before integration.

\subsubsection*{Why landscapes rather than barcodes}
Averaging requires a choice of setting. Barcodes form a metric space,
but not a vector space: there is no ordinary operation of adding two
barcodes and dividing by two. One can instead seek a Fr\'echet mean,
a barcode minimizing the sum of squared distances to the observations~\cite{mileyko-mukherjee-harer,turner-frechet}.
Such a mean, however, need not be unique~\cite{turner-frechet}, and
may depend discontinuously on the data~\cite{munch-frechet}.

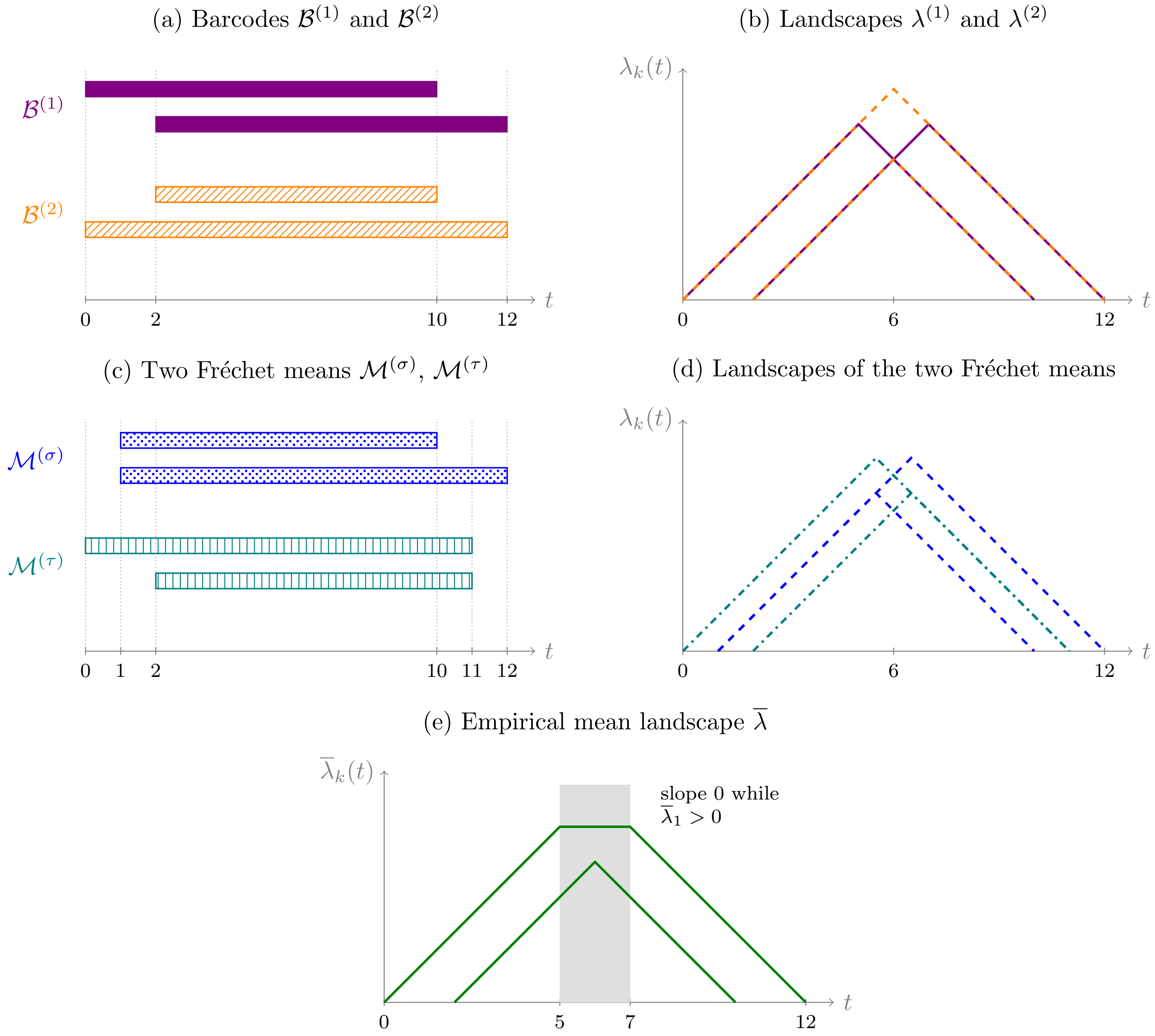
\begin{figure}[htbp]
\centering
\makebox[\linewidth][c]{%
\begin{tikzpicture}[
  x=0.5cm,y=0.5cm,
  lay/.style={line width=1pt},
  ax/.style={->,gray},
  guide/.style={gray!70,densely dotted},
  barB1/.style={fill=violet,draw=violet},
  barB2/.style={draw=orange,line width=0.6pt,pattern=north east lines,pattern color=orange},
  barM1/.style={draw=blue,line width=0.6pt,pattern=crosshatch dots,pattern color=blue},
  barM2/.style={draw=teal,line width=0.6pt,pattern=vertical lines,pattern color=teal},
]
\newcommand{\drawbar}[4]{\draw[#1] (#2,#4-0.22) rectangle (#3,#4+0.22);}
\begin{scope}[shift={(0,0)}]
  \node[font=\small] at (6,8.0) {(a) Barcodes \(\mathcal B^{(1)}\) and \(\mathcal B^{(2)}\)};
  \draw[ax] (0,0)--(12.8,0) node[right,font=\small]{\(t\)};
  \foreach \x in {0,2,10,12}{\draw[gray] (\x,-0.1)--(\x,0.1); \node[below,font=\scriptsize] at (\x,-0.1){\x}; \draw[guide] (\x,0.1)--(\x,6.6);}
  \drawbar{barB1}{0}{10}{6.0}
  \drawbar{barB1}{2}{12}{5.0}
  \drawbar{barB2}{2}{10}{3.0}
  \drawbar{barB2}{0}{12}{2.0}
  \node[left,font=\small,violet] at (-0.3,5.5) {\(\mathcal B^{(1)}\)};
  \node[left,font=\small,orange] at (-0.3,2.5) {\(\mathcal B^{(2)}\)};
\end{scope}
\begin{scope}[shift={(17,0)}]
  \node[font=\small] at (6,8.0) {(b) Landscapes \(\lambda^{(1)}\) and \(\lambda^{(2)}\)};
  \draw[ax] (0,0)--(12.8,0) node[right,font=\small]{\(t\)};
  \draw[ax] (0,0)--(0,6.6) node[left,font=\small]{\(\lambda_k(t)\)};
  \foreach \x in {0,6,12}{\draw[gray] (\x,-0.1)--(\x,0.1); \node[below,font=\scriptsize] at (\x,-0.1){\x};}
  \draw[lay,violet] (0,0)--(5,5)--(6,4)--(7,5)--(12,0);
  \draw[lay,violet] (2,0)--(6,4)--(10,0);
  \draw[lay,orange,dashed] (0,0)--(6,6)--(12,0);
  \draw[lay,orange,dashed] (2,0)--(6,4)--(10,0);
\end{scope}
\begin{scope}[shift={(0,-10)}]
  \node[font=\small] at (6,8.0) {(c) Two Fr\'echet means \(\mathcal M^{(\sigma)}\), \(\mathcal M^{(\tau)}\)};
  \draw[ax] (0,0)--(12.8,0) node[right,font=\small]{\(t\)};
  \foreach \x in {0,1,2,10,11,12}{\draw[gray] (\x,-0.1)--(\x,0.1); \draw[guide] (\x,0.1)--(\x,6.6);}
  \foreach \x in {0,2,10,12}{\node[below,font=\scriptsize] at (\x,-0.1){\x};}
  \foreach \x in {1,11}{\node[below,font=\scriptsize] at (\x,-0.1){\x};}
  \drawbar{barM1}{1}{10}{6.0}
  \drawbar{barM1}{1}{12}{5.0}
  \drawbar{barM2}{0}{11}{3.0}
  \drawbar{barM2}{2}{11}{2.0}
  \node[left,font=\small,blue] at (-0.3,5.5) {\(\mathcal M^{(\sigma)}\)};
  \node[left,font=\small,teal] at (-0.3,2.5) {\(\mathcal M^{(\tau)}\)};
\end{scope}
\begin{scope}[shift={(17,-10)}]
  \node[font=\small] at (6,8.0) {(d) Landscapes of the two Fr\'echet means};
  \draw[ax] (0,0)--(12.8,0) node[right,font=\small]{\(t\)};
  \draw[ax] (0,0)--(0,6.6) node[left,font=\small]{\(\lambda_k(t)\)};
  \foreach \x in {0,6,12}{\draw[gray] (\x,-0.1)--(\x,0.1); \node[below,font=\scriptsize] at (\x,-0.1){\x};}
  \draw[lay,blue,dashed] (1,0)--(6.5,5.5)--(12,0);
  \draw[lay,blue,dashed] (1,0)--(5.5,4.5)--(10,0);
  \draw[lay,teal,dash dot] (0,0)--(5.5,5.5)--(11,0);
  \draw[lay,teal,dash dot] (2,0)--(6.5,4.5)--(11,0);
\end{scope}
\begin{scope}[shift={(8.5,-20)}]
  \node[font=\small] at (6,8.0) {(e) Empirical mean landscape \(\overline\lambda\)};
  \draw[ax] (0,0)--(12.8,0) node[right,font=\small]{\(t\)};
  \draw[ax] (0,0)--(0,6.6) node[left,font=\small]{\(\overline\lambda_k(t)\)};
  \foreach \x in {0,5,7,12}{\draw[gray] (\x,-0.1)--(\x,0.1); \node[below,font=\scriptsize] at (\x,-0.1){\x};}
  \fill[gray!25] (5,0) rectangle (7,6.2);
  \draw[lay,green!50!black] (0,0)--(5,5)--(7,5)--(12,0);
  \draw[lay,green!50!black] (2,0)--(6,4)--(10,0);
  \node[font=\scriptsize,align=left,anchor=west] at (7.6,5.6) {slope \(0\) while\\\(\overline\lambda_1>0\)};
\end{scope}
\end{tikzpicture}
}\par
\caption{
\textbf{Averaging barcodes and averaging landscapes.}
(a) The barcodes \({\color{violet}\mathcal B^{(1)}}\) (solid \textcolor{violet}{violet}
bars) and \({\color{orange}\mathcal B^{(2)}}\) (\textcolor{orange}{orange} bars hatched
with diagonal lines).
The dotted guides mark the endpoints \(0\), \(2\), \(10\), \(12\):
the matching \({\color{blue}\sigma}\) pairs bars with equal death times, and the
matching \({\color{teal}\tau}\) pairs bars with equal birth times.
(b) The landscapes of \({\color{violet}\mathcal B^{(1)}}\) (solid
\textcolor{violet}{violet}) and \({\color{orange}\mathcal B^{(2)}}\) (dashed
\textcolor{orange}{orange}); their second layers coincide.
(c) The two Fr\'echet means: \({\color{blue}\mathcal M^{(\sigma)}}\)
(\textcolor{blue}{blue} bars filled with dots), obtained from
\({\color{blue}\sigma}\) by averaging birth times, and \({\color{teal}\mathcal M^{(\tau)}}\)
(\textcolor{teal}{teal} bars hatched with vertical lines), obtained
from \({\color{teal}\tau}\) by averaging death times.
(d) The landscapes of \({\color{blue}\mathcal M^{(\sigma)}}\) (dashed
\textcolor{blue}{blue}) and \({\color{teal}\mathcal M^{(\tau)}}\) (dash-dotted
\textcolor{teal}{teal}) differ from each other and from
\({\color{green!50!black}\overline\lambda}\).
(e) The unique empirical mean landscape
\({\color{green!50!black}\overline\lambda}=\frac12({\color{violet}\lambda^{(1)}}+{\color{orange}\lambda^{(2)}})\)
(solid \textcolor{green!50!black}{dark green}).
On the \textcolor{gray}{grey} strip its first layer is positive
with slope \(0\), so \({\color{green!50!black}\overline\lambda}\) is not the landscape of any
barcode.
All horizontal axes represent the filtration value \(t\), and all
panels share the same horizontal scale.
}
\label{fig:frechet-vs-landscape}
\end{figure}

Landscapes admit a simpler operation: average corresponding layers
pointwise. The result is a unique element of \(\mathcal H\), although
it need not be the landscape of any barcode. Figure~\ref{fig:frechet-vs-landscape}
displays both phenomena: barcode averaging can give different answers,
while landscape averaging gives one answer that may lie outside the
collection of barcode landscapes. This distinction is essential in
Section~\ref{subsec:disorder}, where we measure the velocity of the
average itself. The single-barcode formula of Theorem~\ref{prop:closed}
does not automatically apply to this average; Theorem~\ref{prop:mean-speed}
provides the required extension.

\subsubsection*{Statistics with persistence landscapes}
Because landscapes live in the separable Hilbert space
\(\mathcal H\) (and, more generally, in the Banach spaces
\(\ell^p(\mathbb N;L^p(\mathbb R))\)), they inherit the tools of
probability and statistics in such spaces~\cite{bubenik,bubenik-properties}.
The empirical mean landscape satisfies a strong law of large
numbers and a central limit theorem, which yield confidence
intervals and hypothesis tests for real-valued functionals of
landscapes~\cite{bubenik}.
Uniform central limit theorems and bootstrap confidence bands for the
population mean landscape were obtained in~\cite{chazal-stochastic,chazal-bootstrap},
and averaging the landscapes of many subsamples provides a stable and
computationally cheaper estimator for large data
sets~\cite{chazal-subsampling}.
Two-sample and permutation tests can be built from distances between
mean landscapes~\cite{bubenik,robinson-turner,kovacev-nikolic},
landscapes are stable with respect to the bottleneck
distance~\cite{bubenik}, and they can be fed directly into
machine-learning methods such as support vector machines and
principal component analysis~\cite{bubenik,kovacev-nikolic}, or
used as differentiable layers in neural networks~\cite{pllay}.
For the related confidence sets for persistence diagrams,
see~\cite{fasy-confidence}; for other functional summaries, such as
silhouettes, see~\cite{chazal-stochastic,berry-functional}; for
multiparameter landscapes, see~\cite{vipond-multiparameter-landscapes}.
Efficient software for computing landscapes, their averages, and
distances between them is available~\cite{bubenik-dlotko}.

\subsection{The landscape velocity: a new topological summary}
\label{subsec:velocity}

\paragraph{Motivation from physics.}
In statistical and quantum physics, phase transitions are detected
through response functions, that is, derivatives of observables with
respect to a control parameter; for instance, the magnetic
susceptibility is the derivative of the magnetization with respect to
the applied field. Response functions peak, or diverge in the
thermodynamic limit, where the system reorganizes most rapidly.
A particularly close analogue of our construction is the fidelity
susceptibility of a family of ground states
\(|\psi_0(g)\rangle\)~\cite{zanardi-paunkovic,you-li-gu,gu-fidelity}:
it is the squared norm of the component of
\(\partial_g|\psi_0(g)\rangle\) orthogonal to \(|\psi_0(g)\rangle\),
that is, the squared speed of the ground state with respect to the
quantum metric~\cite{provost-vallee}.
The landscape velocity mimics this construction, with the
persistence landscape in the role of the state: \(\xi(g)^2\) is the
squared norm of the derivative \(\partial_g\lambda(g)\) in the
Hilbert space \(\mathcal H\), so \(\xi\) is large when the whole
landscape reorganizes rapidly as \(g\) varies.
Two features of the physical setting guide our choices.
First, in disordered systems physical predictions are disorder
averages; the summary must therefore be compatible with averaging,
which is why we work with landscapes rather than barcodes
(Section~\ref{subsec:landscapes}).
Second, the physically meaningful response of a disorder average is
the derivative of the averaged observable, which here is the velocity
of the empirical mean landscape;
Section~\ref{subsec:disorder} computes it
directly from the individual barcodes.

To quantify how the persistence landscape varies with
the parameter \(g\), we introduce the following scalar
measure of variation, which we call the
\emph{landscape velocity}.

\begin{definition}[Landscape velocity]
\label{def:velocity}
Let \(g\in I\) be a parameter value such that, for every
\(k\geq1\), the derivative
\(\partial_g\lambda_k(g;t)\) exists for almost every
\(t\in\mathbb{R}\), and
\[
\sum_{k=1}^{\infty}
\int_{\mathbb{R}}
\bigl|\partial_g\lambda_k(g;t)\bigr|^2\,dt
<\infty.
\]
We define the landscape velocity at \(g\) by
\[
\xi(g)
=
\left(
\sum_{k=1}^{\infty}
\int_{\mathbb R}
\bigl(\partial_g\lambda_k(g;t)\bigr)^2\,dt
\right)^{1/2}.
\]
\end{definition}

No normalization constant enters the definition. If, in some
computation, a per-bar quantity is preferable (for instance, when
comparing barcodes with different numbers \(M\) of finite bars, such as in Section \ref{sec:application}),
we shall write it explicitly, e.g.\ \(\xi/\sqrt M\); since this is
multiplication by a constant, nothing else changes.
Throughout, \(\xi(g)\) is a nonnegative scalar, not an element
of \(\mathcal H\). 
When the curve is differentiable in \(\mathcal H\),
we denote by \(\partial_g\lambda(g)\) its vector-valued derivative,
with components \(\partial_g\lambda_k(g;t)\).
The definition above uses pointwise
layer derivatives; by itself it does not assert differentiability
in the \(\mathcal H\) norm. Theorem~\ref{prop:l2-speed}
establishes this stronger property for the common-birth families
considered below.

We use \(\partial_g f\) to denote differentiation with
respect to \(g\), and \(\partial_g f(g_0)\) to denote
the derivative evaluated at \(g=g_0\).
In particular,
\[
\partial_g\lambda_k(g_0;t)
=
\left.
\frac{\partial\lambda_k(g;t)}{\partial g}
\right|_{g=g_0},
\]
with \(t\) held fixed.
We do not use primes for derivatives with respect to \(g\):
in figures whose horizontal axis is the filtration value \(t\),
a prime could be mistaken for a derivative with respect to \(t\).
One-sided derivatives are written
\(\displaystyle \partial_g^\pm f(g_0)=\lim_{h\to0^\pm}\frac{f(g_0+h)-f(g_0)}{h}\).

\begin{remark}[Terminology: derivative, velocity, and speed]
\label{rem:terminology}
We distinguish three objects.
The \emph{derivative} \(\partial_g\lambda(g)\in\mathcal H\) is a
vector, defined where the landscape curve
\(g\mapsto\lambda(g)\) is differentiable in \(\mathcal H\).
The \emph{landscape velocity} \(\xi(g)\) of
Definition~\ref{def:velocity} is a nonnegative scalar; we keep the
word \emph{velocity}, used in the physics literature that motivates
this work, as the name of this statistic.
The \emph{metric speed} of the curve is
\(\lim_{h\to0}\|\lambda(g+h)-\lambda(g)\|_{\mathcal H}/|h|\),
whenever the limit exists.
Where the curve is differentiable in \(\mathcal H\), the three are
related by
\(\xi(g)=\|\partial_g\lambda(g)\|_{\mathcal H}=\) metric speed
(Theorem~\ref{prop:l2-speed}); at a crossing of bars the
derivative may fail to exist while the metric speed still exists
(Remark~\ref{rem:metric-speed}).
\end{remark}

\subsection{State of the art}
\label{subsec:state-of-art}

This subsection places the landscape velocity in the context of
two bodies of work: summaries of persistent homology, and
applications of topological data analysis to quantum physics.

\subsubsection{Persistent homology summaries}
\label{subsubsec:sota-summaries}

Barcodes and persistence diagrams are complete invariants of
pointwise finite-dimensional persistence modules
\cite{crawley-boevey,zomorodian-carlsson}, and they are stable
with respect to perturbations of the input
\cite{cohen-steiner-stability,chazal-proximity,chazal-structure}.
However, the space of persistence diagrams with the bottleneck or
Wasserstein distance is not a vector space: Fr\'echet means exist
but need not be unique \cite{mileyko-mukherjee-harer,turner-frechet}.
This has motivated a large family of \emph{summaries}, which map
barcodes to a vector space or a Hilbert space; see the surveys of
Berry et al.\ \cite{berry-functional} on functional summaries and
of Ali et al.\ \cite{ali-survey} on vectorizations.

Among the functional summaries, the persistence landscape of
Bubenik \cite{bubenik,bubenik-properties} is the one used in this
paper, and it comes with a software toolbox
\cite{bubenik-dlotko}.
Related functional summaries include the weighted silhouettes of
Chazal et al.\ \cite{chazal-stochastic}, the persistence curves
of Chung and Lawson \cite{chung-lawson}, which contain Betti curves
as a special case, and the stable persistence images of Adams
et al.\ \cite{adams-images}.
Coordinate-based summaries include the ring of algebraic
functions on barcodes of Adcock, Carlsson and Carlsson
\cite{adcock-carlsson-carlsson} and the tropical coordinates of
Kali\v{s}nik \cite{kalisnik-tropical}.
Scalar summaries include persistent entropy
\cite{chintakunta-entropy}, whose stability was studied by
Atienza, Gonzalez-D\'iaz and Soriano-Trigueros
\cite{atienza-entropy}.
Kernel methods avoid an explicit vectorization: examples are the
scale-space kernel of Reininghaus et al.\ \cite{reininghaus-kernel},
the persistence weighted Gaussian kernel of Kusano, Fukumizu and
Hiraoka \cite{kusano-kernel}, and the sliced Wasserstein kernel of
Carri\`ere, Cuturi and Oudot \cite{carriere-sliced}.
Finally, learnable layers such as PersLay \cite{perslay} and the
landscape-based layer PLLay \cite{pllay} incorporate summaries
into neural networks.
For multiparameter persistence
\cite{botnan-lesnick,lesnick-wright}, the landscape construction
has been extended by Vipond \cite{vipond-multiparameter-landscapes}.

The statistical theory of landscapes is particularly well
developed: it includes laws of large numbers and central limit
theorems \cite{bubenik}, uniform convergence and confidence bands
\cite{chazal-stochastic}, bootstrap and subsampling methods
\cite{chazal-bootstrap,chazal-subsampling}, and hypothesis tests
\cite{robinson-turner}; see also \cite{fasy-confidence} for
confidence sets for diagrams and \cite{chazal-michel} for an
introduction.
Landscapes have also been used as inputs for statistical learning,
for example in protein binding analysis \cite{kovacev-nikolic}.

Parameter-dependent and time-dependent barcodes have also been
studied.
Vineyards track barcodes along one-parameter families of
filtrations \cite{vineyards}; CROCKER plots display Betti numbers as
functions of both the filtration and the time parameter
\cite{topaz-crocker}; and Munch et al.\ define probabilistic
Fr\'echet means for time-varying diagrams \cite{munch-frechet}.
Leygonie, Oudot and Tillmann develop a differential calculus for
maps into the space of barcodes \cite{leygonie-calculus}, and
Gameiro, Hiraoka and Obayashi use the derivatives of birth and
death times to continue point clouds along prescribed paths in
the space of diagrams \cite{gameiro-continuation}.
The landscape velocity is close in spirit to these works, since it
differentiates a summary with respect to an external parameter.
It differs from them in that the derivative is taken in the
Hilbert space \(\mathcal H\) of landscapes, so it is a linear
object that commutes with averaging, and this allows us to give a
closed formula for the velocity of an
empirical mean landscape
(Theorem~\ref{prop:mean-speed}).
The scalar velocity with respect to an external parameter differs from the square-root-velocity (SRVF) parameterization used in prior work to analyze persistence-landscape curves as geometric objects \cite{matuk-kurtek-bharath}.

\subsubsection{Topological data analysis in quantum physics}
\label{subsubsec:sota-quantum}

Topological data analysis has been used increasingly to detect and
characterize phase transitions; see the review of Leykam and
Angelakis \cite{leykam-angelakis}.
In classical statistical mechanics, persistent homology has been
used to analyse the phase transitions of the mean-field XY model and
of the \(\phi^4\) lattice model \cite{donato-phase}, to find hidden order in frustrated spin models
\cite{olsthoorn-hidden}, to build interpretable order parameters
\cite{cole-loges-shiu}, and to locate the transition in
two-dimensional XY models quantitatively \cite{sale-xy}.
Sehayek and Melko studied \(\mathbb Z_2\) gauge theories with
persistent homology \cite{sehayek-melko}.
Sun et al.\ combined network science and persistent homology to
characterize networks built from Monte Carlo snapshots of the
two-dimensional classical Ising model across its phase transition
\cite{sun-bianconi}.

In quantum many-body physics, Tran, Chen and Hasegawa build point
clouds from correlations in states.
They detect the Berezinskii--Kosterlitz--Thouless transition of
the classical XY model and the quantum phase transitions of the
transverse Ising and Bose--Hubbard models without prior knowledge of
the phases \cite{tran-chen-hasegawa}.
Tirelli and Costa learn quantum phase transitions of interacting
fermions from persistent homology of quantum Monte Carlo snapshots
\cite{tirelli-costa}.
He et al.\ applied persistent homology to localization in a
generalized Aubry--Andr\'e--Harper model \cite{he-aah}.
Spitz et al.\ found self-similar behaviour in far-from-equilibrium
quantum many-body dynamics \cite{spitz-self-similar}.
Komalan proposed ``quantum barcodes'' for symmetry-protected
topological phases and illustrated them on the SSH model
\cite{komalan-quantum-barcodes}.
The entanglement structure of multiqubit states has also been
studied with persistent homology
\cite{mengoni-multiqubit,hamilton-leditzky}.
Quantum algorithms for estimating Betti numbers were proposed by
Lloyd, Garnerone and Zanardi \cite{lloyd-quantum-algorithms}.

The transverse-field Ising model (TFIM), which is the model of the
application in this paper, has been studied in several of these
works.
Olsthoorn builds a Vietoris--Rips filtration on the sites of a spin
chain, with a dissimilarity that decreases with the quantum mutual
information between sites, and computes it on ground states of the
TFIM \cite{olsthoorn-entanglement}.
With periodic boundary conditions, a persistent \(H_1\) bar appears
at the critical point \(g=1\) (Olsthoorn's \(h=1\)).
With open boundary conditions, the transition is detected in
\(H_0\): since all \(H_0\) bars are born at \(0\), the information
lies in the death times, and the difference between the maximum and
minimum death has a peak below \(g=1\), which he relates to the transition (we also examine the corner that follows the peak, which is closer to \(g=1\)) \cite{olsthoorn-entanglement}.
The same paper studies the many-body localization transition of
the disordered XXZ chain by averaging over disorder realizations
\cite{olsthoorn-entanglement}.
This is precisely the common-birth \(H_0\) setting of
Section~\ref{sec:common-birth}, in which our closed formula
(Theorem~\ref{prop:closed}) expresses the landscape velocity
through the death times and their derivatives only, and in which
Theorem~\ref{prop:mean-speed} handles disorder averages.
Closely related are the complex networks of quantum mutual
information, whose network measures locate the critical points of
the quantum Ising and Bose--Hubbard models \cite{valdez-networks}
and describe thermal states of the Ising chain \cite{sundar-networks}.
Related network descriptions have been given for the entanglement
structure of quantum spin chains \cite{sokolov-networks} and for wave
functions of quantum many-body systems \cite{mendes-santos-network}.
The TFIM itself is a textbook model: it can be solved through
the Jordan--Wigner mapping to free fermions
\cite{lieb-schultz-mattis,pfeuty}, and its quantum phase transition is
described in the monographs \cite{sachdev,dutta-tfim}.
Its disordered version, the random transverse-field Ising chain, was
analysed by Fisher \cite{fisher-prl,fisher-prb}.
Fidelity and the fidelity susceptibility are geometric detectors of
the same transitions \cite{zanardi-paunkovic,you-li-gu,gu-fidelity},
based on the quantum geometric tensor \cite{provost-vallee}.
A related approach, combining fidelity-based measures with persistent homology, has been explored in the study of quantum critical phenomena \cite{fidelity-persistence}; the two frameworks complement each other, with landscape velocity offering an alternative topological route that operates on dissimilarity data rather than directly on quantum states.

\section{The closed form for fixed-birth-zero
  \texorpdfstring{\(H_0\)}{H0} barcodes}
\label{sec:common-birth}
In this section, we will prove that when persistent homology is being computed for $H_0$, for filtrations of type Vietoris--Rips, \v{C}ech or alpha complex filtrations, among others, we can compute the velocity $\xi$ without the need to even compute landscapes, as we can find a closed formula (Theorem \ref{prop:closed}) that depends only on the death times of the bars in the $H_0$ barcode.

\begin{assumptions}[Common-birth \(H_0\) setting]
\label{ass:common-birth}
Let \(X=\{1,\dots,L\}\), with \(L\geq2\), and let
\(I\subseteq\mathbb R\) be an open interval.
For each \(g\in I\), let
\(\mathcal{K}_g=\{K_g(t)\}_{t\geq0}\)
be an increasing, right-continuous filtration
(Definition~\ref{def:filtration}) of simplicial complexes on the
fixed vertex set \(X\), such that:
\begin{enumerate}[label=(A\arabic*),nosep]
\item\label{ass:isolated} \(K_g(0)\) consists of the \(L\) isolated
vertices of \(X\);
\item\label{ass:connected} \(K_g(t)\) is connected for all
sufficiently large \(t\);
\item\label{ass:field} ordinary \(H_0\) is computed with
coefficients in a fixed field, and landscapes are constructed from
the finite bars only.
\end{enumerate}
\end{assumptions}
The common-birth condition (item~(A1) above) is essential: for barcodes with variable birth times, the tent-nesting structure collapses, and the closed formula no longer applies. As seen in Section \ref{subsec:3filtrations}, examples of filtrations that satisfy these conditions are Vietoris--Rips, \v{C}ech and alpha complex filtrations.

Every connected component is therefore born at \(0\),
and components can subsequently merge but cannot split.
Since the filtration eventually becomes connected,
its \(H_0\) barcode consists of one infinite bar
\([0,\infty)\) and \(L-1\) finite bars.
Write the finite part as
\[
\mathcal{B}_g^{\mathrm{fin}}
=
\bigl\{[0,d_k(g))\bigr\}_{k=1}^{L-1},
\qquad
d_1(g)\geq d_2(g)\geq\cdots\geq d_{L-1}(g)>0.
\]
The death times are listed with multiplicity, meaning that each finite bar contributes one entry to this list, so if
several bars have the same death time, that value is repeated
once for each such bar.

More generally, we shall use the following setting, which does
not refer to any filtration.
\begin{assumptions}[Common-birth family]
\label{ass:common-birth-family}
There is a fixed integer \(M\geq1\) such that, for every
\(g\in I\), the finite part of the barcode is
\(\mathcal{B}_g^{\mathrm{fin}}=\{[0,d_k(g))\}_{k=1}^{M}\), with
\(d_1(g)\geq\cdots\geq d_M(g)>0\), and the landscape is constructed
from these finite bars only.
\end{assumptions}
\noindent By the discussion above,
~Assumption~\ref{ass:common-birth} imply
~Assumption~\ref{ass:common-birth-family} with \(M=L-1\). 

Throughout the rest of section \ref{sec:common-birth}, the landscape is constructed
from these finite bars only. Thus
\[
\xi^2(g)
=
\sum_{k=1}^{L-1}
\int_0^\infty
\bigl(\partial_g\lambda_k(g;t)\bigr)^2\,dt.
\]
The integral over \(\mathbb{R}\) in the general
definition reduces to \([0,\infty)\), since the
landscape vanishes for \(t<0\).

Everything we prove in this paper for ~Assumption~\ref{ass:common-birth} applies also to the more general setting of
~Assumption~\ref{ass:common-birth-family}

\begin{remark}[Positive death times]
\label{rem:right-continuity}
For an increasing filtration of simplicial complexes, right-continuity ensures that no edge appears immediately after $0$ without being present at $0$ itself. Together with the assumption that $K_g(0)$ consists of isolated vertices, this ensures that all finite death times are strictly positive: the value of $\varepsilon$ such that no edge is present for $0\leq t\leq\varepsilon$ is positive by the finiteness of possible edges.
\end{remark}

\subsection{Filtrations covered by the common-birth setting}
\label{subsec:3filtrations}

We describe several constructions of filtrations
\[
\mathcal K=\{K(t)\}_{t\geq0}
\]
satisfying 
Assumption~\ref{ass:common-birth}.
For clarity, we describe each construction for a single
choice of its input data, without displaying the family
parameter \(g\).
A parameter-dependent family
\[
\mathcal K_g=\{K_g(t)\}_{t\geq0},
\qquad g\in I,
\]
is obtained by allowing the input data to depend on \(g\)
and applying the construction separately at each
parameter value.
Provided the stated assumptions on the input data hold
for every \(g\), each filtration \(\mathcal K_g\)
satisfies the common-birth hypotheses.
Regularity with respect to \(g\), including the
differentiability of the sorted death times required
in Theorem~\ref{prop:closed}, is a separate condition.

Throughout, \(X=\{1,\dots,L\}\), with \(L\geq2\), is a
fixed labeled vertex set, and \(t\in[0,\infty)\) is the
filtration parameter.
We describe a simplicial complex on \(X\) as a collection
of nonempty subsets of \(X\), called \emph{simplices}.
A subset with one vertex represents a point, a subset
with two vertices represents an edge, a subset with
three vertices represents a filled triangle, and so on.
The collection must contain every nonempty subset of
each of its simplices: for example, including a filled
triangle requires including its three edges and three
vertices.
Some conventions also include the empty set as a simplex;
we omit it from our lists, since this does not affect
the connected components or the ordinary \(H_0\)
persistence considered here.

\begin{enumerate}[label=(\roman*)]

\item \emph{Vietoris--Rips filtrations.}
Let
\[
D:X\times X\longrightarrow[0,\infty)
\]
be a finite-valued symmetric dissimilarity satisfying
\[
D(i,i)=0,
\qquad
D(i,j)>0
\quad\text{for }i\neq j.
\]
No triangle inequality is needed.
Define
\[
K(t)
=
\mathrm{VR}_t(X,D)
=
\left\{
\varnothing\neq\sigma\subseteq X:
D(i,j)\leq t
\text{ for all }i,j\in\sigma
\right\}.
\]

Equivalently, \(K(t)\) is the clique or flag complex of the
graph with vertex set \(X\) and edges
\[
\bigl\{\{i,j\}:i\neq j,\ D(i,j)\leq t\bigr\}.
\]
Here the clique complex of a graph consists of all
nonempty vertex subsets whose distinct vertices
are pairwise joined by edges.

This is the diameter-threshold convention for the
Vietoris--Rips filtration.
When \(D\) is a metric, it gives the usual
Vietoris--Rips complex; the definition itself does
not require the triangle inequality
\cite{edelsbrunner-harer,chazal-desilva-oudot}.
The complex goes back to Vietoris \cite{vietoris}.

Every singleton belongs to \(K(0)\), and no edge
does, so \(K(0)\) consists of \(L\) isolated vertices.
Moreover,
\[
K(t)
=
\{\varnothing\neq\sigma\subseteq X\}
\qquad
\text{whenever }
t\geq\max_{i,j\in X}D(i,j).
\]
Thus the filtration eventually becomes a full simplex
and, in particular, becomes connected.

\item \emph{\v{C}ech filtrations.}
Let
\[
p:X\longrightarrow\mathbb{R}^{n}
\]
be injective, so that \(p(1),\dots,p(L)\) are
distinct Euclidean points.
For \(i\in X\) and \(t\geq0\), define the closed ball
\[
B_{i}(t)
=
\left\{
z\in\mathbb{R}^{n}:
\|z-p(i)\|\leq t
\right\}.
\]
The \v{C}ech complex at radius \(t\) is
\[
K(t)
=
\check{C}(t)
=
\left\{
\varnothing\neq\sigma\subseteq X:
\bigcap_{i\in\sigma}B_{i}(t)\neq\varnothing
\right\}.
\]
Thus \(\check{C}(t)\) is the nerve of the indexed
family of closed balls
\(\{B_{i}(t)\}_{i\in X}\)
\cite{edelsbrunner-harer,chazal-desilva-oudot}.
Since finite intersections of balls are convex, the nerve
theorem \cite{borsuk-nerve,bjorner-topological} shows that
\(\check{C}(t)\) is homotopy equivalent to the union of the balls.

At \(t=0\), the balls are distinct singleton sets,
so \(K(0)\) consists of isolated vertices.
Let
\[
R
=
\max_{i\in X}\|p(i)-p(1)\|.
\]
For every \(t\geq R\), the point \(p(1)\) belongs
to all the balls \(B_{i}(t)\).
Hence \(K(t)\) is a full simplex and is connected.

\item \emph{Unweighted alpha filtrations.}
Let \(p:X\to\mathbb{R}^{n}\) be injective as above.
For each \(i\in X\), define its closed Voronoi cell by
\[
V_{i}
=
\left\{
z\in\mathbb{R}^{n}:
\|z-p(i)\|
\leq
\|z-p(j)\|
\text{ for every }j\in X
\right\}.
\]
Its restriction at radius \(t\) is
\[
A_{i}(t)
=
V_{i}\cap B_{i}(t).
\]
The alpha complex at radius \(t\) is the nerve of
these restricted Voronoi cells:
\[
K(t)
=
\operatorname{Alpha}(t)
=
\left\{
\varnothing\neq\sigma\subseteq X:
\bigcap_{i\in\sigma}A_{i}(t)\neq\varnothing
\right\}.
\]
This is the radius-parameterized, unweighted alpha
filtration
\cite{edelsbrunner-harer,carlsson-carlsson-alpha}.
Alpha shapes were introduced by Edelsbrunner and M\"ucke
\cite{edelsbrunner-mucke}, and the description of the alpha complex
as the nerve of restricted Voronoi cells is due to Edelsbrunner
\cite{edelsbrunner-union}.
The term \emph{unweighted} means that the input points
carry no additional weights: the construction uses
ordinary Euclidean Voronoi cells and balls with a
common radius \(t\).
In particular,
\[
A_i(0)=\{p(i)\},
\]
so distinct input points give isolated vertices at
filtration value \(0\).
Weighted alpha filtrations instead use power-distance
functions associated with weighted points; they need
not satisfy this initial condition and are not
included here without additional assumptions.

Since each \(V_{i}\) is fixed as \(t\) varies
and the balls increase with \(t\), the complexes
\(\operatorname{Alpha}(t)\) form a filtration.

For every \(t\geq0\), the restricted cells cover
the same set as the balls:
\[
\bigcup_{i\in X}A_{i}(t)
=
\bigcup_{i\in X}B_{i}(t).
\]
Indeed, a point in the union of the balls lies in
the restricted Voronoi cell of any nearest input
point.
All nonempty finite intersections of restricted
cells are convex.
Consequently, the nerve theorem \cite{borsuk-nerve,bjorner-topological} identifies the
homotopy type of \(\operatorname{Alpha}(t)\)
with that of the union of the balls
\cite{edelsbrunner-harer,carlsson-carlsson-alpha}.

For sufficiently large \(t\), the balls have a
common point, and their union is connected.
It follows that \(\operatorname{Alpha}(t)\)
is eventually connected.
Unlike the \v{C}ech complex, it need not eventually
be a full simplex.

The alpha and \v{C}ech filtrations have isomorphic
persistent homology when parameterized by the
same radius \cite{edelsbrunner-union,bauer-edelsbrunner}.
If squared radius \(s\) is used instead, the
corresponding filtration is
\[
\operatorname{Alpha}^{\mathrm{sq}}(s)
=
\operatorname{Alpha}(\sqrt{s}),
\qquad s\geq0.
\]
The death times in the closed formula 
in Theorem \ref{prop:closed}
must be
expressed in the parameter actually used.

\item \emph{Edge-threshold filtrations of connected graphs.}
Let \(G=(X,E)\) be a finite connected simple
undirected graph.
Let
\[
w:E\longrightarrow(0,\infty)
\]
assign a finite positive weight to every edge.
Define
\[
E(t)
=
\{e\in E:w(e)\leq t\},
\qquad
G(t)=(X,E(t)).
\]
Viewing each graph as a one-dimensional simplicial
complex, define
\[
K(t)
=
\bigl\{\{i\}:i\in X\bigr\}
\cup E(t).
\]
Equivalently, this is the sublevel filtration 
\(\bigl\{f^{-1}\big( [0,t] \big) \bigr\}_{t\geq0}\)
of
the simplex-wise function
\[
f(\sigma)
=
\begin{cases}
0, & \sigma=\{i\},\\
w(e), & \sigma=e\in E.
\end{cases}
\]
The face inequality
\[
f(\tau)\leq f(\sigma)
\qquad\text{whenever }\tau\subseteq\sigma
\]
ensures that each sublevel set is a simplicial
complex
\cite{edelsbrunner-harer}.

All vertices are present at \(0\), and no edge is.
For
\[
t\geq\max_{e\in E}w(e),
\]
the threshold graph is \(G\), which is connected.

Alternatively, one may use its clique complex,
\[
\widehat{K}(t)
=
\left\{
\varnothing\neq\sigma\subseteq X:
\{i,j\}\in E(t)
\text{ for all distinct }i,j\in\sigma
\right\}.
\]
The inclusion
\[
K(t)\hookrightarrow\widehat{K}(t)
\]
induces an isomorphism on \(H_0\), naturally in
\(t\), because the two complexes have the same
vertices and edges.
Thus these two filtrations have the same
\(H_0\) persistence barcode.

\end{enumerate}

Figure~\ref{fig:rips-cech-alpha} compares the
Vietoris--Rips, \v{C}ech, and unweighted alpha
constructions on the same point cloud, illustrating
how different simplices can produce the same
connected-component structure at corresponding scales.

\begingroup

\newcommand{\sixpointvertices}{%
  \foreach \point/\pointcolour/\position in {
    P/violet/above left,
    Q/orange/above right,
    R/green!40!black/above,
    S/blue/below left,
    T/red/below right,
    U/teal/right}{
    \fill[\pointcolour]
      (\point) circle[radius=1.6pt];
    \node[
      point label,
      \position=3pt,
      text=\pointcolour
    ] at (\point) {$\point$};
  }%
}

\newcommand{\sixpointcoordinates}{%
  \coordinate (P) at (0,0);
  \coordinate (Q) at (1,0);
  \coordinate (R) at (0.5,{sqrt(3)/2});
  \coordinate (S) at (0,-1);
  \coordinate (T) at (1,-1);
  \coordinate (U) at (1.2,-0.7);
}

\newcommand{\sixpointedges}{%
  \draw[complex edge]
    (P) -- (Q) -- (R) -- cycle;
  \draw[complex edge]
    (P) -- (S) -- (T) -- (U) -- (Q);

  \ifnum\mode=3\relax
  \else
    \draw[complex edge] (Q) -- (T);
  \fi
}

\begin{figure}[htbp]
\centering

\begingroup
\setlength{\fboxsep}{6pt}
\makebox[\linewidth][c]{%
  \fcolorbox{gray!65}{white}{%
    \small
    \(
      {\color{violet}P=(0,0)},\;
      {\color{orange}Q=(1,0)},\;
      {\color{green!40!black}
        R=(\frac12,\frac{\sqrt3}{2})},\;
      {\color{blue}S=(0,-1)},\;
      {\color{red}T=(1,-1)},\;
      {\color{teal}U=(1.2,-0.7)},\;
      \rho=0.51
    \)%
  }%
}
\endgroup

\par\vspace{5mm}

\makebox[\linewidth][c]{%
\begin{tikzpicture}[
  x=1.2cm,
  y=1.2cm,
  point P boundary/.style={
    violet,
    line width=1.2pt,
    solid
  },
  point Q boundary/.style={
    orange,
    line width=1.2pt,
    dash pattern=on 4pt off 2pt
  },
  point R boundary/.style={
    green!40!black,
    line width=1.2pt,
    dash pattern=on 3pt off 1.5pt on 1pt off 1.5pt
  },
  point S boundary/.style={
    blue,
    line width=0.75pt,
    line cap=round,
    dash pattern=on 0pt off 3.5pt
  },
  point T boundary/.style={
    red,
    line width=1.2pt,
    dash pattern=on 6pt off 3pt
  },
  point U boundary/.style={
    teal,
    line width=1.2pt,
    dash pattern=on 5pt off 1.5pt on 1pt off 1.5pt
  },
  complex edge/.style={
    black!80,
    line width=1.1pt,
    line join=round
  },
  point label/.style={
    font=\small,
    fill=white,
    inner sep=1pt
  }
]

%
%
%
%

\foreach \xshift/\mode/\heading/\threshold in {
  0/1/{(a) Vietoris--Rips}/{Edge threshold \(2\rho=1.02\)},
  5/2/{(b) \v{C}ech}/{Ball radius \(\rho=0.51\)},
  10/3/{(c) Unweighted alpha}/{Restricted-ball radius \(\rho=0.51\)}
}{
  \begin{scope}[shift={(\xshift cm,0)}]

    \sixpointcoordinates

    \node[font=\small,align=center]
      at (0.6,2.50) {\heading};

    \node[font=\scriptsize,align=center]
      at (0.6,1.85) {\threshold};


    \ifnum\mode=1\relax
      \sixpointedges
    \fi

    \ifnum\mode=2\relax

      \foreach \point in {P,Q,R,S,T,U}{
        \fill[gray!8]
          (\point) circle[radius=0.51];
      }

      \foreach \point in {P,Q,R,S,T,U}{
        \draw[point \point\space boundary]
          (\point) circle[radius=0.51];
      }

    \fi

    \ifnum\mode=3\relax

      %

      \begin{scope}
        \clip
          (0.5,-0.5)
          -- (0.5,{sqrt(3)/6})
          -- (-0.7,{1.7/sqrt(3)})
          -- (-0.7,-0.5)
          -- cycle;

        \fill[
          pattern=north east lines,
          pattern color=violet!75
        ]
          (P) circle[radius=0.51];

        \draw[point P boundary]
          (P) circle[radius=0.51];
      \end{scope}

      \begin{scope}
        \clip
          (1.9,{-17/140})
          -- (1.9,{1.9/sqrt(3)})
          -- (0.5,{sqrt(3)/6})
          -- (0.5,-0.5)
          -- (0.575,-0.5)
          -- cycle;

        \fill[
          pattern=north west lines,
          pattern color=orange!85!black
        ]
          (Q) circle[radius=0.51];

        \draw[point Q boundary]
          (Q) circle[radius=0.51];
      \end{scope}

      \begin{scope}
        \clip
          (0.5,{sqrt(3)/6})
          -- (1.9,{1.9/sqrt(3)})
          -- (1.9,1.5)
          -- (-0.7,1.5)
          -- (-0.7,{1.7/sqrt(3)})
          -- cycle;

        \fill[
          pattern=dots,
          pattern color=green!40!black
        ]
          (R) circle[radius=0.51];

        \draw[point R boundary]
          (R) circle[radius=0.51];
      \end{scope}

      \begin{scope}
        \clip
          (-0.7,-1.7)
          -- (0.5,-1.7)
          -- (0.5,-0.5)
          -- (-0.7,-0.5)
          -- cycle;

        \fill[
          pattern=horizontal lines,
          pattern color=blue!70
        ]
          (S) circle[radius=0.51];

        \draw[point S boundary]
          (S) circle[radius=0.51];
      \end{scope}

      \begin{scope}
        \clip
          (0.5,-1.7)
          -- (1.9,-1.7)
          -- (1.9,{-83/60})
          -- (0.575,-0.5)
          -- (0.5,-0.5)
          -- cycle;

        \fill[
          pattern=vertical lines,
          pattern color=red!70
        ]
          (T) circle[radius=0.51];

        \draw[point T boundary]
          (T) circle[radius=0.51];
      \end{scope}

      \begin{scope}
        \clip
          (0.575,-0.5)
          -- (1.9,{-83/60})
          -- (1.9,{-17/140})
          -- cycle;

        \fill[
          pattern=grid,
          pattern color=teal!70
        ]
          (U) circle[radius=0.51];

        \draw[point U boundary]
          (U) circle[radius=0.51];
      \end{scope}


      \draw[landscape guide]
        (0.5,-0.5)
        -- (0.5,{sqrt(3)/6});

      \draw[landscape guide]
        (0.5,{sqrt(3)/6})
        -- (-0.65,{1.65/sqrt(3)});

      \draw[landscape guide]
        (0.5,{sqrt(3)/6})
        -- (1.75,{1.75/sqrt(3)});

      \draw[landscape guide]
        (-0.65,-0.5) -- (0.5,-0.5);

      \draw[landscape guide]
        (0.5,-0.5) -- (0.575,-0.5);

      \draw[landscape guide]
        (0.5,-0.5) -- (0.5,-1.65);

      \draw[landscape guide]
        (0.575,-0.5)
        -- (1.75,{-23/140});

      \draw[landscape guide]
        (0.575,-0.5)
        -- (1.75,{-77/60});

    \fi

    \sixpointvertices


    \begin{scope}[shift={(0,-4.4cm)}]

      \sixpointcoordinates

      \ifnum\mode=1\relax
        \node[font=\small]
          at (0.6,1.50)
          {Resulting complex for \(2\rho\)};
      \else
        \node[font=\small]
          at (0.6,1.50)
          {Resulting complex for \(\rho\)};
      \fi

      \ifnum\mode=1\relax
        \fill[gray!25]
          (P) -- (Q) -- (R) -- cycle;
      \fi

      \ifnum\mode=3\relax
      \else
        \fill[gray!25]
          (Q) -- (T) -- (U) -- cycle;
      \fi

      \sixpointedges
      \sixpointvertices

    \end{scope}

  \end{scope}
}

%
%

\foreach \xshift/\mode in {
  0/1,
  5/2,
  10/3
}{
  \begin{scope}[shift={(\xshift cm,-7.4cm)}]
    \begin{scope}[
      shift={(-1.125cm,0)},
      x=4.5cm,
      y=1cm
    ]

      \node[font=\small]
        at (0.41,0.65)
        {Combined barcode};

      \draw[landscape guide]
        (0.5,0.16) -- (0.5,-3.18);


      \node[
        anchor=east,
        font=\small
      ] at (-0.035,-0.875) {$H_0$};

      \draw[
        ->,
        gray!75,
        thick,
        dash pattern=on 4pt off 2pt
      ]
        (0,0) -- (0.80,0);

      \fill[gray!75]
        (0,0) circle[radius=1.2pt];

      \foreach \death/\y in {
        {sqrt(13)/20}/-0.35,
        {sqrt(53)/20}/-0.70,
        0.5/-1.05,
        0.5/-1.40,
        0.5/-1.75
      }{
        \draw[black!75,line width=1.2pt]
          (0,\y) -- ({\death},\y);

        \fill[black!75]
          (0,\y) circle[radius=1.2pt];

        \draw[
          draw=black!75,
          fill=white,
          line width=0.8pt
        ]
          ({\death},\y) circle[radius=1.4pt];
      }

      \draw[gray!45,thin]
        (0,-2.05) -- (0.82,-2.05);


      \node[
        anchor=east,
        font=\small
      ] at (-0.035,-2.725) {$H_1$};

      \ifnum\mode=1\relax

        \node[
          font=\scriptsize,
          text=gray!75,
          align=center
        ] at (0.41,-2.45)
          {No positive-length triangular bar};

      \else

        \draw[
          violet,
          line width=1.4pt,
          solid
        ]
          (0.5,-2.45)
          -- ({1/sqrt(3)},-2.45);

        \fill[violet]
          (0.5,-2.45) circle[radius=1.2pt];

        \draw[
          draw=violet,
          fill=white,
          line width=0.8pt
        ]
          ({1/sqrt(3)},-2.45)
          circle[radius=1.4pt];

        \node[
          anchor=east,
          font=\scriptsize
        ] at (0.47,-2.45)
          {$
            {\color{violet}P}
            {\color{orange}Q}
            {\color{green!40!black}R}
          $};

      \fi

      \ifnum\mode=1\relax
        \def\lowerdeath{sqrt(193)/20}
      \else
        \def\lowerdeath{1/sqrt(2)}
      \fi

      \draw[
        orange,
        line width=1.3pt,
        dash pattern=on 4pt off 2pt
      ]
        (0.5,-3.00)
        -- ({\lowerdeath},-3.00);

      \fill[orange]
        (0.5,-3.00) circle[radius=1.2pt];

      \draw[
        draw=orange,
        fill=white,
        line width=0.8pt
      ]
        ({\lowerdeath},-3.00)
        circle[radius=1.4pt];

      \ifnum\mode=3\relax
        \node[
          anchor=east,
          font=\scriptsize
        ] at (0.47,-3.00)
          {$
            {\color{violet}P}
            {\color{orange}Q}
            {\color{teal}U}
            {\color{red}T}
            {\color{blue}S}
          $};
      \else
        \node[
          anchor=east,
          font=\scriptsize
        ] at (0.47,-3.00)
          {$
            {\color{violet}P}
            {\color{orange}Q}
            {\color{red}T}
            {\color{blue}S}
          $};
      \fi


      \draw[landscape axis]
        (0,-3.50) -- (0.82,-3.50)
        node[right] {$\rho$};

      \draw[gray!80]
        (0,-3.50) -- (0,-3.56);
      \node[below=3pt,font=\scriptsize]
        at (0,-3.56) {$0$};

      \draw[gray!80]
        ({sqrt(13)/20},-3.50)
        -- ({sqrt(13)/20},-3.56);
      \node[below=3pt,font=\scriptsize]
        at ({sqrt(13)/20},-3.56)
        {$\frac{\sqrt{13}}{20}$};

      \draw[gray!80]
        ({sqrt(53)/20},-3.50)
        -- ({sqrt(53)/20},-3.56);
      \node[below=15pt,font=\scriptsize]
        at ({sqrt(53)/20},-3.56)
        {$\frac{\sqrt{53}}{20}$};

      \draw[gray!80]
        (0.5,-3.50) -- (0.5,-3.56);
      \node[below=3pt,font=\scriptsize]
        at (0.5,-3.56) {$\frac12$};

      \ifnum\mode=1\relax
      \else
        \draw[gray!80]
          ({1/sqrt(3)},-3.50)
          -- ({1/sqrt(3)},-3.56);
        \node[below=15pt,font=\scriptsize]
          at ({1/sqrt(3)},-3.56)
          {$\frac1{\sqrt3}$};
      \fi

      \draw[gray!80]
        ({\lowerdeath},-3.50)
        -- ({\lowerdeath},-3.56);

      \ifnum\mode=1\relax
        \node[below=3pt,font=\scriptsize]
          at ({\lowerdeath},-3.56)
          {$\frac{\sqrt{193}}{20}$};
      \else
        \node[below=3pt,font=\scriptsize]
          at ({\lowerdeath},-3.56)
          {$\frac1{\sqrt2}$};
      \fi

    \end{scope}
  \end{scope}
}

\end{tikzpicture}%
}

\caption{
Vietoris--Rips, \v{C}ech, and unweighted alpha constructions
on the same six-point cloud, at edge threshold
\(2\rho=1.02\) for Vietoris--Rips and radius
\(\rho=0.51\) for \v{C}ech and alpha.
Vertex colours identify the corresponding ball boundaries
and restricted Voronoi cells:
\(
  {\color{violet}P},
  {\color{orange}Q},
  {\color{green!40!black}R},
  {\color{blue}S},
  {\color{red}T},
  {\color{teal}U}
\).
Boundary styles and cell patterns provide additional distinctions;
grey dashed lines show Voronoi boundaries.
Vietoris--Rips fills the triangles
\(
  {\color{violet}P}
  {\color{orange}Q}
  {\color{green!40!black}R}
\)
and
\(
  {\color{orange}Q}
  {\color{red}T}
  {\color{teal}U}
\);
\v{C}ech fills only the latter, while alpha contains neither
and omits the edge
\({\color{orange}Q}{\color{red}T}\).
All three complexes are connected, with one independent
one-dimensional hole for Vietoris--Rips and two for
\v{C}ech and alpha.
The barcodes share the radius parameter \(\rho\),
with Vietoris--Rips interpreted as \(\mathrm{VR}_{2\rho}\).
Each \(H_0\) barcode has one infinite bar and finite bars
\([0,\sqrt{13}/20)\), \([0,\sqrt{53}/20)\), and
three copies of \([0,1/2)\).
The triangular \(H_1\) bar \([1/2,1/\sqrt3)\) occurs
only for \v{C}ech and alpha.
The remaining \(H_1\) bar is
\([1/2,1/\sqrt2)\) for \v{C}ech and alpha. However, it is
\([1/2,\sqrt{193}/20)\) for Vietoris--Rips, because when \(2\rho=\sqrt{193}/10\), the triangles
\({\color{violet}P}{\color{orange}Q}{\color{teal}U}\),
\({\color{orange}Q}{\color{red}T}{\color{teal}U}\),
\({\color{blue}S}{\color{red}T}{\color{teal}U}\), and
\({\color{violet}P}{\color{blue}S}{\color{teal}U}\)
are all present in the Vietoris--Rips complex and together
fill the hole bounded by the cycle
\({\color{violet}P}{\color{orange}Q}{\color{red}T}{\color{blue}S}\).
Cycle labels identify representatives at the displayed parameter.
Filled and open endpoints indicate included births and excluded
deaths; all barcode plots use the same horizontal scale.
}
\label{fig:rips-cech-alpha}
\end{figure}
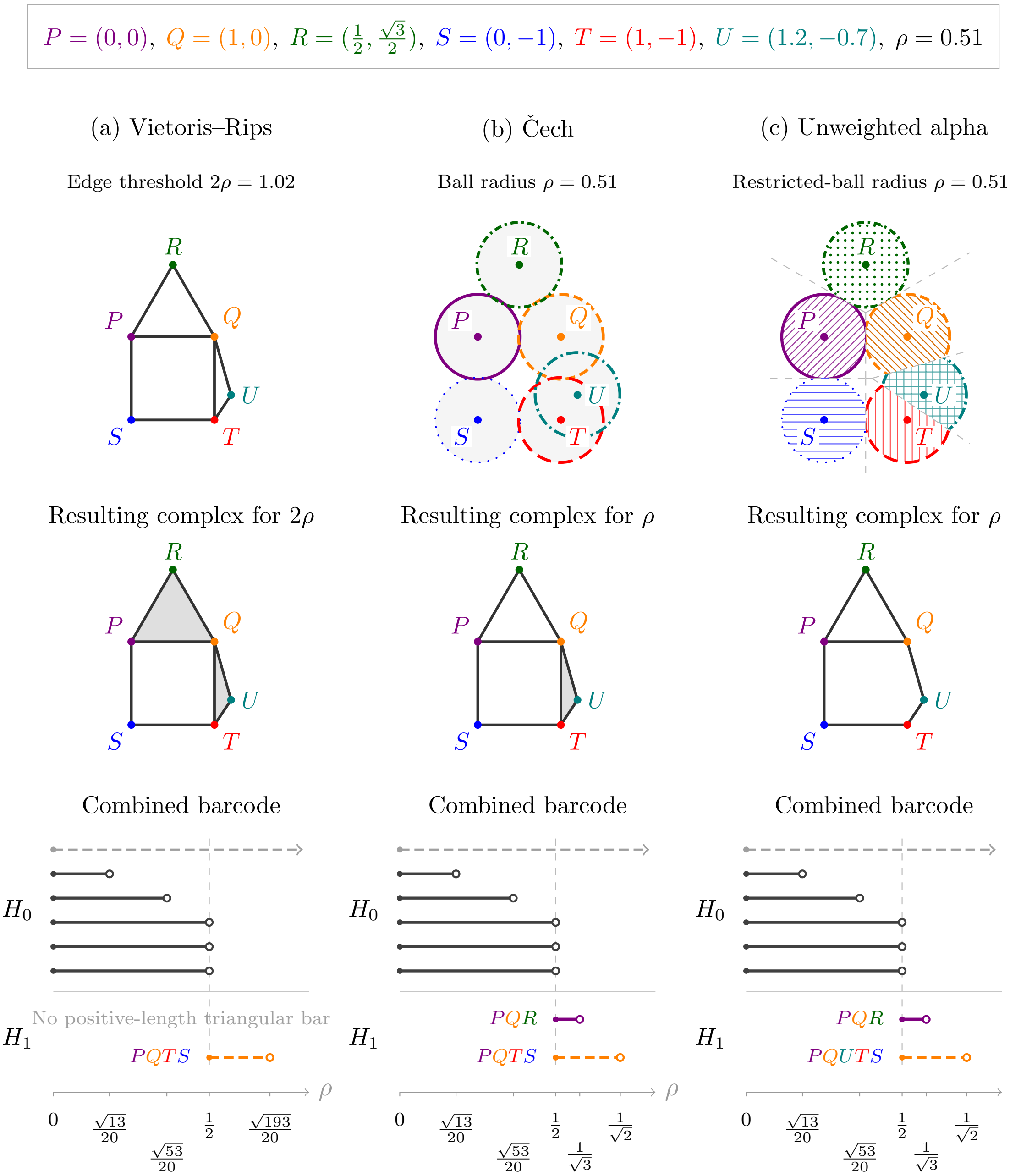

\endgroup

All four constructions give increasing,
right-continuous filtrations that start with the
\(L\) isolated vertices and eventually become connected.
Right-continuity means that a simplex is included at
the threshold where it first appears.
For the nerve constructions, this follows from the
attainment of the relevant minimum-radius values.
Consequently, their ordinary \(H_0\) barcodes consist
of one infinite bar \([0,\infty)\) and \(L-1\)
positive-length finite bars, all born at \(0\).

Allowing the input data to depend on \(g\) gives
families \(\mathcal K_g\) satisfying these common-birth
hypotheses at every parameter value for which the
input assumptions hold.
No uniform bound in \(g\) on the filtration value
at which connectivity is reached is required.
To apply Theorem~\ref{prop:closed} at a particular
\(g_0\), the sorted death times must additionally
be differentiable there.

The edge-threshold graph construction includes the
filtered \(1\)-skeleton of a Vietoris--Rips filtration
as the case of a complete weighted graph.
Passing to its clique complex gives the full
Vietoris--Rips filtration without changing \(H_0\).
More generally, \(H_0\) persistence depends only on
the filtered vertices and edges: adding
higher-dimensional simplices does not change the
connected components or their merger times.
Thus those higher-dimensional simplices do not enter
the closed formula.

The numerical death times, and hence the landscape
velocity, depend on the convention used for the
filtration parameter \(t\).
For distinct Euclidean points, an edge appears in
the Vietoris--Rips filtration at its Euclidean length,
whereas it appears in the \v{C}ech filtration at half
that length when \(t\) denotes the ball radius.
Their \(H_0\) death times therefore differ by a factor
of \(2\).
Likewise, parameterizing an alpha filtration by
squared radius rather than radius squares its death
times.
Theorem~\ref{prop:closed} applies in each convention,
provided the death times and their derivatives with
respect to \(g\) are computed consistently using
that convention.

The stated assumptions are essential to this setting.
Nonzero or unequal vertex birth times, or truncation
before the filtration becomes connected, can prevent
the barcode from having the prescribed form.
Coincident points or zero dissimilarities between
distinct vertices can likewise produce connections
already at \(t=0\), reducing the number of
positive-length finite bars below \(L-1\).
Such cases are not covered by the present
\(L-1\)-bar formulation without modifying its
assumptions.

\subsection{A closed formula: computing landscape velocity from death times}

This subsection develops the key technical result: a closed formula for landscape velocity that depends only on the death times of the finite bars, without requiring explicit computation of the landscape itself. This simplification is crucial for both theoretical understanding and practical computation.

\begin{lemma}[Landscape layers for common-birth barcodes]
\label{lem:nested}
Under ~Assumption~\ref{ass:common-birth}, the tents
associated with the finite \(H_0\) bars are pointwise
ordered:
\[
\Lambda_{(0,d_1(g))}(t)
\geq
\Lambda_{(0,d_2(g))}(t)
\geq\cdots\geq
\Lambda_{(0,d_{L-1}(g))}(t)
\]
for every \(t\in\mathbb{R}\).
Consequently,
\[
\lambda_k(g;t)
=
\Lambda_{(0,d_k(g))}(t)
=
\max\{0,\min\{t,d_k(g)-t\}\},
\qquad
k=1,\dots,L-1.
\]
\end{lemma}

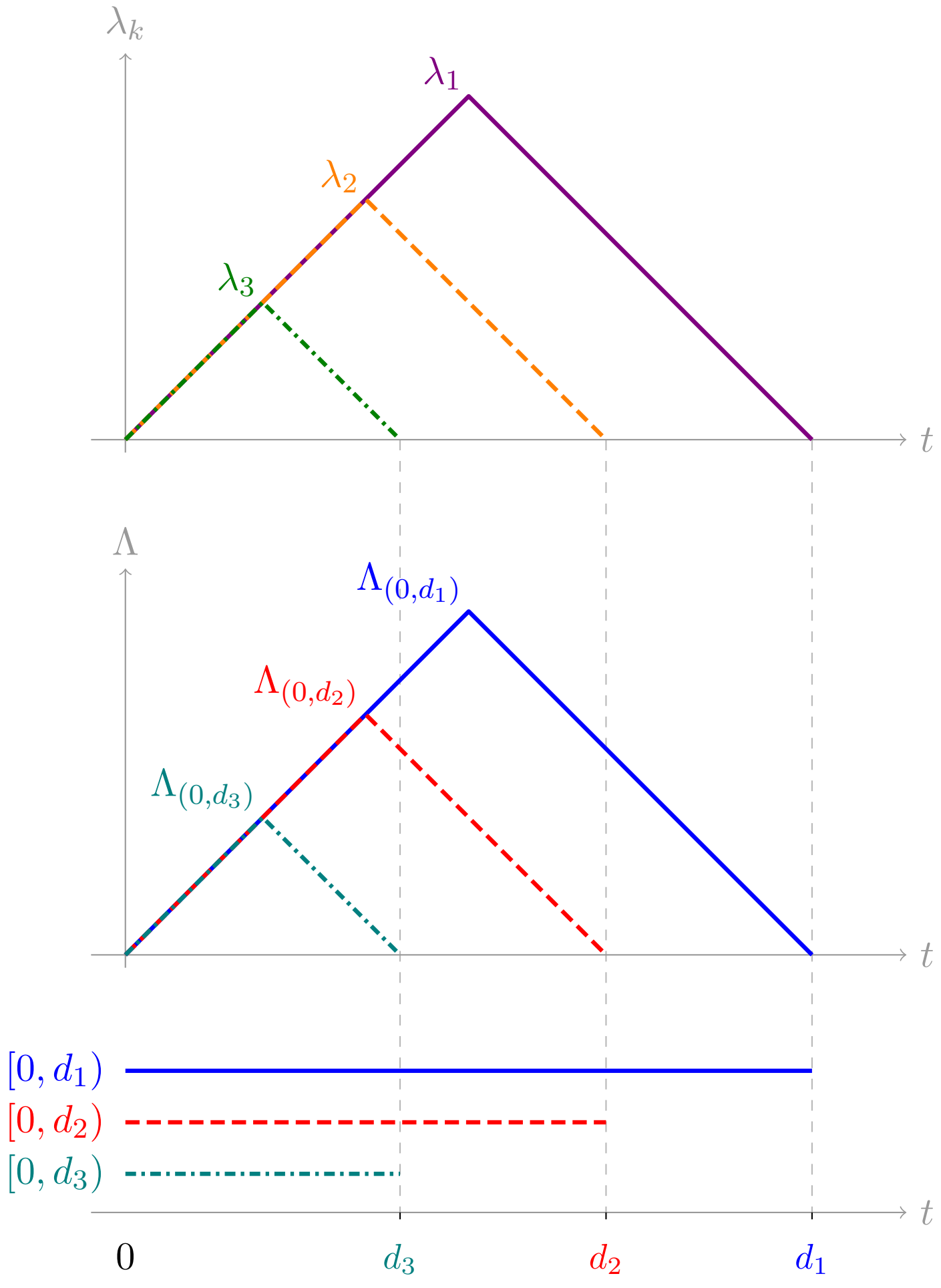
\begin{figure}[htbp]
\centering
\begin{tikzpicture}[x=0.85cm,y=0.85cm]

  \tikzset{
    texture blue/.style={solid},
    texture red/.style={dash pattern=on 4pt off 2pt},
    texture teal/.style={dash pattern=on 3pt off 1.5pt on 1pt off 1.5pt},
    texture violet/.style={solid},
    texture orange/.style={dash pattern=on 4pt off 2pt},
    texture green/.style={dash pattern=on 3pt off 1.5pt on 1pt off 1.5pt},
  }

  %

  \foreach \x in {3.2,5.6,8}{
    \draw[landscape guide] (\x,6) -- (\x,-3);
  }

  \begin{scope}[shift={(0,6)}]

    \draw[landscape axis]
      (-0.4,0) -- (9.1,0) node[right] {$t$};

    \draw[landscape axis]
      (0,-0.15) -- (0,4.5) node[above] {$\lambda_k$};

    \draw[landscape curve,violet,texture violet]
      (0,0) -- (4,4) -- (8,0);

    \draw[landscape curve,orange,texture orange]
      (0,0) -- (2.8,2.8) -- (5.6,0);

    \draw[landscape curve,green!50!black,texture green]
      (0,0) -- (1.6,1.6) -- (3.2,0);

    \node[violet,above left=4pt] at (4.2,3.8)
      {$\lambda_1$};

    \node[orange,above left=4pt] at (3,2.6)
      {$\lambda_2$};

    \node[
      text=green!50!black,
      above left=4pt
    ] at (1.8,1.4)
      {$\lambda_3$};

  \end{scope}

  \begin{scope}

    \draw[landscape axis]
      (-0.4,0) -- (9.1,0) node[right] {$t$};

    \draw[landscape axis]
      (0,-0.15) -- (0,4.5) node[above] {$\Lambda$};

    \draw[landscape curve,blue,texture blue]
      (0,0) -- (4,4) -- (8,0);

    \draw[landscape curve,red,texture red]
      (0,0) -- (2.8,2.8) -- (5.6,0);

    \draw[landscape curve,teal,texture teal]
      (0,0) -- (1.6,1.6) -- (3.2,0);

    \node[blue,above left=4pt] at (4.2,3.8)
      {$\Lambda_{(0,d_1)}$};

    \node[red,above left=4pt] at (3,2.6)
      {$\Lambda_{(0,d_2)}$};

    \node[teal,above left=4pt] at (1.8,1.4)
      {$\Lambda_{(0,d_3)}$};

  \end{scope}

  \begin{scope}[shift={(0,-3)}]

    \draw[landscape axis]
      (-0.4,0) -- (9.1,0) node[right] {$t$};

    \draw[landscape curve,blue,texture blue]
      (0,1.65) -- (8,1.65);

    \draw[landscape curve,red,texture red]
      (0,1.05) -- (5.6,1.05);

    \draw[landscape curve,teal,texture teal]
      (0,0.45) -- (3.2,0.45);

    \node[
      anchor=east,
      text=blue,
      inner sep=2pt
    ] at (-0.15,1.65)
      {$[0,d_1)$};

    \node[
      anchor=east,
      text=red,
      inner sep=2pt
    ] at (-0.15,1.05)
      {$[0,d_2)$};

    \node[
      anchor=east,
      text=teal,
      inner sep=2pt
    ] at (-0.15,0.45)
      {$[0,d_3)$};

    \node[below=5pt] at (0,0) {$0$};

    \foreach \x/\lab/\col in {
      3.2/{d_3}/teal,
      5.6/{d_2}/red,
      8/{d_1}/blue}{
      \draw (\x,0) -- (\x,-0.08);
      \node[
        below=5pt,
        font=\small,
        text=\col
      ] at (\x,0) {$\lab$};
    }

  \end{scope}

\end{tikzpicture}
\caption{
Illustration of Lemma~\ref{lem:nested} for the \(H_0\)
barcode of a Vietoris--Rips, \v{C}ech or alpha complex  filtration of distinct Euclidean points, or another
filtration satisfying ~Assumption~\ref{ass:common-birth}.
Only the finite part of the barcode is shown;
the single infinite bar \([0,\infty)\) is omitted
and does not contribute to the landscape considered here.
The 3 finite bars (bottom) determine individual tent
functions (middle), whose pointwise ordered values
form the landscape layers (top).
Since all bars are born at \(0\) and
\({\color{blue}d_1}\geq {\color{red}d_2}\geq {\color{teal}d_3}\), the tents are nested.
Consequently, \(\lambda_k(t)=\Lambda_{(0,d_k)}(t)\):
each layer coincides with a single tent, unlike
the distinct-birth example.
The common rising sides coincide; the colors distinguish
the tents and layers along their falling sides.
\textcolor{violet}{Violet}, \textcolor{orange}{orange}, and \textcolor{green!50!black}{dark green} identify the landscape
layers, while \textcolor{blue}{blue}, \textcolor{red}{red}, and \textcolor{teal}{teal} identify the original
bars, tents, and corresponding axis labels.
Each data series also has a distinct line texture:
the first interval and first landscape layer are solid,
the second are dashed, and the third are dash-dotted.
All panels use the same horizontal scale, and the
continuous dashed guides align their death times.
}
\label{fig:common-birth-barcode-to-landscape}
\end{figure}

\noindent For an illustration of Lemma \ref{lem:nested}, see Figure~\ref{fig:common-birth-barcode-to-landscape}.

\begin{proof}
Let \(0\leq d'\leq d\).
For every \(t\in\mathbb{R}\), we have
\[
d'-t\leq d-t.
\]
Since taking the minimum with \(t\), and then taking
the maximum with \(0\), both preserve inequalities,
\[
\max\{0,\min\{t,d'-t\}\}
\leq
\max\{0,\min\{t,d-t\}\}.
\]
Thus
\[
\Lambda_{(0,d')}(t)\leq\Lambda_{(0,d)}(t)
\]
for every \(t\).
Applying this to the sorted death times shows that
the tents are already ordered pointwise.
Their \(k\)-th largest value is therefore
\(\Lambda_{(0,d_k(g))}(t)\).
\end{proof}

\begin{theorem}[Landscape velocity from death times
and their derivatives]
\label{prop:closed}
Under Assumption~\ref{ass:common-birth}, let \(g_0\in I\)
be a point at which every sorted death time \(d_k(g)\)
is differentiable.
Then
\begin{equation}\label{eq:closed}
\boxed{
\xi(g_0)
=
\left(
\frac12
\sum_{k=1}^{L-1}
d_k(g_0)
\bigl(\partial_g d_k(g_0)\bigr)^2
\right)^{1/2}.
}
\end{equation}
\end{theorem}

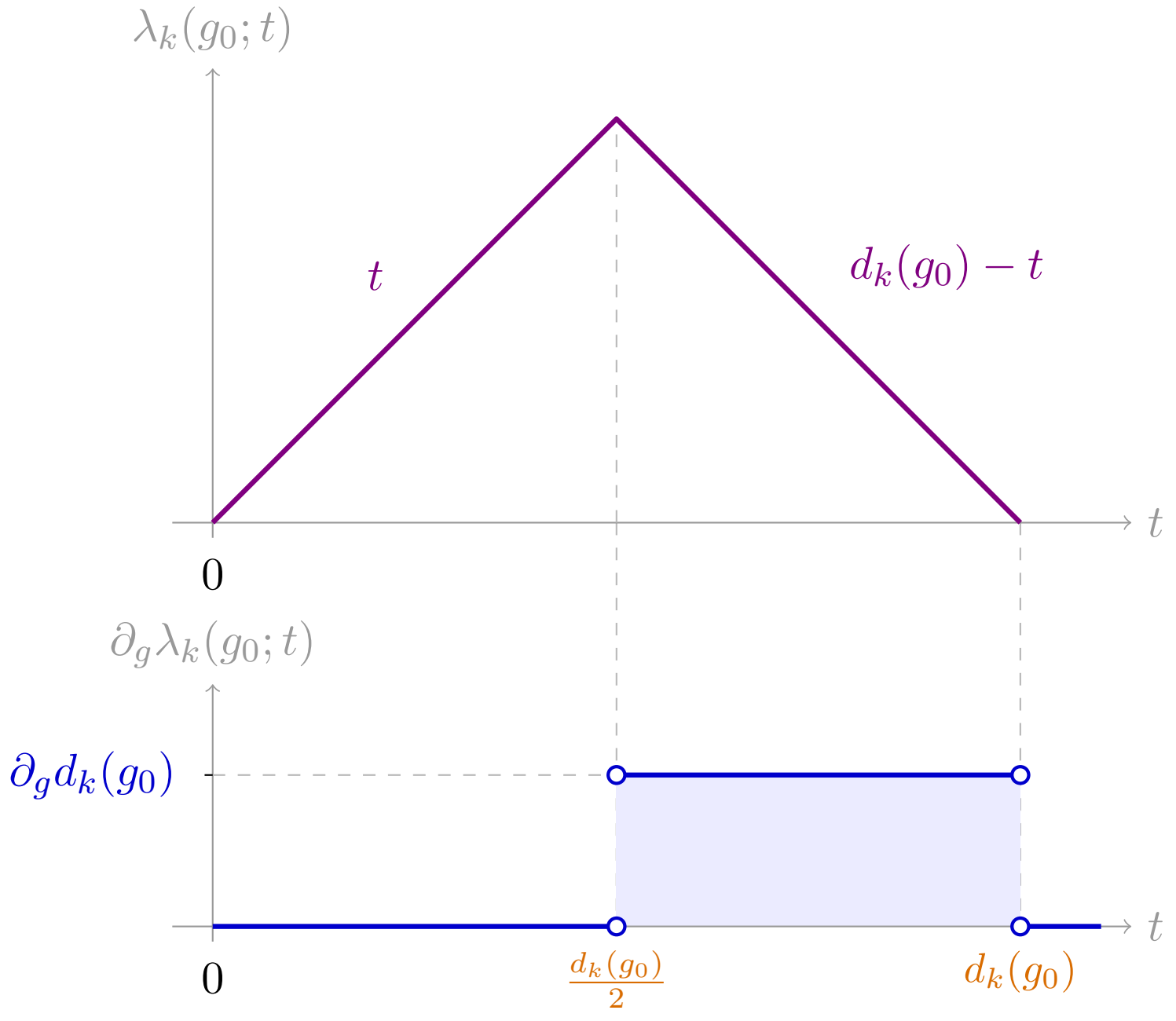
\begin{figure}[!htb]
\centering
\begin{tikzpicture}[x=0.85cm,y=0.85cm]


  \foreach \x in {4,8}{
    \draw[landscape guide] (\x,0) -- (\x,-4);
  }

  \begin{scope}

    \draw[landscape axis]
      (-0.4,0) -- (9.1,0)
      node[right,text=gray!80] {$t$};

    \draw[landscape axis]
      (0,-0.15) -- (0,4.5)
      node[above] {
        $\lambda_k(g_0;{\color{gray!80}t})$
      };

    \draw[landscape guide] (4,0) -- (4,4);

    \draw[landscape curve,violet]
      (0,0) -- (4,4) -- (8,0);

    \path (0,0) -- (4,4)
      node[midway,above left=5pt,text=violet] {$t$};

    \path (4,4) -- (8,0)
      node[midway,above right=5pt,text=violet] {$d_k(g_0)-t$};

    \node[below=5pt] at (0,0) {$0$};


  \end{scope}

  \begin{scope}[shift={(0,-4)}]

    \fill[blue!8] (4,0) rectangle (8,1.5);

    \draw[landscape axis]
      (-0.4,0) -- (9.1,0)
      node[right,text=gray!80] {$t$};

    \draw[landscape axis]
      (0,-0.15) -- (0,2.4)
      node[above] {
        $\partial_g\lambda_k(g_0;{\color{gray!80}t})$
      };

    \draw[landscape guide] (0,1.5) -- (4,1.5);

    \draw (0,1.5) -- (-0.08,1.5);
    \node[left=5pt,text=blue!80!black]
      at (0,1.5) {$\partial_g d_k(g_0)$};

    \draw[landscape curve,blue!80!black]
      (0,0) -- (4,0);
    \draw[landscape curve,blue!80!black]
      (4,1.5) -- (8,1.5);
    \draw[landscape curve,blue!80!black]
      (8,0) -- (8.8,0);

    \foreach \x/\y in {4/0,4/1.5,8/0,8/1.5}{
      \draw[
        draw=blue!80!black,
        fill=white,
        thick
      ] (\x,\y) circle[radius=2pt];
    }

    \node[below=5pt] at (0,0) {$0$};

    \node[
      anchor=north,
      font=\small,
      text=orange!85!black,
      inner sep=2pt
    ] at (4,-0.15) {$\frac{d_k(g_0)}{2}$};

    \node[
      anchor=north,
      font=\small,
      text=orange!85!black,
      inner sep=2pt
    ] at (8,-0.15) {$d_k(g_0)$};

  \end{scope}

\end{tikzpicture}
\caption{
Illustration of equation~\eqref{eq:landscape-derivative}:
a landscape layer (above, \textcolor{violet}{violet})
and its derivative with respect to \(g\)
(below, \textcolor{blue!80!black}{blue}).
The rising branch equals \({\color{violet}t}\),
independent of \(g\), so its derivative is zero.
The falling branch equals
\({\color{violet}d_k(g_0)-t}\),
so its parameter derivative is
\({\color{blue!80!black}\partial_g d_k(g_0)}\)
on
\(\displaystyle
  \left(
    {\color{orange!85!black}\frac{d_k(g_0)}{2}},
    {\color{orange!85!black}d_k(g_0)}
  \right)
\).
Outside this interval, the parameter derivative vanishes.
At the two exceptional breakpoint values
\(
  {\color{gray!80}t}
  ={\color{orange!85!black}\frac{d_k(g_0)}{2}}
\)
and
\(
  {\color{gray!80}t}
  ={\color{orange!85!black}d_k(g_0)}
\),
no value of the derivative is asserted.
The lower panel illustrates
\({\color{blue!80!black}\partial_g d_k(g_0)>0}\);
a negative derivative gives a rectangle below the axis.
The continuous dashed guides align the two panels.
}
\label{fig:landscape-derivative}
\end{figure}

\noindent For an illustration of the proof of Theorem \ref{prop:closed}, see Figure \ref{fig:landscape-derivative}.

\begin{proof}
Fix \(k\in\{1,\dots,L-1\}\), and write
\[
d=d_k(g_0),
\qquad
d_h=d_k(g_0+h).
\]
By Lemma~\ref{lem:nested},
\[
\lambda_k(g;t)
=
\max\{0,\min\{t,d_k(g)-t\}\}.
\]
Differentiability of \(d_k\) at \(g_0\) implies
continuity there, so \(d_h\to d\) as \(h\to0\).

Fix \(t>0\) with \(t\notin\{d/2,d\}\).
Because \(t\) is strictly separated from these two
limiting breakpoints, its position relative to
\(d_h/2\) and \(d_h\) agrees, for all sufficiently
small \(h\), with its position relative to
\(d/2\) and \(d\).
We may therefore compute the difference quotient
within a single branch of the tent.

If \(0<t<d/2\), then \(t<d_h/2\) for all sufficiently
small \(h\). Hence
\[
\lambda_k(g_0+h;t)=\lambda_k(g_0;t)=t,
\]
and
\[
\partial_g\lambda_k(g_0;t)=0.
\]

If \(d/2<t<d\), then \(d_h/2<t<d_h\) for all
sufficiently small \(h\). Hence
\[
\lambda_k(g_0+h;t)=d_h-t,
\qquad
\lambda_k(g_0;t)=d-t,
\]
and
\[
\partial_g\lambda_k(g_0;t)
=
\lim_{h\to0}
\frac{d_h-d}{h}
=
\partial_g d_k(g_0).
\]

If \(t>d\), then \(t>d_h\) for all sufficiently
small \(h\). Both tent values vanish, so
\[
\partial_g\lambda_k(g_0;t)=0.
\]

At \(t=d/2\) and \(t=d\), the derivative need not
exist. These two exceptional values form a
Lebesgue-null set.
Consequently, for almost every \(t>0\),
\begin{equation}\label{eq:landscape-derivative}
\partial_g\lambda_k(g_0;t)
=
\partial_g d_k(g_0)\,
\mathbf 1_{\left\{
\frac{d_k(g_0)}{2}<t<d_k(g_0)
\right\}}.
\end{equation}
It follows that
\[
\begin{aligned}
\int_0^\infty
\bigl(\partial_g\lambda_k(g_0;t)\bigr)^2\,dt
&=
\int_{d_k(g_0)/2}^{d_k(g_0)}
\bigl(\partial_g d_k(g_0)\bigr)^2\,dt\\
&=
\frac{d_k(g_0)}{2}
\bigl(\partial_g d_k(g_0)\bigr)^2,
\end{aligned}
\]
where the last equality follows because the
integrand does not depend on \(t\).
Summing over \(k\) and taking
the nonnegative square root proves the formula.
\end{proof}

\begin{corollary}[Constant-birth extension]
If every bar is $[b,d_k(g))$ with $b$ constant in $g$, then translating the tents gives
$$
\xi^2=\frac12\sum_k(d_k-b)(d_k')^2.
$$
\end{corollary}

The result follows from Theorem~\ref{prop:closed} by translating
the filtration coordinate by \(-b\), which preserves the \(L^2\)
landscape velocity and replaces each \(d_k\) by \(d_k-b\),
while leaving \(d_k'\) unchanged, since \(b\) is independent of \(g\).

\begin{remark}[Interpretation and scope]
\label{rem:closed-form}
\begin{enumerate}[label=(\roman*)]

\item
The closed formula depends only on the sorted death
times and their derivatives.
No additional information about the simplices or
the combinatorial structure of the filtration enters
the formula. We do not even need to compute the landscape, as we can extract the death times directly from the barcode.

If sorted death times fail to be differentiable at
a crossing, Theorem~\ref{prop:closed} does not
apply there as stated.
If the relevant one-sided derivatives exist, the same
proof gives the corresponding one-sided formula.
The sum is invariant under relabeling, but this
symmetry alone does not imply that a two-sided
derivative exists.

\item
Each bar contributes
\[
\frac{d_k(g_0)}{2}
\bigl(\partial_g d_k(g_0)\bigr)^2
\]
to $\xi^2(g_0)$.
The factor \(d_k(g_0)/2\) is the length of the landscape's falling
side's projection onto the \(t\)-axis.
Thus, for equal death-time derivatives, longer bars
make larger contributions, linearly in their
death times.

\item
The proof uses the common-birth property, rather
than the homological degree or the particular
filtration construction.
More precisely, it applies to any family of finite
barcodes
\[
\mathcal{B}_g^{\mathrm{fin}}
=
\{[0,d_k(g))\}_{k=1}^{M}
\]
with a fixed positive number \(M\) of bars (~Assumption~\ref{ass:common-birth-family}), at
parameter values where all sorted death times
are differentiable.
In that setting, the formula reads
\[
\xi^2(g_0)
=
\frac12
\sum_{k=1}^{M}
d_k(g_0)
\bigl(\partial_g d_k(g_0)\bigr)^2.
\]
The \(H_0\) filtrations considered here provide
a natural class of examples, with \(M=L-1\).

For barcodes with unequal or moving birth times,
the tents need not be nested, and a landscape
layer can switch between different tents.
The argument above therefore does not, in general,
give the same closed formula.

\end{enumerate}
\end{remark}

When does the differentiability assumption in Theorem~\ref{prop:closed} hold?
The following Proposition \ref{prop:regularity} uses Lemma~\ref{lem:rips-stability} to show that for Vietoris--Rips filtrations (and more generally, any family of dissimilarities that depend continuously on the parameter), the sorted death times are differentiable almost everywhere.

\begin{lemma}[Stability of sorted death times]
\label{lem:rips-stability}
For Vietoris--Rips filtrations, let \(D\) and \(\widetilde D\) be two dissimilarity
matrices satisfying the assumptions of
Subsection~\ref{subsec:rips-mst}, on the same vertex set.
Let
\[
d_1\geq\cdots\geq d_{L-1},
\qquad
\widetilde d_1\geq\cdots\geq\widetilde d_{L-1}
\]
be their sorted finite \(H_0\) death times.
Then
\begin{equation}
\label{eq:stability_d_k_vs_D}
    |d_k-\widetilde d_k|
    \leq
    \|D-\widetilde D\|_\infty,
    \qquad k=1,\dots,L-1,  
\end{equation}
where
\[
\|D-\widetilde D\|_\infty
=
\max_{i,j\in X}|D(i,j)-\widetilde D(i,j)|.
\]
\end{lemma}

\begin{proof}
Write
\[
\delta=\|D-\widetilde D\|_\infty.
\]
Every edge present at threshold \(t\) for \(D\)
is present at threshold \(t+\delta\) for
\(\widetilde D\).
The reverse statement holds after interchanging
the matrices.
Hence
\[
\mathrm{VR}_t(X,D)
\subseteq
\mathrm{VR}_{t+\delta}(X,\widetilde D),
\qquad
\mathrm{VR}_t(X,\widetilde D)
\subseteq
\mathrm{VR}_{t+\delta}(X,D).
\]


Let \(c_D(t)\) denote the number of connected
components of \(\mathrm{VR}_t(X,D)\).
Since the two filtrations have the same vertex set,
the inclusions above imply
\begin{equation}
\label{eq:interleaving_components}
c_{\widetilde D}(t+\delta)\leq c_D(t).    
\end{equation}

All \(H_0\) bars are born at zero, and a bar
with death time \(d_j\) is alive at threshold
\(t\) precisely when \(t<d_j\). Therefore
\[
c_D(t)
=
1+\#\{j\in\{1,\dots,L-1\}:d_j>t\},
\qquad t\geq0,
\]
where the additional component corresponds
to the infinite bar.
Since the finite death times are sorted
in nonincreasing order, this gives
\begin{equation}
\label{eq:components_at_threshold_d_k}
    c_D(t)\leq k
    \quad\Longleftrightarrow\quad
    d_k\leq t,
    \qquad k=1,\dots,L-1.    
\end{equation}
The same equivalence holds for \(\widetilde D\).

Taking \(t=d_k\), we obtain
\[
c_{\widetilde D}(d_k+\delta)
\overset{\eqref{eq:interleaving_components}}{\leq}
c_D(d_k)
\overset{\eqref{eq:components_at_threshold_d_k}}{\leq}
k.
\]
Applying \eqref{eq:components_at_threshold_d_k} for \(\widetilde D\)
at the inequality just proved, $c_{\widetilde D}(d_k+\delta) \leq k$,
yields
\[
\widetilde d_k\leq d_k+\delta.
\]
Interchanging \(D\) and \(\widetilde D\) gives
\(d_k\leq\widetilde d_k+\delta\).
Hence
\[
|d_k-\widetilde d_k|\leq\delta,
\]
as claimed.
\end{proof}

\begin{remark}[Relation to known stability results]
\label{rem:stability-literature}
The inclusions in the first step of the proof say that the
Vietoris--Rips filtrations of \(D\) and \(\widetilde D\) are
\(\delta\)-interleaved.
This is the argument behind the stability of Vietoris--Rips
persistence \cite{chazal-desilva-oudot}, where the algebraic
stability theorem \cite{chazal-proximity,chazal-structure}
bounds the bottleneck distance between the barcodes, in every
homological degree, by \(\delta\).
Our proof is similar, but it replaces the algebraic stability
theorem by an elementary count of connected components, which is
specific to \(H_0\).
It gives the sorted bound in Equation~\eqref{eq:stability_d_k_vs_D},
which compares the \(k\)-th deaths without matching bars to the
diagonal, and therefore implies the bottleneck bound for the
\(H_0\) barcodes.
For the corresponding stability of single-linkage clustering, see
\cite{carlsson-memoli}.
\end{remark}

\begin{proposition}[Almost-everywhere differentiability of
Vietoris--Rips death times]
\label{prop:regularity}
Fix \(X=\{1,\dots,L\}\), with \(L\geq2\), and an open interval
\(I\subseteq\mathbb R\). For each \(g\in I\), let \(D_g\) be a
finite-valued symmetric dissimilarity on \(X\), with zero diagonal
and strictly positive off-diagonal entries, as in
Subsection~\ref{subsec:rips-mst}.
Assume that every entry \(D_g(i,j)\) is locally Lipschitz
as a function of \(g\).

Then, for each \(k\), the \(k\)-th sorted death
\(d_k(g)\) is locally Lipschitz on \(I\).
Consequently, \(d_k\) is differentiable at
almost every \(g\in I\).
\end{proposition}

\begin{proof}
Fix \(g_*\in I\). Since the matrix has finitely many entries,
there is a common neighbourhood \(U\subseteq I\) of \(g_*\)
and a constant \(C_{\mathrm{Lip}}>0\) such that
\[
|D_g(i,j)-D_{\widetilde g}(i,j)|
\leq C_{\mathrm{Lip}}|g-\widetilde g|
\qquad
\text{for all }
g,\widetilde g\in U
\]
for every \(i,j\). Hence
\[
\|D_g-D_{\widetilde g}\|_\infty
\leq C_{\mathrm{Lip}}|g-\widetilde g|.
\]

By Lemma~\ref{lem:rips-stability},
\[
|d_k(g)-d_k(\widetilde g)|
\leq
\|D_g-D_{\widetilde g}\|_\infty
\leq
C_{\mathrm{Lip}}|g-\widetilde g|.
\]
Thus \(d_k\) is locally Lipschitz on \(I\).

By Rademacher's theorem \cite{rademacher}, \cite[Section~3.1]{evans-gariepy}, every locally Lipschitz
real-valued function on an open subset of
\(\mathbb R\) is differentiable almost everywhere.
Therefore each \(d_k\) is differentiable at
almost every \(g\in I\).
\end{proof}

Thus, under the local Lipschitz assumption on the
dissimilarity entries, the closed-form formula in
Theorem~\ref{prop:closed} applies at
almost every \(g_0\in I\), simultaneously for all deaths,
since there are only finitely many of them.
The exceptional parameters are those at which at least
one sorted death fails to be differentiable.


\subsection{Landscape velocity from minimum-spanning-tree edge weights}
\label{subsec:rips-mst}

The closed formula Theorem~\ref{prop:closed} reduces the computation of landscape
velocity to the finite death times and their parameter
derivatives.
For Vietoris--Rips filtrations and homology degree 0, these deaths can be obtained
directly from the input dissimilarities through the classical
correspondence between zero-dimensional persistence and
minimum spanning trees (MST); see
\cite[Sections~2.3 and 5.1]{cultrera2024chromatic} and
\cite[Lemma~20]{beers2024fiberspointcloudpersistence}.

Let \(X=\{1,\dots,L\}\), and let \(D_g\) be a finite-valued
symmetric dissimilarity matrix satisfying
\[
D_g(i,i)=0,
\qquad
D_g(i,j)>0
\quad\text{for }i\neq j.
\]
No triangle inequality is required.
We use the diameter-threshold convention
\[
\mathrm{VR}_t(X,D_g)
=
\left\{
\varnothing\neq\sigma\subseteq X:
D_g(i,j)\leq t
\text{ for all }i,j\in\sigma
\right\}.
\]
Consider the complete graph on \(X\), with edge
\(\{i,j\}\) assigned weight \(D_g(i,j)\).
A minimum spanning tree is a connected, cycle-free subgraph
containing all vertices whose total edge weight is minimal.
With the convention above, the finite \(H_0\) death times
are precisely the edge weights of any minimum spanning tree,
counted once for each tree edge.

Indeed, Kruskal's algorithm processes edges in nondecreasing
order of weight and accepts an edge exactly when it joins
two previously distinct components \cite{kruskal}.
Each accepted edge therefore records one component merger
and one finite \(H_0\) death.
Ties may be processed in any order: the resulting tree need
not be unique, but its multiset of edge weights is unchanged.
These death times are also the merge values of single-linkage
hierarchical clustering, which can be recovered from a
minimum spanning tree \cite{gower-ross};
see \cite{carlsson-memoli} for a characterization and
stability results.
Figure~\ref{fig:mst-barcode-landscape} illustrates the
correspondence between accepted edges, component mergers,
finite deaths, and landscape layers.

\begin{figure}[htbp]
\centering
\makebox[\linewidth][c]{%
\begin{tikzpicture}[
  x=1cm,
  y=1cm,
  merger first/.style={
    violet,
    line width=1.6pt,
    solid,
    line join=round
  },
  merger second/.style={
    orange,
    line width=1.3pt,
    dash pattern=on 4pt off 2pt,
    line join=round
  },
  merger third/.style={
    green!40!black,
    line width=1.3pt,
    dash pattern=on 3pt off 1.5pt on 1pt off 1.5pt,
    line join=round
  },
  graph vertex/.style={
    circle,
    draw=gray!80,
    fill=white,
    minimum size=5mm,
    inner sep=0pt,
    font=\small
  },
  edge weight/.style={
    font=\small,
    fill=white,
    inner sep=2pt
  }
]

  %
  %
  \foreach \xshift/\mode/\heading in {
    0/0/{(a) Weighted complete graph},
    7/1/{(b) Minimum spanning tree}
  }{

    \begin{scope}[shift={(\xshift cm,0)}]

      \node[font=\small,align=center]
        at (1.4,2.75) {\heading};

      \coordinate (vP) at (0,1.8);
      \coordinate (vQ) at (2.8,1.8);
      \coordinate (vR) at (2.8,0);
      \coordinate (vS) at (0,0);

      \ifnum\mode=0

        \draw[gray!75,thick]
          (vP) -- (vQ)
          node[midway,above=3pt,edge weight] {$10$};

        \draw[gray!75,thick]
          (vQ) -- (vR)
          node[midway,right=3pt,edge weight] {$20$};

        \draw[gray!75,thick]
          (vR) -- (vS)
          node[midway,below=3pt,edge weight] {$30$};

        \draw[gray!75,thick]
          (vP) -- (vR)
          node[pos=0.27,edge weight] {$40$};

        \draw[gray!75,thick]
          (vQ) -- (vS)
          node[pos=0.27,edge weight] {$50$};

        \draw[gray!75,thick]
          (vP) -- (vS)
          node[midway,left=3pt,edge weight] {$60$};

      \else

        \draw[merger third]
          (vP) -- (vQ)
          node[
            midway,
            above=3pt,
            edge weight,
            text=green!40!black
          ] {$10$};

        \draw[merger second]
          (vQ) -- (vR)
          node[
            midway,
            right=3pt,
            edge weight,
            text=orange
          ] {$20$};

        \draw[merger first]
          (vR) -- (vS)
          node[
            midway,
            below=3pt,
            edge weight,
            text=violet
          ] {$30$};

      \fi

      \node[graph vertex] at (vP) {$P$};
      \node[graph vertex] at (vQ) {$Q$};
      \node[graph vertex] at (vR) {$R$};
      \node[graph vertex] at (vS) {$S$};

      \ifnum\mode=0

        \node[
          font=\scriptsize,
          align=center,
          text width=5.1cm
        ] at (1.4,-1.10)
          {Vertex positions and drawn edge lengths
           have no meaning.\\
           Each pair's dissimilarity is given solely
           by its edge weight.};

      \else

        \node[font=\scriptsize]
          at (1.4,-1.00)
          {%
            \renewcommand{\arraystretch}{1.2}%
            \begin{tabular}{@{}c|c|c|c|c@{}}
              Threshold \(t\)
                & \(0\) & \(10\) & \(20\) & \(30\)\\
              \hline
              Number of components
                & \(4\) & \(3\) & \(2\) & \(1\)
            \end{tabular}%
          };

      \fi

    \end{scope}
  }

  %
  \begin{scope}[
    shift={(-0.1cm,-5.7cm)},
    x=0.095cm,
    y=0.95cm
  ]

    \node[font=\small,align=center]
      at (15,3.20) {(c) Ordinary \(H_0\) barcode};

    \node[font=\scriptsize]
      at (15,2.65)
      {All four components are born at \(t=0\).};

    \draw[landscape axis]
      (-1,0) -- (36.5,0)
      node[right] {$t$};

    \foreach \x in {0,5,...,30}{
      \draw[gray!80]
        (\x,0) -- (\x,-0.05);
    }

    \foreach \x in {0,10,20,30}{
      \node[below=3pt,font=\scriptsize]
        at (\x,0) {$\x$};
    }

    \foreach \x in {10,20,30}{
      \draw[landscape guide]
        (\x,0) -- (\x,1.70);
    }

    \draw[
      ->,
      gray!75,
      thick,
      dash pattern=on 4pt off 2pt
    ]
      (0,2.0) -- (35,2.0);

    \node[
      anchor=east,
      font=\small,
      text=gray!75
    ] at (-1.2,2.0) {$[0,\infty)$};

    \draw[merger first]
      (0,1.50) -- (30,1.50);

    \draw[merger second]
      (0,1.00) -- (20,1.00);

    \draw[merger third]
      (0,0.50) -- (10,0.50);

    \node[
      anchor=east,
      font=\small,
      text=violet
    ] at (-1.2,1.50) {$[0,d_1)$};

    \node[
      anchor=east,
      font=\small,
      text=orange
    ] at (-1.2,1.00) {$[0,d_2)$};

    \node[
      anchor=east,
      font=\small,
      text=green!40!black
    ] at (-1.2,0.50) {$[0,d_3)$};

    \node[font=\small]
      at (15,-0.75)
      {$
        {\color{violet}d_1=30},
        \quad
        {\color{orange}d_2=20},
        \quad
        {\color{green!40!black}d_3=10}
      $};

  \end{scope}

  %
  \begin{scope}[
    shift={(7cm,-5.7cm)},
    x=0.095cm,
    y=0.095cm
  ]

    \node[font=\small,align=center]
      at (15,32) {(d) Landscape of the finite bars};

    \node[font=\scriptsize]
      at (15,26.5)
      {The infinite bar is omitted.};

    \draw[landscape axis]
      (-1,0) -- (36.5,0)
      node[right] {$t$};

    \draw[landscape axis]
      (0,-1) -- (0,19)
      node[above,font=\small] {$\lambda_k$};

    \foreach \x in {0,5,...,30}{
      \draw[gray!80]
        (\x,0) -- (\x,-0.5);
    }

    \foreach \x in {0,10,20,30}{
      \node[
        below=3pt,
        font=\scriptsize,
        fill=white,
        inner sep=1pt
      ] at (\x,0) {$\x$};
    }

    \draw[gray!80]
      (0,10) -- (-0.5,10);
    \node[left=3pt,font=\scriptsize]
      at (0,10) {$10$};

    \foreach \x in {10,20,30}{
      \draw[landscape guide]
        (\x,0) -- (\x,17.5);
    }

    \draw[merger first]
      (0,0) -- (15,15) -- (30,0);

    \draw[merger second]
      (0,0) -- (10,10) -- (20,0);

    \draw[merger third]
      (0,0) -- (5,5) -- (10,0);

    \node[
      anchor=south west,
      inner sep=0pt,
      font=\small,
      text=violet
    ] at (20.5,9.5) {$\lambda_1$};

    \node[
      anchor=south west,
      inner sep=0pt,
      font=\small,
      text=orange
    ] at (13.5,6.5) {$\lambda_2$};

    \node[
      anchor=south west,
      inner sep=0pt,
      font=\small,
      text=green!40!black
    ] at (6.5,3.5) {$\lambda_3$};

    \node[font=\small]
      at (15,-7.5)
      {$\lambda_k(t)=\Lambda_{(0,d_k)}(t)$};

  \end{scope}

\end{tikzpicture}%
}
\caption{
From a weighted graph to its \(H_0\) persistence landscape.
(a) A complete graph on vertices \(P,Q,R,S\), with
dissimilarities specified solely by the edge weights;
the drawn geometry has no meaning.
(b) Kruskal's algorithm accepts the edges of weights
{\color{green!40!black}\(10\)}, {\color{orange}\(20\)}, {\color{violet}\(30\)},
each joining two previously
distinct components.
The number of components therefore decreases from
\(4\) to \(1\) at these thresholds.
(c) The finite barcode deaths are the accepted
edge weights, sorted as
\({\color{violet}d_1=30}\),
\({\color{orange}d_2=20}\), and
\({\color{green!40!black}d_3=10}\).
The infinite bar is shown as a grey arrow.
(d) Omitting that infinite bar gives three nested
landscape layers.
Across (b)--(d), solid \textcolor{violet}{violet} identifies death {\color{violet}\(30\)},
dashed \textcolor{orange}{orange} death {\color{orange}\(20\)}, and dash-dotted \textcolor{green!50!black}{dark green}
death {\color{green!40!black}\(10\)}.
The lower panels share the same horizontal scale;
the landscape panel uses equal horizontal and
vertical unit lengths.
}
\label{fig:mst-barcode-landscape}
\end{figure}

Consequently, at every parameter value where the sorted
death times are differentiable,
\textbf{Theorem~\ref{prop:closed} expresses the landscape
velocity entirely in terms of the sorted minimum-spanning-tree
edge weights and their derivatives.}
This gives a computational route that avoids constructing
landscapes and interprets their speed in terms of changing
component-merger thresholds.
The minimum spanning tree itself need not remain unchanged
as \(g\) varies; the differentiability requirement concerns
the sorted death times, not a fixed choice of tree edges.

\textbf{The same correspondence holds for the edge-threshold
filtration of any fixed connected graph} with positive edge
weights and all vertices present at filtration value \(0\):
the complete graph is not essential.

For finite Euclidean point clouds, \textbf{the MST description also
extends to \v{C}ech and unweighted alpha filtrations.}
When Vietoris--Rips is indexed by edge length and \v{C}ech
by ball radius, their finite \(H_0\) death times are,
respectively, the MST edge lengths and half those lengths
\cite[Sec.~2.2 and Lemma~20]{beers2024fiberspointcloudpersistence}.
Moreover, alpha and \v{C}ech filtrations have identical
persistence barcodes when indexed by the same radius
\cite{fraser-filtered-covers}.
Hence, under the same regularity assumptions, the MST-based
argument and the common-birth landscape-velocity formula
apply to both filtrations, with the corresponding
rescaling of death times.

\subsection{Landscape velocity as a Hilbert-space speed}
\label{subsec:l2-speed}

Differentiating a landscape at each filtration value gives a candidate
velocity vector. But does this vector describe the motion of the whole
landscape? Pointwise convergence alone does not answer the question:
the approximation must also hold in the Hilbert-space norm.

We therefore seek a first-order expansion
\[
    \lambda(g_0+h)
    =
    \lambda(g_0)+h\,\partial_g\lambda(g_0)+r_h,
    \qquad
    \lim_{h\to0}\frac{\|r_h\|_{\mathcal H}}{|h|}=0.
\]
Theorem~\ref{prop:l2-speed} proves this expansion under the hypotheses
of Theorem~\ref{prop:closed}. Consequently, the quantity computed from
the death times is the speed of the landscape curve:
\[
    \xi(g_0)
    =\|\partial_g\lambda(g_0)\|_{\mathcal H}
    =\lim_{h\to0}
    \frac{\|\lambda(g_0+h)-\lambda(g_0)\|_{\mathcal H}}{|h|}.
\]

\begin{theorem}[Hilbert-space differentiability and
identification of the speed]
\label{prop:l2-speed}
Under the hypotheses of Theorem~\ref{prop:closed},
the map \(g\mapsto\lambda(g)\) is differentiable
at \(g_0\) in the \(\mathcal H\) norm.
Its derivative is the sequence
\[
\partial_g\lambda(g_0)=\bigl(\partial_g\lambda_k(g_0;\cdot)\bigr)_{k\geq1},
\]
where, for almost every \(t\),
\[
\partial_g\lambda_k(g_0;t)
=
\begin{cases}
\displaystyle
\partial_g d_k(g_0)\,
\mathbf 1_{\left\{
\frac{d_k(g_0)}{2}<t<d_k(g_0)
\right\}},
& 1\leq k\leq L-1,\\[2ex]
0, & k\geq L.
\end{cases}
\]
In particular, the landscape velocity is the
 \(L^2\) (metric) speed:
\[
\xi(g_0)
=
\|\partial_g\lambda(g_0)\|_{\mathcal H}
=
\lim_{h\to0}
\frac{
\|\lambda(g_0+h)-\lambda(g_0)\|_{\mathcal H}
}{
\,|h|
}.
\]
\end{theorem}

\begin{proof}
Let $w=(w_k)_{k\geq1}$ be the sequence defined
in the statement, regarded as a candidate derivative.
Each nonzero \(w_k\) is bounded and supported
on a bounded interval, and only finitely many
components are nonzero.
Thus \(w\in\mathcal H\).

We will prove that
\[
\lim_{h\to0}
\left\|
\frac{\lambda(g_0+h)-\lambda(g_0)}{h}
-w
\right\|_{\mathcal H}
=0,
\]
which establishes that \(\partial_g\lambda(g_0)=w\).

For all sufficiently small nonzero \(h\), write
\[
q_{k,h}(t)
=
\frac{\lambda_k(g_0+h;t)-\lambda_k(g_0;t)}{h}.
\]
The proof of Theorem~\ref{prop:closed} shows that
\[
q_{k,h}(t)\longrightarrow w_k(t)
\qquad\text{for almost every }t.
\]
Moreover, the tent function is \(1\)-Lipschitz
in its death parameter (compare with the \(1\)-Lipschitz property
of landscape layers in \cite{bubenik}), so
\[
|q_{k,h}(t)|
\leq
\left|
\frac{d_k(g_0+h)-d_k(g_0)}{h}
\right|.
\]
Differentiability of the finitely many death times
at \(g_0\) therefore gives a constant \(C>0\)
such that, for all sufficiently small nonzero \(h\),
\[
|q_{k,h}(t)|\leq C,
\qquad
|w_k(t)|\leq C,
\qquad k=1,\dots,L-1.
\]
Continuity of these death times also gives \(R>0\)
such that all nearby tents, and all the \(w_k\),
are supported in \([0,R]\).
Consequently,
\[
|q_{k,h}(t)-w_k(t)|^2
\leq
4C^2\,\mathbf 1_{[0,R]}(t).
\]
The right-hand side is integrable.
Dominated convergence, followed by summation
over the finitely many nonzero layers, yields
\[
\lim_{h\to0}
\sum_{k=1}^{L-1}
\int_{\mathbb R}
|q_{k,h}(t)-w_k(t)|^2\,dt
=0.
\]
Equivalently,
\[
\lim_{h\to0}
\left\|
\frac{\lambda(g_0+h)-\lambda(g_0)}{h}
-w
\right\|_{\mathcal H}
=0.
\]
Thus \(\partial_g\lambda(g_0)=w\) in \(\mathcal H\).
By the definition of the landscape velocity,
\[
\xi(g_0)
=
\|w\|_{\mathcal H}.
\]
Since the difference quotients converge to \(w\) in
\(\mathcal H\), the reverse triangle inequality gives
\[
\left|
\frac{\|\lambda(g_0+h)-\lambda(g_0)\|_{\mathcal H}}{|h|}
-\|w\|_{\mathcal H}
\right|
\leq
\left\|
\frac{\lambda(g_0+h)-\lambda(g_0)}{h}-w
\right\|_{\mathcal H}
\longrightarrow 0.
\]
Thus the instantaneous speed of the landscape curve at
\(g_0\) is \(\|w\|_{\mathcal H}\).
Combining these identities proves the speed
interpretation.
Hence, under the hypotheses of
Theorem~\ref{prop:closed}, \(\xi(g_0)\) is
the  \(L^2\) speed of the landscape
curve at \(g_0\).

\end{proof}

Thus, under the hypotheses of Theorem~15, landscape velocity is
the \(L^2\) speed of the landscape curve in \(\mathcal{H}\).
In other words, it measures distance traveled in \(\mathcal{H}\)
per unit change in \(g\).

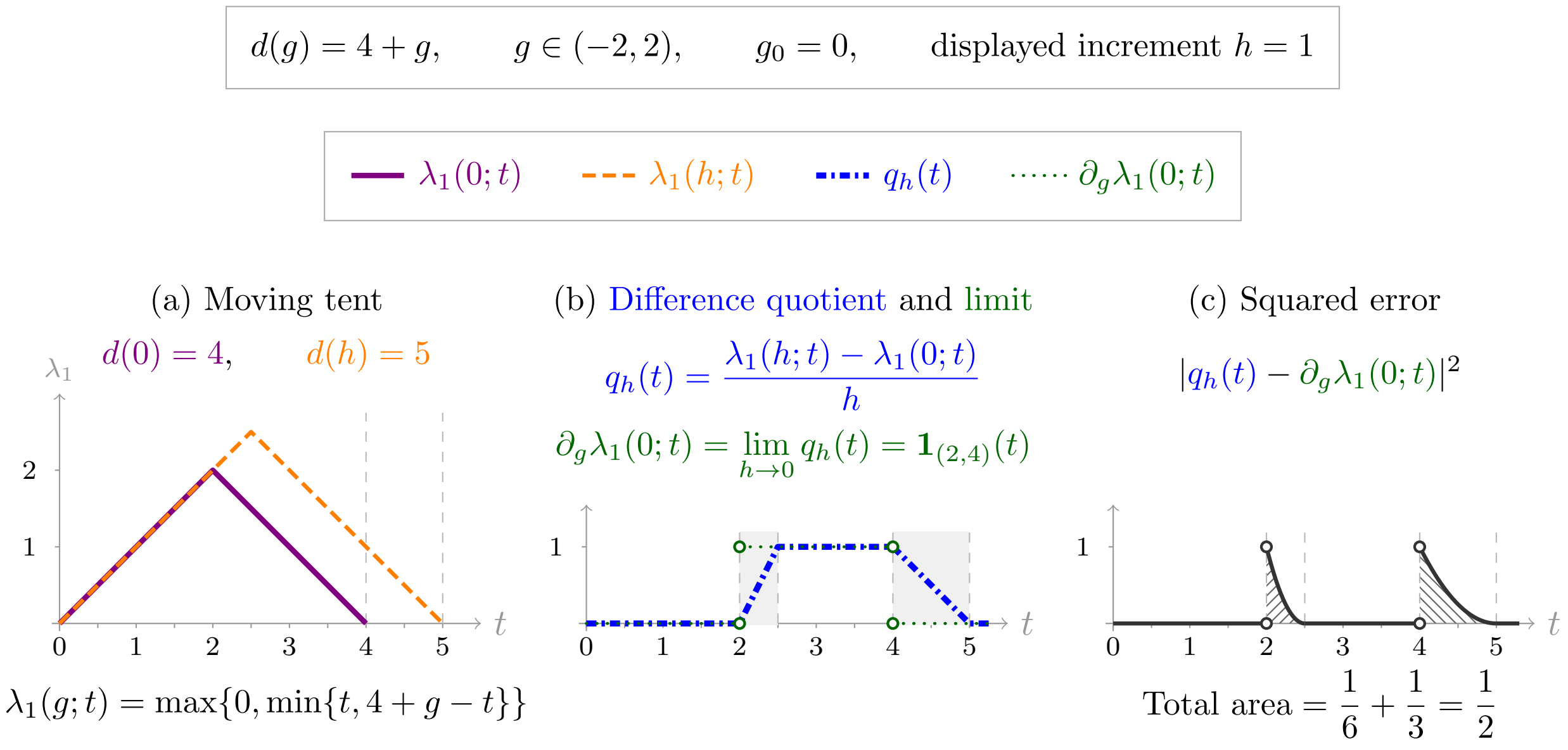
\begin{figure}[htbp]
\centering

\begingroup
\setlength{\fboxsep}{7pt}
\fcolorbox{gray!65}{white}{%
  \small
  \(
    d(g)=4+g,\qquad
    g\in(-2,2),\qquad
    g_0=0,\qquad
    \text{displayed increment }h=1
  \)%
}
\endgroup

\par\vspace{4mm}

\begingroup
\setlength{\fboxsep}{7pt}
\fcolorbox{gray!65}{white}{%
  \small
  \begin{tabular}{
    @{}l@{\hspace{6mm}}
       l@{\hspace{6mm}}
       l@{\hspace{6mm}}
       l@{}
  }

    \tikz[baseline=-0.6ex]{
      \draw[violet,line width=1.6pt,solid]
        (0,0) -- (0.55,0);
    }
    \({\color{violet}\lambda_1(0;t)}\)

    &

    \tikz[baseline=-0.6ex]{
      \draw[
        orange,
        line width=1.2pt,
        dash pattern=on 4pt off 2pt
      ]
        (0,0) -- (0.55,0);
    }
    \({\color{orange}\lambda_1(h;t)}\)

    &

    \tikz[baseline=-0.6ex]{
      \draw[
        blue,
        line width=1.8pt,
        dash pattern=on 3pt off 1.5pt on 1pt off 1.5pt
      ]
        (0,0) -- (0.55,0);
    }
    \({\color{blue}q_h(t)}\)

    &

    \tikz[baseline=-0.6ex]{
      \draw[
        green!40!black,
        line width=0.9pt,
        line cap=round,
        dash pattern=on 0pt off 3pt
      ]
        (0,0) -- (0.55,0);
    }
    \({\color{green!40!black}
      \partial_g\lambda_1(0;t)}\)

  \end{tabular}%
}
\endgroup

\par\vspace{5mm}

\makebox[\linewidth][c]{%
\begin{tikzpicture}[
  x=0.80cm,
  y=0.80cm,
  original tent/.style={
    violet,
    line width=1.6pt,
    solid,
    line join=round
  },
  perturbed tent/.style={
    orange,
    line width=1.2pt,
    dash pattern=on 4pt off 2pt,
    line join=round
  },
  difference quotient/.style={
    blue,
    line width=1.8pt,
    dash pattern=on 3pt off 1.5pt on 1pt off 1.5pt,
    line join=round
  },
  limiting derivative/.style={
    green!40!black,
    line width=0.9pt,
    line cap=round,
    dash pattern=on 0pt off 3pt
  },
  squared error/.style={
    black!80,
    line width=1.2pt
  }
]

  \foreach \xshift/\mode/\heading in {
    0/1/{(a) Moving tent},
    5.5/2/{(b) {\color{blue}Difference quotient} and {\color{green!40!black}limit}},
    11/3/{(c) Squared error}
  }{

    \begin{scope}[shift={(\xshift cm,0)}]

      \node[font=\small,align=center]
        at (2.7,4.20) {\heading};

      \draw[landscape axis]
        (-0.10,0) -- (5.50,0)
        node[right] {$t$};

      \ifnum\mode=1

        \draw[landscape axis]
          (0,-0.10) -- (0,3.00)
          node[above,font=\scriptsize] {$\lambda_1$};

        \foreach \y in {1,2}{
          \draw[gray!80]
            (0,\y) -- (-0.06,\y);
          \node[left=3pt,font=\scriptsize]
            at (0,\y) {$\y$};
        }

      \else

        \draw[landscape axis]
          (0,-0.10) -- (0,1.55);

        \draw[gray!80]
          (0,1) -- (-0.06,1);
        \node[left=3pt,font=\scriptsize]
          at (0,1) {$1$};

      \fi

      \foreach \x in {0,0.5,...,5}{
        \draw[gray!80]
          (\x,0) -- (\x,-0.06);
      }

      \foreach \x in {0,1,...,5}{
        \node[
          below=3pt,
          font=\scriptsize,
          fill=white,
          inner sep=1pt
        ] at (\x,0) {$\x$};
      }

    \end{scope}
  }

  \begin{scope}

    \node[font=\small]
      at (2.7,3.50)
      {$
        {\color{violet}d(0)=4},
        \qquad
        {\color{orange}d(h)=5}
      $};

    \draw[landscape guide]
      (4,0) -- (4,2.75);
    \draw[landscape guide]
      (5,0) -- (5,2.75);

    \draw[original tent]
      (0,0) -- (2,2) -- (4,0);

    \draw[perturbed tent]
      (0,0) -- (2.5,2.5) -- (5,0);

    \node[font=\small]
      at (2.7,-1.05)
      {$\lambda_1(g;t)
        =\max\{0,\min\{t,4+g-t\}\}$};

  \end{scope}

  \begin{scope}[shift={(5.5cm,0)}]

    \node[font=\small,text=blue]
      at (2.7,3.25)
      {$q_h(t)
        =\dfrac{\lambda_1(h;t)-\lambda_1(0;t)}{h}$};

    \fill[gray!12]
      (2,-0.02) rectangle (2.5,1.20);

    \fill[gray!12]
      (4,-0.02) rectangle (5,1.20);

    \foreach \x in {2,2.5,4,5}{
      \draw[landscape guide]
        (\x,0) -- (\x,1.25);
    }

    \draw[difference quotient]
      (0,0)
      -- (2,0)
      -- (2.5,1)
      -- (4,1)
      -- (5,0)
      -- (5.30,0);

    \draw[limiting derivative]
      (0,0) -- (2,0);

    \draw[limiting derivative]
      (2,1) -- (4,1);

    \draw[limiting derivative]
      (4,0) -- (5.30,0);

    \foreach \x/\y in {2/0,2/1,4/1,4/0}{
      \draw[
        draw=green!40!black,
        fill=white,
        line width=0.8pt
      ]
        (\x,\y) circle[radius=1.5pt];
    }

    \node[
      font=\small,
      text=green!40!black
    ] at (2.7,2.2)
      {$\displaystyle\partial_g\lambda_1(0;t) = \lim_{h\to0}q_h(t)
        =\mathbf1_{(2,4)}(t)$};

  \end{scope}

  \begin{scope}[shift={(11cm,0)}]

    \node[font=\small]
      at (2.7,3.25)
      {$|{\color{blue}q_h(t)}-{\color{green!40!black}\partial_g\lambda_1(0;t)}|^2$};

    \path[
      pattern=north east lines,
      pattern color=black!60
    ]
      (2,0)
      -- (2,1)
      -- plot[
        domain=2:2.5,
        samples=35
      ] (\x,{(5-2*\x)^2})
      -- cycle;

    \path[
      pattern=north west lines,
      pattern color=black!60
    ]
      (4,0)
      -- (4,1)
      -- plot[
        domain=4:5,
        samples=45
      ] (\x,{(5-\x)^2})
      -- cycle;

    \foreach \x in {2,2.5,4,5}{
      \draw[landscape guide]
        (\x,0) -- (\x,1.20);
    }

    \draw[squared error]
      (0,0) -- (2,0);

    \draw[squared error]
      plot[
        domain=2:2.5,
        samples=35
      ] (\x,{(5-2*\x)^2});

    \draw[squared error]
      (2.5,0) -- (4,0);

    \draw[squared error]
      plot[
        domain=4:5,
        samples=45
      ] (\x,{(5-\x)^2});

    \draw[squared error]
      (5,0) -- (5.30,0);

    \foreach \x/\y in {2/0,2/1,4/0,4/1}{
      \draw[
        draw=black!80,
        fill=white,
        line width=0.8pt
      ]
        (\x,\y) circle[radius=1.5pt];
    }

    \node[font=\small]
      at (2.7,-1.05)
      {$\displaystyle
        \text{Total area}
        =\frac16+\frac13=\frac12$};

  \end{scope}

\end{tikzpicture}%
}
\caption{
Difference-quotient convergence for the common-birth
bar with death \(d(g)=4+g\), evaluated at \(g=0\).
The displayed increment is \(h=1\).
(a) The original tent is solid \textcolor{violet}{violet}; the perturbed
tent is dashed \textcolor{orange}{orange}.
(b) The difference quotient \(q_h\) is dash-dotted \textcolor{blue}{blue}.
For almost every fixed \(t\), its two-sided limit is
\(\lim_{h\to0}q_h(t)=\partial_g\lambda_1(0;t)=\mathbf1_{(2,4)}(t)\),
shown in dotted \textcolor{green!50!black}{dark green}; the two exceptional
filtration values are \(t=2\) and \(t=4\).
The displayed positive increment illustrates the
approach from \(h>0\).
The \textcolor{gray}{grey} strips mark the only regions where the
displayed difference quotient and its limit differ,
apart from breakpoint values.
For general \(0<h<2\), these regions are
\((2,2+h/2)\) and \((4,4+h)\).
(c) The integral of the squared error is
\(h/6+h/3=h/2\), which tends to zero as \(h\to0^+\).
The two contributions are distinguished by diagonal
lines sloping in opposite directions.
The left-sided difference quotients converge to the
same function in \(L^2\), establishing the two-sided
\(L^2\) derivative.
Open circles omit breakpoint values, which do not
affect the norm.
All panels use the same coordinate scales.
}
\label{fig:landscape-l2-differentiability}
\end{figure}

Figure~\ref{fig:landscape-l2-differentiability}
illustrates the difference-quotient convergence
in Theorem~\ref{prop:l2-speed}:
the error is confined to shrinking neighbourhoods
of the moving tent breakpoints, and its squared
\(L^2\) norm tends to zero.

\subsection{Metric speed without differentiability at a crossing}
\label{subsec:crossings}

\begin{figure}[htbp]
\centering
\begin{tikzpicture}[
  crossing heading/.style={font=\small,align=center},
  crossing tick/.style={font=\scriptsize},
  crossing legend/.style={font=\small,anchor=west},
  texture blue/.style={solid},
  texture red/.style={dash pattern=on 4pt off 2pt},
  texture teal/.style={
    dash pattern=on 3pt off 1.5pt on 1pt off 1.5pt
  },
  texture violet/.style={solid},
  texture orange/.style={dash pattern=on 4pt off 2pt},
  texture green/.style={
    dash pattern=on 3pt off 1.5pt on 1pt off 1.5pt
  },
  texture black/.style={dash pattern=on 1pt off 1.2pt}
]

\begin{scope}[x=1.2cm,y=1.2cm]
  \node[crossing heading] at (0,2.55) {(a) Unsorted deaths};

  \draw[landscape axis]
    (-1.2,0)--(1.35,0) node[right] {$g$};
  \draw[landscape axis]
    (-1.2,-0.05)--(-1.2,2.2);
  \draw[landscape guide]
    (0,0)--(0,2.15);

  \foreach \x in {-1,0,1}{
    \draw (\x,0)--(\x,-0.05);
    \node[below=3pt,crossing tick] at (\x,0) {$\x$};
  }
  \foreach \y in {0,1,2}{
    \draw (-1.2,\y)--(-1.25,\y);
    \node[left=3pt,crossing tick] at (-1.2,\y) {$\y$};
  }

  \draw[landscape curve,blue,texture blue]
    (-1,2)--(1,0);
  \draw[landscape curve,red,texture red]
    (-1,0)--(1,2);
  \draw[fill=white,thick]
    (0,1) circle[radius=1.8pt];

  \node[font=\small,text=blue]
    at (0,-0.57) {$a(g)=1-g$};
  \node[font=\small,text=red]
    at (0,-0.92) {$b(g)=1+g$};
\end{scope}

\begin{scope}[shift={(6.8cm,0)},x=1.2cm,y=1.2cm]
  \node[crossing heading] at (0,2.55) {(b) Sorted deaths};

  \draw[landscape axis]
    (-1.2,0)--(1.35,0) node[right] {$g$};
  \draw[landscape axis]
    (-1.2,-0.05)--(-1.2,2.2);
  \draw[landscape guide]
    (0,0)--(0,2.15);

  \foreach \x in {-1,0,1}{
    \draw (\x,0)--(\x,-0.05);
    \node[below=3pt,crossing tick] at (\x,0) {$\x$};
  }
  \foreach \y in {0,1,2}{
    \draw (-1.2,\y)--(-1.25,\y);
    \node[left=3pt,crossing tick] at (-1.2,\y) {$\y$};
  }

  \draw[landscape curve,violet,texture violet]
    (-1,2)--(0,1)--(1,2);
  \draw[landscape curve,orange,texture orange]
    (-1,0)--(0,1)--(1,0);
  \draw[fill=white,thick]
    (0,1) circle[radius=1.8pt];

  \node[font=\small,text=violet]
    at (0,-0.57) {$d_1(g)=1+|g|$};
  \node[font=\small,text=orange]
    at (0,-0.92) {$d_2(g)=1-|g|$};
\end{scope}

\begin{scope}[shift={(-1cm,-4.6cm)},x=2.1cm,y=2.1cm]
  \node[crossing heading]
    at (0.55,1.15) {(c) Barcode snapshots};

  \draw[landscape axis]
    (-0.05,-0.18)--(1.4,-0.18) node[right] {$t$};

  \foreach \x in {0,1}{
    \draw (\x,-0.18)--(\x,-0.22);
    \node[below=2pt,crossing tick]
      at (\x,-0.22) {$\x$};
  }

  \draw[landscape guide]
    (1,-0.08)--(1,0.98);

  \node[anchor=east,font=\small]
    at (-0.12,0.85) {$g=-\tfrac14$};
  \draw[landscape curve,blue,texture blue]
    (0,0.90)--(1.25,0.90);
  \draw[landscape curve,red,texture red]
    (0,0.80)--(0.75,0.80);

  \node[anchor=east,font=\small]
    at (-0.12,0.45) {$g=0$};
  \draw[landscape curve,blue,texture blue]
    (0,0.50)--(1,0.50);
  \draw[landscape curve,red,texture red]
    (0,0.40)--(1,0.40);

  \node[anchor=east,font=\small]
    at (-0.12,0.05) {$g=\tfrac14$};
  \draw[landscape curve,blue,texture blue]
    (0,0.10)--(0.75,0.10);
  \draw[landscape curve,red,texture red]
    (0,0)--(1.25,0);
\end{scope}

\begin{scope}[shift={(5.7cm,-4.6cm)},x=2.1cm,y=2.1cm]
  \node[crossing heading]
    at (0.55,1.15) {(d) Landscape snapshots};

  \draw[landscape axis]
    (0,0)--(1.4,0) node[right] {$t$};
  \draw[landscape axis]
    (0,-0.04)--(0,0.8)
    node[above] {$\lambda_k(g;t)$};

  \draw[landscape guide]
    (0.5,0)--(0.5,0.72);
  \draw[landscape guide]
    (1,0)--(1,0.72);

  \foreach \x/\lab in {0/0,0.5/{\frac12},1/1}{
    \draw (\x,0)--(\x,-0.025);
    \node[below=3pt,crossing tick]
      at (\x,0) {$\lab$};
  }

  \draw[landscape curve,black,texture black]
    (0,0)--(0.5,0.5)--(1,0);

  \draw[landscape curve,violet,texture violet]
    (0,0)--(0.625,0.625)--(1.25,0);
  \draw[landscape curve,orange,texture orange]
    (0,0)--(0.375,0.375)--(0.75,0);

  \node[font=\scriptsize] at (0.65,-0.38)
    {dotted: $g=0$; solid/dashed: $g=\pm\tfrac14$};
\end{scope}

%
\foreach \side/\xshift/\sgn/\axisstart/\heading in {
  -/-1/-1/-0.05/{(e) Left derivative $\partial_g^-\lambda(0)$},
  +/5.7/1/0/{(f) Right derivative $\partial_g^+\lambda(0)$}
}{
  \begin{scope}[
    shift={(\xshift cm,-8.6cm)},
    x=2.1cm,
    y=1.3cm
  ]
    \node[crossing heading]
      at (0.55,1.65) {\heading};

    \draw[landscape axis]
      (\axisstart,0)--(1.4,0) node[right] {$t$};
    \draw[landscape axis]
      (0,-1.25)--(0,1.28);

    \draw[landscape guide]
      (0.5,-1.13)--(0.5,1.13);
    \draw[landscape guide]
      (1,-1.13)--(1,1.13);

    \foreach \y in {-1,1}{
      \draw (0,\y)--(-0.028,\y);
      \node[left=3pt,crossing tick]
        at (0,\y) {$\y$};
    }

    \node[below left=2pt,crossing tick]
      at (0,0) {$0$};

    \foreach \x/\lab in {0.5/{\frac12},1/1}{
      \draw (\x,0)--(\x,-0.04);
      \node[
        below=3pt,
        crossing tick,
        fill=white,
        inner sep=1pt
      ] at (\x,-1.3) {$\lab$};
    }

    \draw[landscape curve,violet,texture violet]
      (0,0)--(0.5,0);
    \draw[landscape curve,violet,texture violet]
      (1,0)--(1.35,0);
    \draw[landscape curve,orange,texture orange]
      (0,0)--(0.5,0);
    \draw[landscape curve,orange,texture orange]
      (1,0)--(1.35,0);

    \draw[landscape curve,violet,texture violet]
      (0.5,\sgn)--(1,\sgn);
    \draw[landscape curve,orange,texture orange]
      (0.5,{-\sgn})--(1,{-\sgn});

    \foreach \x in {0.5,1}{
      \draw[draw=violet,fill=white,thick]
        (\x,\sgn) circle[radius=1.6pt];
      \draw[draw=orange,fill=white,thick]
        (\x,{-\sgn}) circle[radius=1.6pt];
      \draw[fill=white,thick]
        (\x,0) circle[radius=1.6pt];
    }

    \node[anchor=west,font=\small,text=violet]
      at (1.08,\sgn)
      {$\partial_g^{\side}\lambda_1(0)$};

    \node[anchor=west,font=\small,text=orange]
      at (1.08,{-\sgn})
      {$\partial_g^{\side}\lambda_2(0)$};
  \end{scope}
}
\end{tikzpicture}
\caption{
Crossing in Example~\ref{ex:crossing}, 
where $a(g)$ and $b(g)$ denote the death times of the two finite bars in the barcode at parameter $g$.
(a) Smooth unsorted deaths cross at \(g=0\).
(b) The sorted deaths develop corners.
(c) Barcode snapshots: {\color{blue}blue} denotes {\color{blue}\([0,a(g))\)}, {\color{red} red
\([0,b(g))\)}; the bar multiset is the same at \(g=\pm1/4\).
(d) At \(g=0\), the two layers coincide (dotted black).
At \(g=\pm1/4\), the common landscape has first layer
\textcolor{violet}{violet} and second layer \textcolor{orange}{orange}.
(e)--(f) The left and right derivative fields at \(g=0\)
differ, although their norms agree. Since $\partial_g^\pm\lambda_k(0;t)
=\partial_g^\pm d_k(0)\,
\mathbf1_{(1/2,1)}(t),$ we can read the one-sided derivatives directly from (b).
\textcolor{violet}{Violet} denotes the first
component, \textcolor{orange}{orange} the second. They are nonzero only for
filtration values \(1/2<t<1\). Open circles omit breakpoint
values, which do not affect the \(L^2\) classes.
Colors are paired with distinct line textures:
\textcolor{blue}{blue} and \textcolor{violet}{violet} are solid, \textcolor{red}{red} and \textcolor{orange}{orange} are dashed,
and the coincident \(g=0\) landscape in (d) is dotted black.
The horizontal axes in (a)--(b) represent \(g\); those in
(c)--(f) represent \(t\). Panel scales need not agree.
}
\label{fig:crossing}
\end{figure}
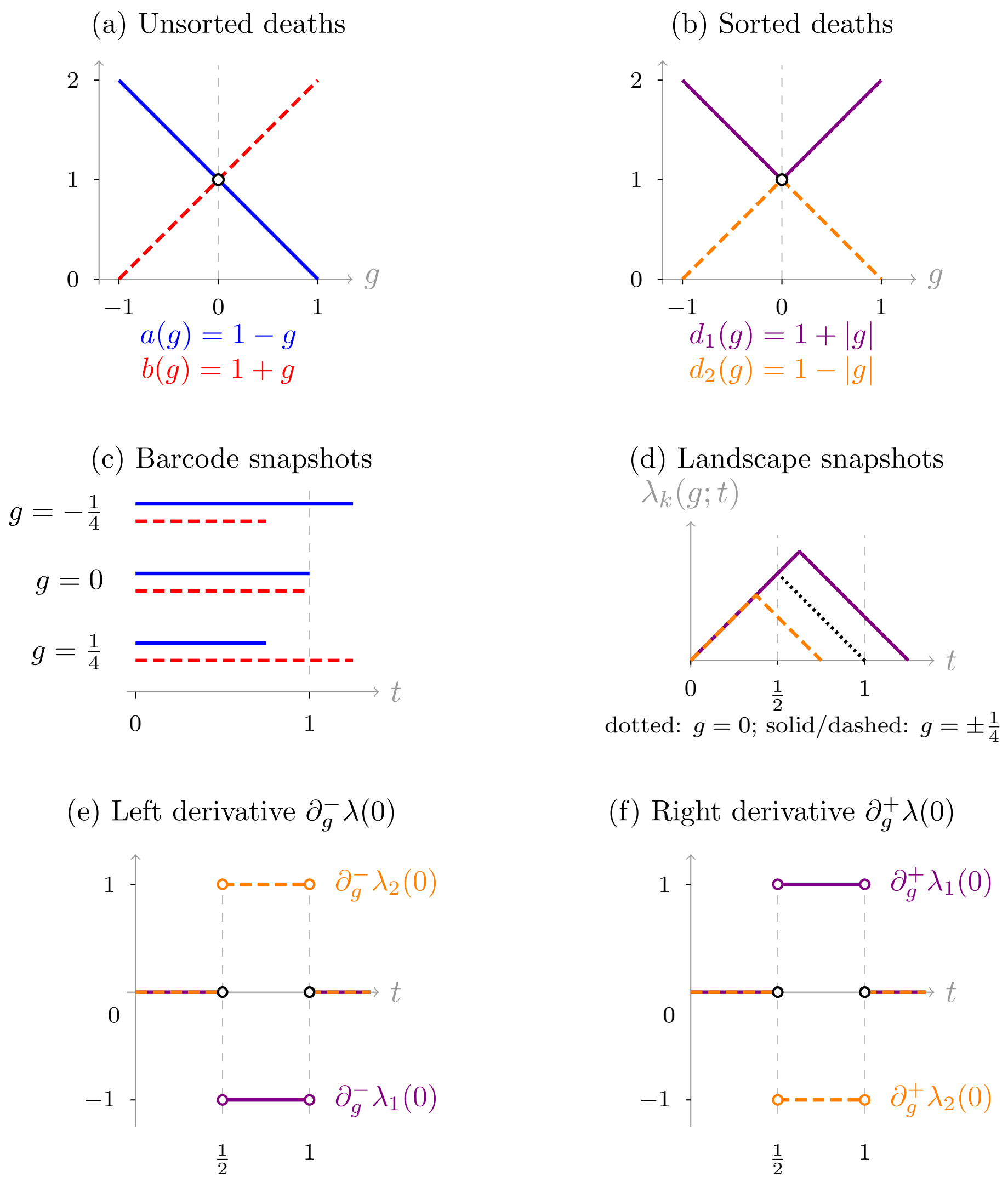

The differentiability assumption in
Theorem~\ref{prop:closed} is not automatic at every
parameter value, even when the barcode admits a labeling
in which all death times vary smoothly.
When two deaths cross (as two vines cross in a vineyard \cite{vineyards}), sorting them can create corners
in the sorted death functions and, consequently, prevent
the landscape curve from having a two-sided derivative
in \(\mathcal H\).
The following example illustrates this obstruction and
shows how the corresponding one-sided derivatives can
still be computed and used to obtain one-sided
landscape velocities through the same closed formula. Figure~\ref{fig:crossing} illustrates the deaths, barcodes,
landscapes, and one-sided derivative fields of Example \ref{ex:crossing}.

\begin{example}[A crossing of barcode bars]
\label{ex:crossing}
Let \(g\in(-1,1)\), and consider the two common-birth bars
\[
\mathcal B_g^{\mathrm{fin}}
=\{[0,a(g)),[0,b(g))\},
\qquad a(g)=1-g,\quad b(g)=1+g.
\]
Here \(g\) parametrizes the barcode family, whereas \(t\)
denotes the filtration variable. With this definition, the death times \(a\) and \(b\)
are smooth and cross at \(g=0\). Sorting them in nonincreasing
order gives
\[
d_1(g)=\max\{a(g),b(g)\}=1+|g|,
\qquad
d_2(g)=\min\{a(g),b(g)\}=1-|g|.
\]
Thus \(d_1=a\), \(d_2=b\) for \(g<0\), whereas
\(d_1=b\), \(d_2=a\) for \(g>0\). Both sorted deaths
have a corner at \(0\). Their one-sided derivatives are
\[
\partial_g^-d_1(0)=-1,\quad \partial_g^+d_1(0)=1,
\qquad
\partial_g^-d_2(0)=1,\quad \partial_g^+d_2(0)=-1,
\]
where, for a scalar function \(f\), we use the convention
\[
\partial_g^\pm f(0)
:=\lim_{h\to0^\pm}\frac{f(h)-f(0)}{h}.
\]
In particular, the denominator of the left difference quotient
is \(h<0\), not \(|h|\).

By the common-birth property, the two nonzero landscape layers are
\[
\lambda_k(g;t)=\max\{0,\min\{t,d_k(g)-t\}\},
\qquad k=1,2,
\]
and \(\lambda_k=0\) for \(k\geq3\).
At \(g=0\), both layers coincide with the tent
\[
\lambda_1(0;t)=\lambda_2(0;t)
=\max\{0,\min\{t,1-t\}\}.
\]
Its rising branch occupies \(0<t<1/2\), and its falling
branch occupies \(1/2<t<1\).

We can compute the one-sided layer derivatives by exactly the
same branchwise argument as in Theorem~\ref{prop:closed}.
Fix \(k\in\{1,2\}\), write \(d=d_k(0)=1\) and
\(d_h=d_k(h)\), and hold \(t\) fixed. If \(d/2<t<d\),
continuity of \(d_k\) implies that \(t\) remains on the
falling branch for all sufficiently small \(|h|\). Hence
\[
\lambda_k(h;t)=d_h-t,\qquad \lambda_k(0;t)=d-t,
\]
and therefore
\[
\partial_g^\pm\lambda_k(0;t)
=\lim_{h\to0^\pm}
\frac{\lambda_k(h;t)-\lambda_k(0;t)}{h}
=\lim_{h\to0^\pm}\frac{d_h-d}{h}
=\partial_g^\pm d_k(0).
\]
This equality applies on the falling branch, not at every
filtration value. For fixed \(0<t<1/2\), the layer remains
on its rising branch and equals \(t\), independently of \(g\).
For fixed \(t<0\) or \(t>1\), it vanishes for all sufficiently
small \(|g|\). Consequently, for almost every \(t\),
\[
\partial_g^\pm\lambda_k(0;t)
=\partial_g^\pm d_k(0)\,
\mathbf1_{(1/2,1)}(t),\qquad k=1,2.
\]
The values at the two breakpoints \(t=1/2\) and \(t=1\)
do not affect the associated \(L^2\) functions.
The intervals \((-1,0)\) and \((0,1)\) describe the two
sides of the crossing in the family parameter \(g\);
the interval \((1/2,1)\) instead specifies the filtration
values \(t\) on the falling side of the tent at \(g=0\).

Define the one-sided derivatives of the full landscape by
\[
\partial_g^\pm\lambda(0)
:=\lim_{h\to0^\pm}\frac{\lambda(h)-\lambda(0)}{h}
\quad\text{in }\mathcal H.
\]
These limits exist by the one-sided version of the
dominated-convergence argument in Theorem~\ref{prop:l2-speed}:
the difference quotients are uniformly bounded, have common
compact support for small \(|h|\), and converge almost everywhere
on each side. With
\[
A(t)=\mathbf1_{(1/2,1)}(t),
\]
their components are the functions
\[
\partial_g^-\lambda_1(0;t)=-A(t),\quad \partial_g^-\lambda_2(0;t)=A(t),
\qquad
\partial_g^+\lambda_1(0;t)=A(t),\quad \partial_g^+\lambda_2(0;t)=-A(t),
\]
with all higher components zero. Equivalently,
\[
\partial_g^-\lambda(0)=(-A,A,0,\dots),\qquad
\partial_g^+\lambda(0)=(A,-A,0,\dots).
\]
They differ as elements of \(\mathcal H\), so the two-sided
\(\mathcal H\) derivative does not exist.

The failure of pointwise differentiability can also be seen
directly. For every fixed \(t\in(1/2,1)\) and all sufficiently
small \(|g|\),
\[
\lambda_1(g;t)=1+|g|-t,\qquad
\lambda_2(g;t)=1-|g|-t.
\]
For example, it suffices to take
\(|g|<\min\{2t-1,1-t\}\). Thus
\[
\frac{\lambda_1(h;t)-\lambda_1(0;t)}{h}
=\frac{|h|}{h},\qquad
\frac{\lambda_2(h;t)-\lambda_2(0;t)}{h}
=-\frac{|h|}{h},
\]
and neither two-sided limit exists. The layer derivatives
therefore fail on a set of filtration values of positive
Lebesgue measure, not merely at the tent breakpoints.
Definition~\ref{def:velocity} consequently does not define
\(\xi(0)\).

Thus smoothness of a labeling of the barcode deaths
does not suffice to apply Theorem~\ref{prop:closed}
at a crossing: the sorted deaths may fail to be
differentiable, and the landscape may have distinct
left and right derivatives.
\textbf{Nevertheless, the proof of the proposition applies
one-sidedly and computes these derivatives explicitly.}
\end{example}

\begin{remark}[Numerical handling of death-time crossings]
\label{rem:crossings-numerical}
In practice, when computing $\xi$ numerically and death times cross at some $g_0$, 
one can use one-sided finite differences (forward on $[g_0, g_0+\epsilon]$ or backward 
on $[g_0-\epsilon, g_0]$) on either side of the crossing. The formula for $\xi$ remains 
well-defined almost everywhere despite the non-differentiability at crossings.
\end{remark}

\begin{remark}[Vector differentiability and metric speed]
\label{rem:metric-speed}The failure of a two-sided landscape derivative does
not necessarily imply the failure of a distance-based
scalar speed.
In Example~\ref{ex:crossing}, the two one-sided derivative
vectors differ, but their norms (the one-sided landscape velocities) agree.
Indeed, since
\(\displaystyle \int_{\mathbb R}|A(t)|^2\,dt=\frac12\),
\[
\|\partial_g^-\lambda(0)\|_{\mathcal H}^2
=\|\partial_g^+\lambda(0)\|_{\mathcal H}^2=1.
\]
Hence both one-sided landscape velocities equal $1$, also obtained from the one-sided
closed formula
\[
\|\partial_g^\pm\lambda(0)\|_{\mathcal H}
=\left(\frac12\sum_{k=1}^2
 d_k(0)\bigl(\partial_g^\pm d_k(0)\bigr)^2\right)^{1/2}
=1.
\]
It follows that the distance-based limit (the metric speed, which is the metric derivative of the curve \(g\mapsto\lambda(g)\) in the sense of \cite[Chapter~1]{ambrosio-gigli-savare}) exists:
\[
\lim_{h\to0}
\frac{\|\lambda(h)-\lambda(0)\|_{\mathcal H}}{\,|h|}
=1,
\]
even though the two-sided vector derivative does not.
Moreover, \(\lambda(h)=\lambda(-h)\), so the centred vector
quotient \((\lambda(h)-\lambda(-h))/(2h)\) is identically zero
and does not recover this metric speed at the crossing.

Therefore, if we denote the  one-sided landscape velocities by
\[
\xi_\pm(0)
:=
\|\partial_g^\pm\lambda(0)\|_{\mathcal H},
\]
then the one-sided version of the closed formula gives
\[
\xi_\pm(0)
=
\left(
\frac12
\sum_{k=1}^{2}
d_k(0)\bigl(\partial_g^\pm d_k(0)\bigr)^2
\right)^{1/2}.
\]
Both equal $1$, although the two-sided landscape
derivative does not exist.
\end{remark}

Example~\ref{ex:crossing} shows why symmetry of the closed-form
sum alone cannot justify a two-sided derivative. In general,
one-sided vector derivatives, when they exist, can differ; their
norms may agree or differ. Theorem~\ref{prop:closed} and
Theorem~\ref{prop:l2-speed} apply at differentiability points,
and Proposition~\ref{prop:regularity} supplies such points almost
everywhere under its hypotheses.

\section{Application: the transverse-field Ising chain}
\label{sec:application}

Can a single topological response locate a quantum phase transition without
changing the homology degree or the summary when the boundary
conditions change? We test this question on the transverse-field
Ising chain, revisiting Olsthoorn's analysis~\cite{olsthoorn-entanglement}
with landscape velocity in place of his summaries.

Olsthoorn uses the most persistent \(H_1\) lifetime for periodic
boundaries but the range of finite \(H_0\) death times for open
boundaries. We use the same prescription throughout: compute \(H_0\)
persistence and estimate the critical point by the position of the
maximum of \(\xi\). We test this prescription with both
mutual-information and correlator dissimilarities.

As an additional candidate to give his method more chances to succeed in locating the critical point,
for open boundaries, we also examine the corner following the
maximum of the death-time range, since it lies closer to the critical
point than the maximum itself.
 
These comparisons distinguish two
questions: whether an estimator lies near the transition at a given
system size, and whether it approaches the critical point as the
system grows. The advantages of landscape velocity are collected in
Section~\ref{subsec:app-advantages}.

\subsection{Model, ground state and distances}
\label{subsec:app-model}

\begin{figure}[htbp]
\centering
\includegraphics[width=\textwidth]{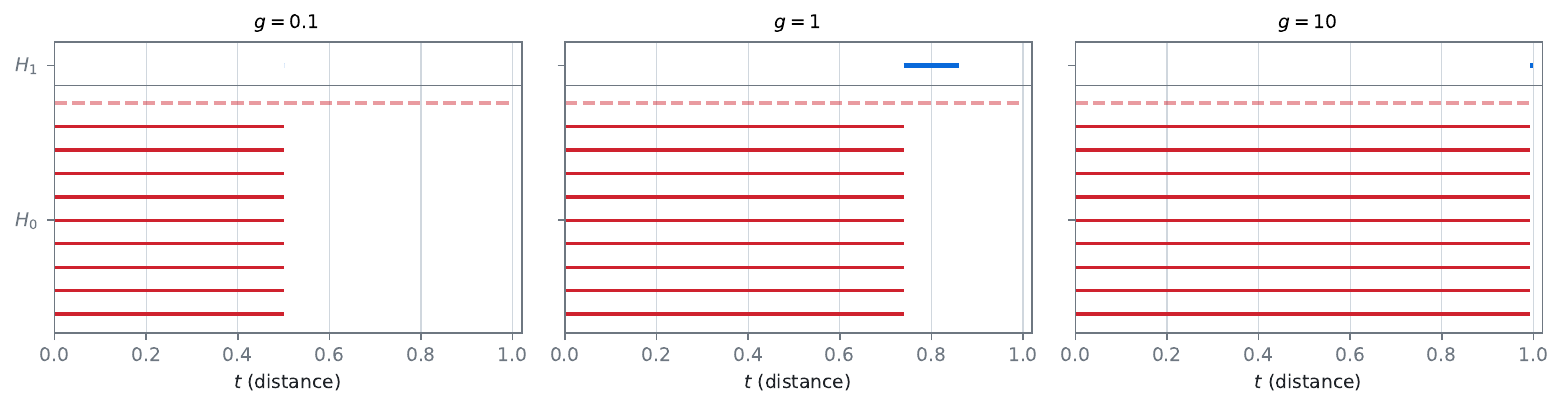}
\caption[Barcodes of the periodic chain]{\textbf{Barcodes of the periodic chain.}
The three panels reproduce Olsthoorn's Figs.~6(a)--(c)
\cite{olsthoorn-entanglement}.
They show the barcodes of the Vietoris--Rips filtration for \(L=10\)
and the distance \(D^{\mathrm{MI}}=1-\frac{M_{ij}}{2\ln2}\), at
\(g=0.1\) (ferromagnet), \(g=1\) (critical point) and \(g=10\)
(paramagnet), from left to right.
The \textcolor{barred}{red} lines are the \(L-1=9\) finite \(H_0\)
bars, the \textcolor{barlightred}{light red, dashed} line is the infinite
\(H_0\) bar, which we draw up to \(\epsilon=1\), and the
\textcolor{cycleblue}{blue} line is the most persistent \(H_1\) bar.
In each panel all nine finite \(H_0\) bars die at the same value, the
nearest-neighbour distance \(d_1(g)\), which equals
\(d_1(0.1)=0.5018\), \(d_1(1)=0.7396\) and \(d_1(10)=0.9937\); the
caption of Figure~\ref{fig:pbc6} explains why.
The \(H_1\) bar is clearly visible only at the critical point, where
it is born at \(d_1(1)=0.7396\) and dies at \(0.8599\).
Its length is \(1.1\times10^{-5}\) at \(g=0.1\) and
\(6.3\times10^{-3}\) at \(g=10\).}
\label{fig:barcodes}
\end{figure}

\paragraph{Hamiltonian.}
We use Olsthoorn's Hamiltonian \cite{olsthoorn-entanglement},
\begin{equation}
\label{eq:app-hamiltonian}
H(g)
=
-J\sum_{i=1}^{L_{\mathrm b}}\sigma^x_i\sigma^x_{i+1}
-h\sum_{i=1}^{L}\sigma^z_i,
\qquad
g=\frac hJ,
\qquad
J=1,
\end{equation}
on a chain of \(L\) sites, where \(h\) is the transverse field, \(L_{\mathrm b}=L-1\) for open boundary conditions (OBC)
and \(L_{\mathrm b}=L\) for periodic boundary conditions (PBC).
Here, $L$ is the number of spin-1/2 particles (qubits) with Hilbert space $(\mathbb{C}^2)^{\otimes L}$, and $\sigma_i^x$ and $\sigma_i^z$ are the Pauli matrices acting on site $i$.
For \(J>0\), the interaction term energetically favours alignment
of neighbouring spins along the \(x\)-direction. The transverse field \(h\),
applied uniformly to every site, favours polarization along \(z\) and
induces quantum fluctuations. The competition between the interaction
and field energy scales is governed by the dimensionless ratio \(g=h/J\);
in this paper we set \(J=1\), so \(g=h\).
In the periodic case, \(L\) is even and we identify
\(\sigma_{L+1}^{x}=\sigma_1^{x}\).
Therefore, the first and the last spins interact exactly as two
consecutive neighbours do, which gives the chain a circular
structure: the sites form a ring, and so do the nearest-neighbour
edges of the Vietoris--Rips filtration.
This ring is what the most persistent \(H_1\) bar detects: it is born
when the nearest-neighbour edges close the cycle, and it dies when
the cycle is filled by triangles (Figure~\ref{fig:barcodes}).
With open boundary conditions (OBC) there is no closing bond, the
chain is a path, and there is no ring to detect, so there is no
\(H_1\) bar to measure.
We take \(g\ge0\), as Olsthoorn does.
This is no loss of generality.
It is a classical result, following from the exact solution
of the model by Pfeuty \cite{pfeuty}, that in the
thermodynamic limit (\(L\to\infty\)) the quantum phase
transition occurs at \(g_c=1\), separating the
ferromagnetic phase (\(g<1\)) from the paramagnetic
phase (\(g>1\))
\cite{pfeuty,sachdev,dutta-tfim}.

For finite \(L\), this transition is rounded rather than
a genuine thermodynamic singularity.
One therefore studies finite-size estimators of the
critical point, often called pseudocritical points.
Their locations depend on the observable and on the
criterion used to define them.
Here we use the position \(g_{\mathrm{max},L}\) of the maximum
of the landscape velocity as such an estimator.
Our aim is to determine whether, and how rapidly,
\(g_{\mathrm{max},L}\) converges to the thermodynamic critical
point \(g_c=1\) as \(L\to\infty\).

This known value is the ground truth against which we evaluate every
estimator of the critical point in this section, ours and
Olsthoorn's.

The location of a peak alone does not establish that an observable
captures the critical behaviour: one should also examine how the
peak changes as the system grows. In this model,  the exact solution exhibits a nonanalyticity at the critical point \cite{pfeuty}, and the second derivative of the ground-state energy per site diverges logarithmically as $g\to1$ (Proposition~\ref{prop:app-log-infinite-energy}). Finite system size rounds this singularity
into a finite response, whose growth with \(L\) can be studied
through finite-size scaling. As we explain below, the landscape
velocity inherits this logarithmic behaviour in the periodic
chain, under the hypotheses stated in
Section~\ref{subsec:app-log-proof}. Thus logarithmic growth of its normalized
peak height provides evidence that the observable captures the
expected critical response, beyond merely exhibiting a maximum
near \(g=1\).

\paragraph{Ground state.}
The state is the ground state of the even-parity sector
\(P=\prod_j\sigma^z_j=+1\).
This is the \(\mathbb Z_2\)-symmetric state given by exact diagonalization at
finite \(L\), the one used by Olsthoorn, and symmetry-broken states are never
used.
The Jordan--Wigner transformation maps \(H\) to a quadratic Hamiltonian of free
fermions, with antiperiodic boundary conditions in the periodic case, so for
both boundary conditions the ground state in the sector \(P=+1\) is a
fermionic Gaussian state that can be computed exactly.
The one-site magnetization and the two-point functions
\(\langle\sigma^a_i\sigma^a_j\rangle\), \(a=x,y,z\), follow from its fermionic
two-point functions.
They determine the two-site reduced density matrices, whose spectrum is
available in closed form, and hence the correlations
\(C_{ij}=\langle\sigma^x_i\sigma^x_j\rangle\) and the von Neumann mutual
information \(M_{ij}\) exactly.
The cost is polynomial in \(L\), and for small \(L\) the results agree with
exact diagonalization of the even sector to machine precision.

\paragraph{Distances and normalization of bar lengths.}
Let \(\rho_i\) and \(\rho_{ij}\) be the one- and two-site reduced
density matrices, and let
\[
M_{ij}=S(\rho_i)+S(\rho_j)-S(\rho_{ij})
\]
be the von Neumann mutual information, computed with natural
logarithms.
For two qubits, \(0\le M_{ij}\le 2\ln 2\).

We also use the equal-time spin correlator
\[
C_{ij}
=
\langle\sigma_i^x\sigma_j^x\rangle
=
\operatorname{Tr}\!\left[
  \rho_{ij}(\sigma^x\otimes\sigma^x)
\right],
\qquad i\ne j,
\]
where the expectation is taken in the ground state at the
chosen value of \(g\).
This measures correlations between the spin components along
the interaction axis.
We use the raw correlator, without subtracting
\(\langle\sigma_i^x\rangle\langle\sigma_j^x\rangle\).

We use two dissimilarities on the vertex set
\(X=\{1,\dots,L\}\) of sites:
\begin{align}
D^{\mathrm{MI}}_{ij}
&=1-\frac{M_{ij}}{2\ln 2}\in[0,1],
\label{eq:app-dmi}\\
D^{C}_{ij}
&=1-C_{ij}\in[0,2],
\label{eq:app-dc}
\end{align}
for \(i\ne j\), with \(D_{ii}=0\).
The first one is the main distance: it is Olsthoorn's distance
\(2\ln2-M_{ij}\) \cite{olsthoorn-entanglement} divided by its
maximum value \(2\ln 2\), so it is measured in units of
\(2\ln2\), as on his axis ``\(\epsilon\,[2\ln 2]\)''.
The second dissimilarity is built from \(C_{ij}\)
and is not rescaled.
We use \(D^C\) for two reasons.
First, it gives a second, independent test of our method.
We apply to \(D^C\) exactly the same procedure as to
\(D^{\mathrm{MI}}\), without any change (\(H_0\) of the
Vietoris--Rips filtration, the landscape velocity \(\xi\)
and the position of its maximum), and we check whether
it still locates the critical point.
Second, \(D^C\) has a simpler measurement protocol than
\(D^{\mathrm{MI}}\).
All the correlators \(C_{ij}\) can be estimated from the same
collection of repeated measurements, with every qubit measured
in the \(x\) basis.
By contrast, evaluating the mutual information from
\(\rho_{ij}\) requires reconstructing the two-site reduced
density matrices, which generally involves measurements
in several bases.
This makes \(D^C\) an attractive choice for noisy
intermediate-scale quantum (NISQ) devices
\cite{preskill-nisq}.
In this line of research, we plan to compute \(\xi\)
from correlators measured on such devices.
These are the only normalizations of bar lengths in this section.
All bars are born at \(0\), so the length of each finite \(H_0\)
bar is its death time, and it is measured in the units of
\eqref{eq:app-dmi} or \eqref{eq:app-dc}.
We do not rescale the barcodes any further, for example by the
largest death time.
The two filtrations are always analysed separately.

The normalization does not affect the position of the maximum of
\(\xi\), but it does affect its height.
If \(D\) is multiplied by a constant \(c>0\), then so are all death
times and their derivatives, and Theorem~\ref{prop:closed}
shows that \(\xi\) is multiplied by \(c^{3/2}\).
This argmax invariance holds for multiplicative unit changes only; nonlinear reparameterizations such as $D\mapsto D^2$ or affine shifts change both the observable and the location of finite-size maximizers.
Hence the arg-max of \(\xi\) does not depend on the choice of
units, whereas the peak heights and the amplitudes \(A\) below
refer to the normalizations in \eqref{eq:app-dmi} and
\eqref{eq:app-dc}. Here \(A\) denotes the logarithmic amplitude,
that is, the coefficient multiplying the leading logarithmic
divergence of the per-bar landscape velocity.

\subsection{Topological observables}
\label{subsec:app-observables}

For each \(g\) we build the Vietoris--Rips filtration
\(\mathrm{VR}_t(X,D_g)\) of Section~\ref{subsec:3filtrations}, with
\(D_g=D^{\mathrm{MI}}\) or \(D_g=D^{C}\) evaluated in the ground
state of \(H(g)\), and compute its persistent homology with
\(\mathbb Z_2\) coefficients.
We retain the notation \(t\) for the filtration threshold
throughout, including in the figures; it corresponds to
\(\epsilon\) in \cite{olsthoorn-entanglement}.
For all the fields considered, both dissimilarities are positive off
the diagonal, so ~Assumption~\ref{ass:common-birth} holds.
By Section~\ref{subsec:rips-mst}, the \(L-1\) finite \(H_0\) death
times are the edge weights of a minimum spanning tree, which can be
computed with Prim's algorithm \cite{prim}; the infinite bar is
discarded.

In all our numerical data, for both dissimilarities and both
boundary conditions, the path \(1-2-\cdots-L\) is a minimum
spanning tree of the complete graph with edge weights
\(D_{ij}(g)\). Thus the tree can be chosen to use only edges
joining consecutive sites of the spin chain.
With periodic boundary conditions, this path omits the closing
edge \(\{L,1\}\) of the nearest-neighbour ring; by translation
invariance, omitting any other edge of the ring gives a minimum
spanning tree with the same total weight.
Consequently, the \(L-1\) finite \(H_0\) death times are the
nearest-neighbour distances \(D_{j,j+1}(g)\),
\(j=1,\dots,L-1\), arranged in decreasing order as
\(d_1(g)\ge\dots\ge d_{L-1}(g)\).
For the comparison with Olsthoorn, \(H_1\) is computed with GUDHI
\cite{maria-gudhi,gudhi-manual} on the \(2\)-skeleton, for
\(L\le128\).

\paragraph{Olsthoorn's topological summaries.}
\begin{itemize}
\item The \emph{\(H_0\) lifetime}, the mean of the \(L-1\) finite
  death times.
\item The \emph{\(H_1\) lifetime}, the length of the most persistent
  \(H_1\) bar. From his Fig.~6(d), we construct the estimator as the position of the argmax of the $H_1$ lifetime with periodic boundaries.
\item The \emph{range of \(H_0\) deaths},
  \(\Delta d=d_1-d_{L-1}\), the difference between the largest and the
  smallest finite death time, which he uses with open boundaries
  \cite[Appendix~E]{olsthoorn-entanglement}.
  His estimator of the transition is the position of the global maximum of
  \(\Delta d\), whose position we denote by
  \(g_{\mathrm{max},L}(\Delta d)\); this is the feature of
  his Fig.~12(b) that he relates to the transition.
  Since this maximum lies far from \(g=1\), 
  we introduce as an additional candidate to give his method more chances to succeed: the \emph{corner} of \(\Delta d\), that is, the
  non-differentiable local minimum that follows the maximum, which is closer to \(g=1\).
  We denote its position by \(g_{\mathrm{corner},L}\).
\end{itemize}
Thus Olsthoorn changes the homology degree and the summary with the boundary conditions: with periodic boundaries he
locates the critical point at the maximum of the \(H_1\) lifetime,
and with open boundaries at the maximum (or, in our second candidate, at the corner) of \(\Delta d\) for \(H_0\).
Our method makes the same choices in all cases: \(H_0\), the
landscape velocity \(\xi\), and the position of its maximum.
In this section the number of finite bars is always \(L-1\).

\paragraph{Landscape velocity and its normalization.}
The landscape velocity \(\xi(g)\) is that of
Definition~\ref{def:velocity}, without any normalization constant.
Since all bars are born at \(0\), Lemma~\ref{lem:nested} shows that
the \(k\)-th layer is the tent of the \(k\)-th largest death time,
and Theorem~\ref{prop:closed} gives
\begin{equation}
\label{eq:app-closed}
\xi(g)^2
=
\sum_{k=1}^{L-1}\int_{\mathbb R}\bigl|\partial_g\lambda_k(g;t)\bigr|^2\,dt
=
\frac12\sum_{k=1}^{L-1}d_k(g)\bigl(\partial_gd_k(g)\bigr)^2 .
\end{equation}
Thus no landscape needs to be computed.
To compare chains of different lengths $L$, we report the per-bar
quantity
\[
\frac{\xi(g)}{\sqrt{L-1}}
=
\Bigl(\frac{1}{L-1}\sum_{k=1}^{L-1}
\int_{\mathbb R}\bigl|\partial_g\lambda_k(g;t)\bigr|^2\,dt\Bigr)^{1/2},
\]
which is the root mean square of the velocities of the individual
layers, as announced after Definition~\ref{def:velocity} with
\(M=L-1\).
The square root is needed because \(\xi^2\), and not \(\xi\), is a sum
of \(L-1\) terms.

Since we aim to detect a change in the quantum phase, it is natural
to look for the parameter value at which the persistence landscape
changes most rapidly. The landscape velocity \(\xi(g)\) measures
this rate of change, motivating the estimator
\[
g_{\mathrm{max},L}\coloneq \operatorname*{arg\,max}_g \frac{\xi(g)}{\sqrt{L-1}}=\operatorname*{arg\,max}_g\xi(g).
\]
Note that these two amounts are the same because our normalization consists of simply dividing by a positive amount.
We then assess this estimator by studying whether \(g_{\mathrm{max},L}\)
converges to the critical point as \(L\to\infty\).
Throughout this section we distinguish two numbers.
The \emph{position of the maximum} is \(g_{\mathrm{max},L}\), the arg-max; it
is our estimator of the critical point.
The \emph{height of the maximum} is the maximum value
\(\max_g\xi(g)/\sqrt{L-1}\); it is not used to locate the critical
point, only to study how the maximum scales with \(L\).

For Olsthoorn's summaries, \(g_{\mathrm{max},L}\) likewise denotes the
position of the maximum of the summary under consideration (for \(\Delta d\) we also write
\(g_{\mathrm{max},L}(\Delta d)\)),
whereas \(g_{\mathrm{corner},L}\) denotes the position of the corner
of \(\Delta d\). These quantities depend on the chain length \(L\)
and are used as estimators of the critical field \(g_c=1\).

With periodic boundaries, translation invariance gives
\(d_k(g)=d_1(g)\) for every \(k\) (see the caption of
Figure~\ref{fig:pbc6}), and \eqref{eq:app-closed} reduces to
\begin{equation}
\label{eq:app-pbc}
\frac{\xi(g)}{\sqrt{L-1}}
=
\Bigl(\int_{\mathbb R}\bigl|\partial_g\lambda_1(g;t)\bigr|^2\,dt\Bigr)^{1/2}
=
\bigl|\partial_gd_1(g)\bigr|\sqrt{\frac{d_1(g)}{2}} .
\end{equation}
In that case Olsthoorn's \(H_0\) lifetime is exactly \(d_1\), and
\(\xi/\sqrt{L-1}\) is its derivative with respect to the field, up
to the bounded factor \(\sqrt{d_1/2}\).

\paragraph{Physical interpretation.}
For the correlator dissimilarity $D^C$ with periodic boundaries,
the normalized landscape velocity is directly related to
the second field derivative of the ground-state energy
per site. We derive this identity and its logarithmic
critical behaviour in
Section~\ref{subsec:app-log-proof}.

\subsection{Numerical method}
\label{subsec:app-numerics}

We use a logarithmic grid of \(401\) points in \([0.1,10]\)
(Olsthoorn's axis), and a linear fine grid in \([0.5,1.5]\)
with \(\Delta g=10^{-3}\), used to determine the positions
and heights of all maxima.
The derivatives \(\partial_gd_k\) are approximated by
second-order centred finite differences of the sorted
death times, without smoothing.
By Proposition~\ref{prop:regularity}, the sorted death
times are differentiable almost everywhere.
Each maximum is refined by fitting a parabola through
three neighbouring grid points.
With periodic boundaries we use \(L=10,14,18,22\)
(Olsthoorn's sizes) and \(L=32,64,128,256,512\).
With open boundaries we use \(L=12,16,20\)
(Olsthoorn's sizes) and \(L=32,64,128,256\).
The largest size, \(L=512\), is considered only with
periodic boundaries because, in our implementation,
the open-boundary computations require substantially
longer runtimes.

\subsection{Periodic boundary conditions}
\label{subsec:app-pbc}

\begin{figure}[htbp]
\centering
\includegraphics[width=\textwidth]{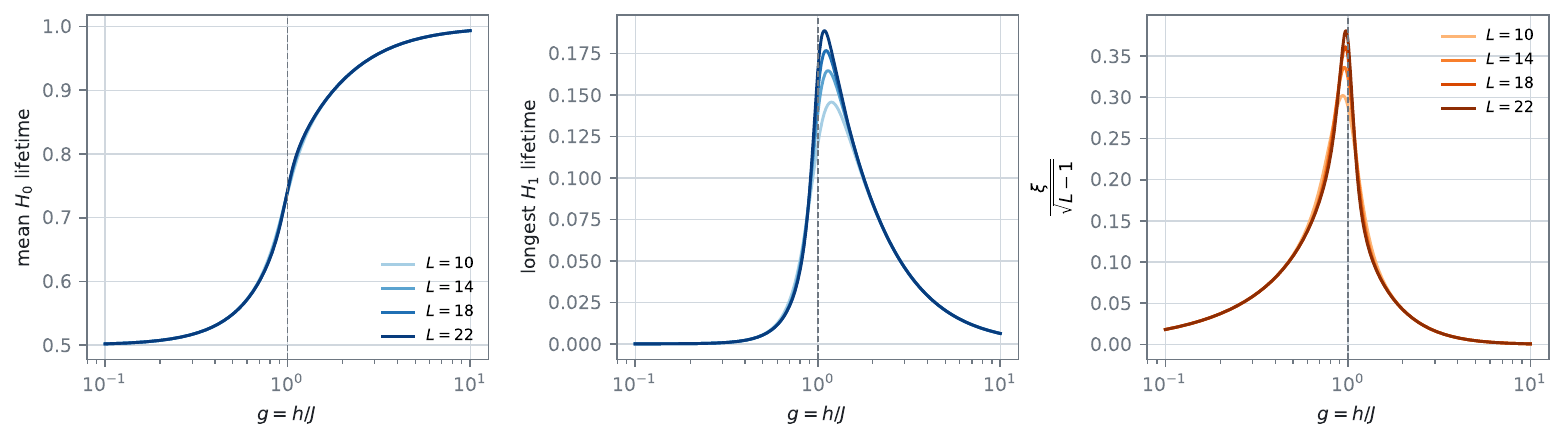}
\caption[Topological summaries of the periodic chain]{\textbf{Topological summaries of the periodic chain as functions of
the field.}
The left and centre panels, in \textcolor{olsthoornblue}{blue},
reproduce Olsthoorn's Fig.~6(d) \cite{olsthoorn-entanglement} and
show his summaries; the right panel, in
\textcolor{velocityorange}{orange}, shows our normalized landscape velocity.
All panels use the distance \(D^{\mathrm{MI}}\) and Olsthoorn's sizes
\(L=10,14,18,22\), with a logarithmic axis in \(g=h/J\); the dashed
line marks \(g_c=1\).
\textbf{Left} (Olsthoorn). The mean \(H_0\) lifetime increases
monotonically from \(0.50\) to \(1\); it shows no visible feature at
\(g_c\) and hardly depends on \(L\).
\textbf{Centre} (Olsthoorn). The length of the most persistent \(H_1\)
bar has a maximum whose height grows with \(L\) and whose position
lies in the paramagnet, at \(g=1.19\) for \(L=10\) and at \(g=1.09\)
for \(L=22\).
\textbf{Right} (ours). The per-bar landscape velocity
\(\xi/\sqrt{L-1}\) has a sharp maximum, whose position is close to
\(g_c\), at \(g=0.938\) for \(L=10\) and at \(g=0.974\) for \(L=22\),
and whose height grows with \(L\).
In Figure~\ref{fig:barcodes} all finite bars have the same length,
and this is no coincidence.
In the clean chain (without disorder) with periodic boundaries,
translation invariance makes all nearest-neighbour pairs equivalent.
Hence every bond has the same two-site state \(\rho_{j,j+1}\), and
therefore the same correlator
\(\langle\sigma^x_j\sigma^x_{j+1}\rangle\), which is minus the energy
of the bond in units of \(J\), and the same mutual information.
Thus all nearest-neighbour distances are equal to \(d_1\), and since
the minimum spanning tree is the nearest-neighbour path, all \(L-1\)
finite death times equal \(d_1\).
As a consequence, all layers \(\lambda_1,\dots,\lambda_{L-1}\) are the
same tent, and to compute \(\xi/\sqrt{L-1}\) we do not need to average
the contributions of the layers: we can simply take the contribution
of the first layer, which by Equation~\eqref{eq:app-pbc} is
\(\bigl(\int_{\mathbb R}|\partial_g\lambda_1(g;t)|^2\,dt\bigr)^{1/2}
=|\partial_gd_1|\sqrt{d_1/2}\).
\textbf{So the same \(H_0\) information that shows nothing in the left panel
locates the transition once it is differentiated with respect to the
field.}}
\label{fig:pbc6}
\end{figure}

\begin{figure}[htbp]
\centering
\includegraphics[width=\textwidth]{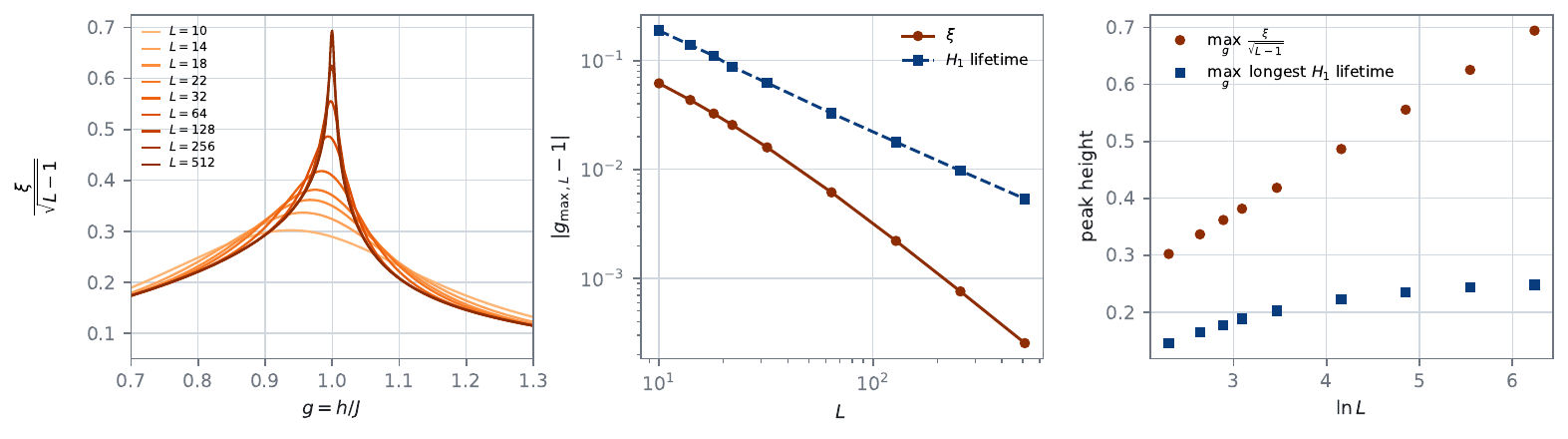}
\caption[Finite-size scaling for the periodic chain]{\textbf{Finite-size scaling of the landscape velocity for the
periodic chain.}
All panels use the distance \(D^{\mathrm{MI}}\) and the fine grid
\(\Delta g=10^{-3}\).
The left and right panels show only our landscape velocity, in
\textcolor{velocityorange}{orange}, and the centre panel compares it
(\textcolor{velocityorange}{orange} circles) with Olsthoorn's \(H_1\)
lifetime (\textcolor{olsthoornblue}{blue} squares).
\textbf{Left.} The per-bar velocity \(\xi/\sqrt{L-1}\) is plotted
against \(g\) for \(L=10\) to \(512\), from light to dark
\textcolor{velocityorange}{orange}.
\textbf{Centre.} The distance \(|g_{\mathrm{max},L}-1|\) between the critical
point and the position of the maximum is plotted against \(L\) on
logarithmic axes, for \(\xi\) and for the \(H_1\) lifetime
(\(L\le128\)).
For \(\xi\) it decays roughly as \(L^{-1.6}\) (a local estimate between
\(L=128\) and \(512\)), and for the \(H_1\) lifetime roughly as
\(L^{-1}\).
At \(L=128\) the two distances are \(2.2\times10^{-3}\) and
\(1.8\times10^{-2}\), respectively.
\textbf{Right.} The height of the maximum, \(\max_g\xi/\sqrt{L-1}\), is
plotted against \(\ln L\).
The points lie approximately along a straight line, suggesting that
the maximum diverges logarithmically with \(L\).}
\label{fig:pbcscaling}
\end{figure}

With periodic boundaries the \(H_0\) lifetime equals \(d_1\) and is monotone in
\(g\), so it shows no feature at the transition
(Figure~\ref{fig:pbc6}, left).
The normalized landscape velocity, in contrast, has a maximum near \(g=1\)
(Figure~\ref{fig:pbc6}, right).
The height of the maximum grows as \(\ln L\), and its position
converges to \(g_c=1\) faster than that of the \(H_1\) lifetime
(Figure~\ref{fig:pbcscaling}).
At \(L=128\) the position of the maximum of the \(H_1\) lifetime is \(8\) times
(\(D^{\mathrm{MI}}\)) and \(52\) times (\(D^C\)) farther from \(g_c\)
than the position of the maximum of \(\xi\).

\subsection{Open boundary conditions}
\label{subsec:app-obc}


\begin{figure}[htbp]
\centering
\includegraphics[width=0.7\textwidth]{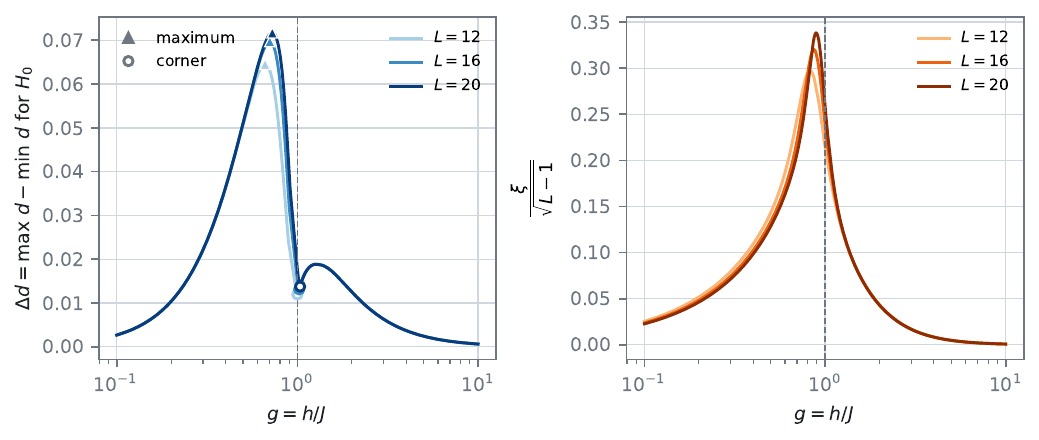}
\caption[Observables of the open chain]{\textbf{Observables of the open chain as functions of the field.}
The left panel, in \textcolor{olsthoornblue}{blue},
reproduce Olsthoorn's Fig.~12 \cite{olsthoorn-entanglement} and shows
his observable; the right panel, in
\textcolor{velocityorange}{orange}, shows our normalized landscape velocity.
All panels use the distance \(D^{\mathrm{MI}}\) and Olsthoorn's sizes
\(L=12,16,20\), with a logarithmic axis in \(g\); the dashed line
marks \(g_c=1\).
\textbf{Left} (Olsthoorn). The range \(\Delta d\) of the \(H_0\)
death times is his open-chain summary, and our curve reproduces his.
It has a main maximum at \(g\approx0.66\)--\(0.73\) (marked by triangles), which is the feature Olsthoorn looks at, then a non-differentiable minimum (marked by open circles) that we are referring to as a corner, and
finally a second maximum near \(g\approx1.3\).
The position of the maximum, \(g_{\mathrm{max},L}(\Delta d)\approx0.66\)--\(0.73\) for \(L=12\), \(16\) and \(20\), is Olsthoorn's estimator of
the critical point.
As a second candidate, which gives his method more chances to succeed because it is closer to \(g_c\), we also consider the position of the corner,
\(g_{\mathrm{corner},L}=0.998\), \(1.019\) and
\(1.029\) for \(L=12\), \(16\) and \(20\).
\textbf{Right} (ours). The per-bar velocity \(\xi/\sqrt{L-1}\) has a
single maximum, whose position approaches \(g_c\) as \(L\) grows:
it is at \(g=0.825\), \(0.867\) and \(0.893\).}
\label{fig:obc12}
\end{figure}

\begin{figure}[htbp]
\centering
\includegraphics[width=\textwidth]{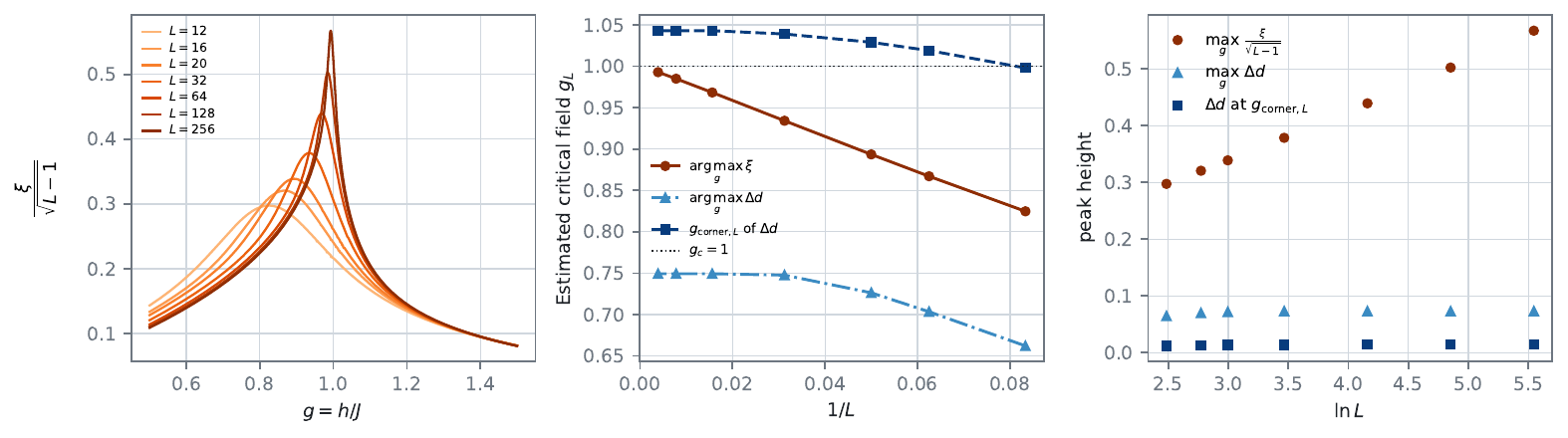}
\caption[Finite-size scaling for the open chain]{\textbf{Finite-size scaling of the landscape velocity for the open
chain.}
All panels use the distance \(D^{\mathrm{MI}}\) and the fine grid.
The left and right panels show only our landscape velocity, in
\textcolor{velocityorange}{orange}, and the centre panel compares the
position of its maximum (\textcolor{velocityorange}{orange} circles)
with Olsthoorn's two candidates (\textcolor{olsthoornblue}{blue}): the position of the maximum of
\(\Delta d\) (triangles) and the position of its corner (squares).
\textbf{Left.} The per-bar velocity \(\xi/\sqrt{L-1}\) is plotted
against \(g\) for \(L=12\) to \(256\), from light to dark
\textcolor{velocityorange}{orange}.
\textbf{Centre.} The estimated critical fields are plotted against \(1/L\):
the position \(g_{\mathrm{max},L}\) of the landscape-velocity maximum and
Olsthoorn's candidates (and our additional corner candidate): the position \(g_{\mathrm{max},L}(\Delta d)\) of the maximum of \(\Delta d\) and
the position \(g_{\mathrm{corner},L}\) of its corner.
The dotted line marks the critical field \(g_c=1\).
The position of the maximum of \(\xi\) approaches \(g_c=1\) roughly as \(1/L\), and it is \(0.9929\) at \(L=256\).
The maximum of \(\Delta d\) is far from \(g_c\) (\(g_{\mathrm{max},L}(\Delta d)\approx0.66\) at \(L=12\)) and stays near \(g\approx0.75\) for \(L\ge32\), so it does not approach the critical point.
The corner of \(\Delta d\) is close to \(g_c\) at \(L=12\)
(\(g_{\mathrm{corner},L}=0.998\)), but it moves away from \(g_c\) as \(L\)
grows and stays at \(g_{\mathrm{corner},L}=1.043\) for \(L\ge64\), so it
does not locate the transition in the thermodynamic limit.
\textbf{Right.} The height of the maximum, \(\max_g\xi/\sqrt{L-1}\),
is plotted against \(\ln L\). For \(D^{\mathrm{MI}}\), the points are
consistent with logarithmic growth; the local slopes are \(0.091\)
between \(L=64\) and \(128\) and \(0.094\) between \(128\) and \(256\).
For \(D^C\), the corresponding local slopes are \(0.127\) and \(0.132\),
close to the proved value \(A_C=0.1357\).}
\label{fig:obcscaling}
\end{figure}

With open boundaries the only Olsthoorn summary that responds to the
transition is \(\Delta d\).
Olsthoorn reports that \(\Delta d\) indicates the transition around
\(g=1\) \cite[Appendix~E]{olsthoorn-entanglement}.
The maximum of \(\Delta d\) lies at \(g\approx0.66\)--\(0.73\) at his
sizes (Figure~\ref{fig:obc12}, left), and it is the feature he relates to the transition.
This maximum is far from \(g_c=1\) and, as the chain grows, it does not approach it: it stays near
\(g\approx0.75\) for \(L\ge32\) (Figure~\ref{fig:obcscaling}, centre).
To give his method more chances to succeed, we also consider the corner that follows the maximum, which is closer to \(g=1\).
At his sizes the corner is indeed close to the critical point:
\(g_{\mathrm{corner},L}=0.998\), \(1.019\) and \(1.029\) for \(L=12\), \(16\)
and \(20\).
For longer chains, however, it moves away from \(g_c\), the opposite of
what an estimator should do, and freezes at \(g_{\mathrm{corner},L}=1.043\)
for \(L\ge64\) (Figure~\ref{fig:obcscaling}, centre).
So the agreement at \(L=12\) is a small-size coincidence.
Thus neither candidate locates the transition in the thermodynamic limit.
Since \(\Delta d\) measures the difference between two extreme death times, 
its features may be a property of the edge and not of the
transition.
Instead, the landscape velocity adds the responses of all bonds.
The position of its maximum converges to \(g_c\), and the height of the maximum grows as \(\ln L\),
as in the periodic chain.

\subsection{The correlator distance}
\label{subsec:app-correlator}

We now repeat the whole analysis with the correlator distance
\(D^C=1-\langle\sigma^x_i\sigma^x_j\rangle\)
.
As explained in Section~\ref{subsec:app-model}, this is a second test
of our method, so we apply exactly the same procedure as for
\(D^{\mathrm{MI}}\): \(H_0\) of the Vietoris--Rips filtration, the
landscape velocity \(\xi\), and the position of its maximum as the
estimator of the critical point.
Nothing has to be adapted to the new distance.

\begin{figure}[htbp]
\centering
\includegraphics[width=\textwidth]{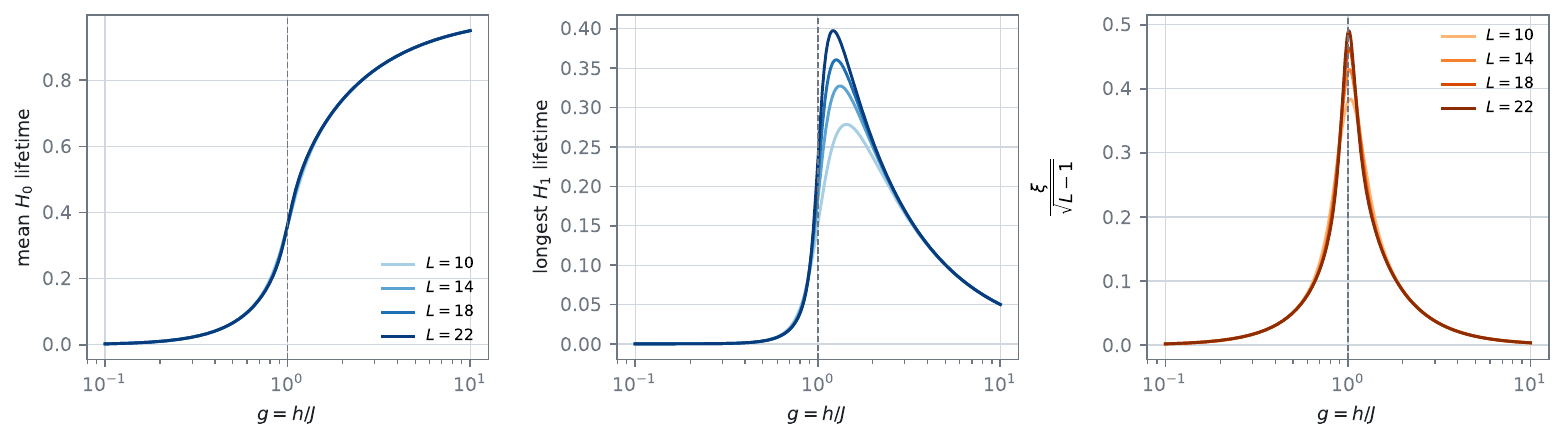}\\[2mm]
\includegraphics[width=\textwidth]{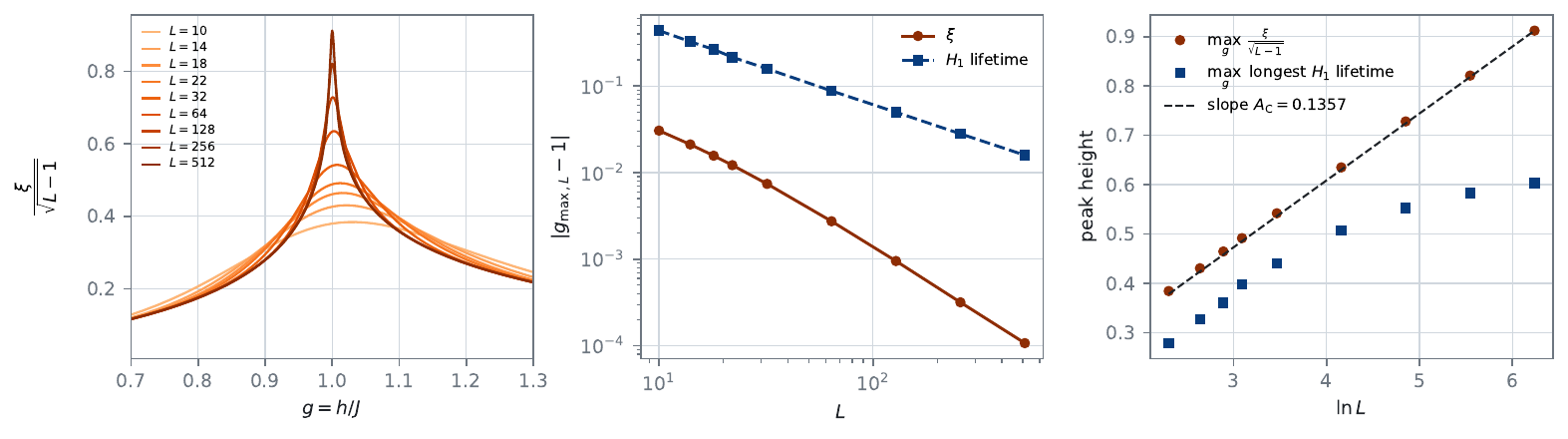}
\caption[The periodic chain with the correlator distance]{\textbf{The periodic chain with the correlator distance \(D^C\).}
The top row shows the panels of Figure~\ref{fig:pbc6} for \(D^C\):
Olsthoorn's summaries in \textcolor{olsthoornblue}{blue} (left and
centre) and our \(\xi/\sqrt{L-1}\) in
\textcolor{velocityorange}{orange} (right).
Here the \(H_0\) lifetime is
\(d_1=1-\langle\sigma^x_j\sigma^x_{j+1}\rangle\), which increases from
about \(0\) to about \(1\).
The bottom row shows the panels of Figure~\ref{fig:pbcscaling} for
\(D^C\), with our landscape velocity in
\textcolor{velocityorange}{orange} and Olsthoorn's \(H_1\) lifetime in
\textcolor{olsthoornblue}{blue} (centre).
At \(L=128\) the position of the maximum of the \(H_1\) lifetime is
\(g=1.050\), whereas the position of the maximum of \(\xi\) is
\(g=1.00095\).
The black dashed line on the right has the theoretical slope
\(A_C=\frac1\pi\sqrt{\frac12-\frac1\pi}\approx0.1357\) of
Section~\ref{subsec:app-log-proof}.
The fitted slope is \(0.1337\), so it differs from the proved
coefficient \(A_{\mathrm C}\approx0.1357\) by less than
\(1.5\%\).}
\label{fig:C1}
\end{figure}

With periodic boundaries the results are the same as with
\(D^{\mathrm{MI}}\).
The position of the maximum of \(\xi\) converges to \(g_c\), now from
above instead of from below, and faster than with \(D^{\mathrm{MI}}\):
at \(L=512\) it is \(g=1.00011\), whereas with \(D^{\mathrm{MI}}\) it is
\(g=0.99975\).
The position of the maximum of Olsthoorn's \(H_1\) lifetime also
converges to \(g_c\), but more slowly: at \(L=128\) it is \(g=1.050\),
whereas that of \(\xi\) is \(g=1.00095\) (Figure~\ref{fig:C1}).
The height of the maximum of \(\xi/\sqrt{L-1}\) again grows as
\(\ln L\), with a slope close to the theoretical value \(A_C\) of
Section~\ref{subsec:app-log-proof}.

\begin{figure}[htbp]
\centering
\includegraphics[width=0.7\textwidth]{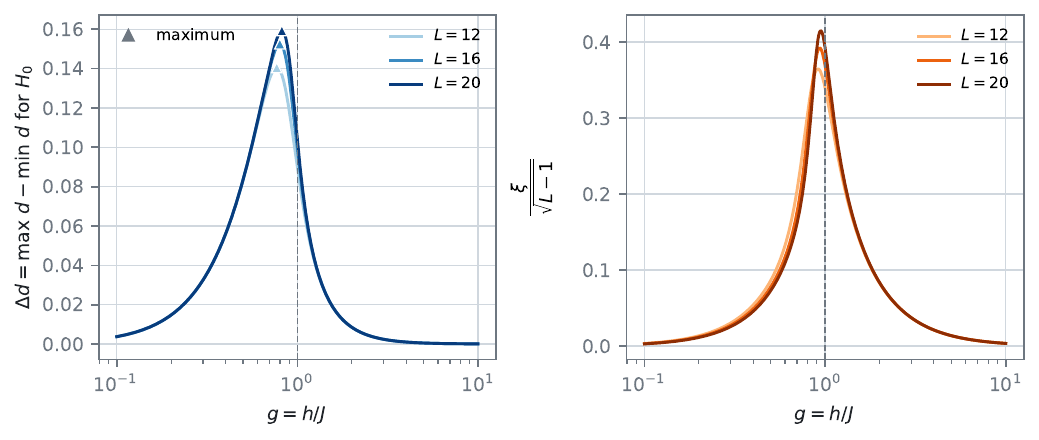}\\[2mm]
\includegraphics[width=\textwidth]{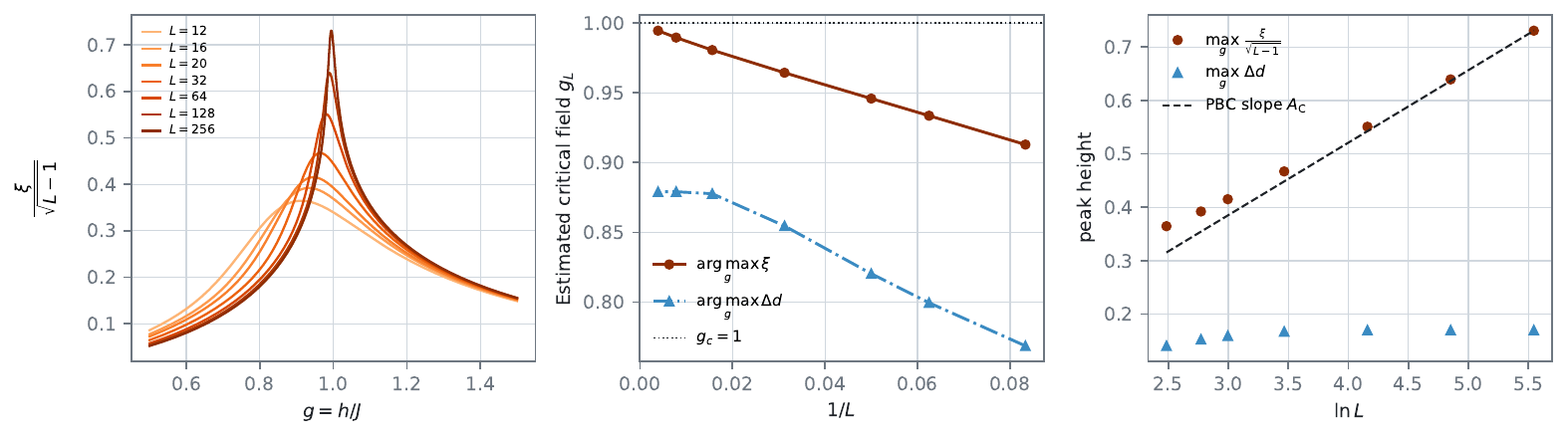}
\caption[The open chain with the correlator distance]{\textbf{The open chain with the correlator distance \(D^C\).}
The top row shows the panels of Figure~\ref{fig:obc12} for \(D^C\):
Olsthoorn's summary \(\Delta d\) in \textcolor{olsthoornblue}{blue} (left) and our \(\xi/\sqrt{L-1}\) in
\textcolor{velocityorange}{orange} (right).
The bottom row shows the panels of Figure~\ref{fig:obcscaling} for
\(D^C\); the left and right panels are ours, in \textcolor{velocityorange}{orange}.
Its centre panel compares the position of the maximum of \(\xi\) (\textcolor{velocityorange}{orange} circles) with the position of the maximum of \(\Delta d\)
(\textcolor{olsthoornblue}{blue} triangles), which is Olsthoorn's estimator.
For \(D^C\) the range \(\Delta d\) has no corner, so the second candidate used for \(D^{\mathrm{MI}}\) does not exist here.
The position of the maximum of \(\Delta d\) stays at \(g=0.8791\) for
\(L\ge128\), whereas the position of the maximum of \(\xi\) reaches
\(g=0.9945\) at \(L=256\).}
\label{fig:C2}
\end{figure}

With open boundaries the range \(\Delta d\) of the \(H_0\) death times
behaves differently for \(D^C\) than for \(D^{\mathrm{MI}}\): after
its maximum it decays smoothly and has no corner
(Figure~\ref{fig:C2}, top left).
Hence only one of Olsthoorn's two candidates is available for \(D^C\), the position of the maximum
of \(\Delta d\), which is the one he actually uses; the corner candidate, which was introduced for \(D^{\mathrm{MI}}\), does not exist.
The position of the maximum of \(\Delta d\) stays at \(g=0.8791\) for \(L\ge128\).
The position of the maximum of \(\xi\), in contrast, converges to
\(g_c\) as it does for \(D^{\mathrm{MI}}\), and it reaches \(g=0.9945\)
at \(L=256\) (Figure~\ref{fig:C2}, bottom centre).


\subsection{Advantages of the landscape velocity}
\label{subsec:app-advantages}

\begin{enumerate}
\item \textbf{One rule for both boundary conditions.}
    Olsthoorn adapts two choices to the boundary conditions: the
  homology degree (\(H_1\) with periodic boundaries, \(H_0\) with open
  ones) and the summary (the \(H_1\) lifetime with periodic boundaries, the range \(\Delta d\) with open
  ones).
  Our rule is always the same: \(H_0\) of the Vietoris--Rips
  filtration, the landscape velocity \(\xi\), and the position of its
  maximum, \(\operatorname*{arg\,max}_g\xi(g)\).
  
  This single rule locates the transition for both boundary
  conditions and for both distances (Figures~\ref{fig:pbc6}
  and~\ref{fig:obc12}, right panels; Figures~\ref{fig:C1}
  and~\ref{fig:C2}, top-right panels).
\item \textbf{Only \(H_0\) is needed.}
  Olsthoorn's periodic-chain summary requires \(H_1\), whereas \(\xi\)
  uses only \(H_0\), and \(H_0\) is much cheaper to compute.
  For a Vietoris--Rips filtration of \(L\) points, the \(H_0\) barcode is
  given by a minimum spanning tree (Section~\ref{subsec:rips-mst}),
  which Prim's algorithm computes in \(O(L^2)\) time from the dense
  distance matrix \cite{prim}, that is, in time linear in the size of
  the input.
  For \(H_1\) one needs at least the \(2\)-skeleton, which has
  \(\binom L3\sim L^3/6\) triangles, and the standard persistence
  algorithm is cubic in the number of simplices in the worst case
  \cite{elz,zomorodian-carlsson,edelsbrunner-harer}.
  Optimized software such as Ripser \cite{ripser} and GUDHI
  \cite{maria-gudhi} is much faster in practice, but it still has to
  process the triangles.
  Our computational cost includes the distance-matrix construction from input data and the numerical derivatives of those distances with respect to $g$. For PBC with translation-invariant correlators, the
two-site reduced density matrices need only be computed
once for each distinct site separation, rather than
separately for every pair of sites. The resulting
correlators and dissimilarities can then be reused for
all pairs with the same separation. This reduces the
number of distinct two-site calculations from
\(O(L^2)\) to \(O(L)\), independently of the subsequent
persistence computation.
  
  In our computations \(H_1\) was affordable only up to \(L=128\),
  while \(\xi\) was computed up to \(L=512\).
  Moreover, whenever the minimum spanning tree is the
  nearest-neighbour path, as in all our data, \(\xi\) needs only the
  \(L-1\) nearest-neighbour distances.
\item \textbf{No landscape is computed.}
  By Theorem~\ref{prop:closed}, \(\xi\) is an explicit function of
  the death times and their derivatives
  \eqref{eq:app-closed}; with periodic boundaries it reduces to a
  
  single function of \(d_1\) \eqref{eq:app-pbc}. 
\item \textbf{Faster convergence to the critical point.}
  The centre panels of
  Figures~\ref{fig:pbcscaling} and~\ref{fig:C1} (periodic case) and 
  the centre panels of Figures~\ref{fig:obcscaling} and~\ref{fig:C2} (open case) show that in all cases, our normalized velocity converges faster to the critical point. Indeed, with periodic boundaries, the distance \(|g_{\mathrm{max},L}-1|\)
  between \(g_c\) and the position of the maximum decays roughly as
  \(L^{-1.6}\) for \(\xi\), against roughly \(L^{-1}\) for the \(H_1\)
  lifetime (Figure~\ref{fig:pbcscaling}).
  With open boundaries, \(g_{\mathrm{max},L}\to1\) roughly as \(1/L\),
  whereas neither of Olsthoorn's candidates converges: the maximum of
  \(\Delta d\) stays near \(g\approx0.75\) for \(L\ge32\), and the
  corner of \(\Delta d\) moves away from \(g_c\), from \(0.998\) at
  \(L=12\) to \(1.043\) for \(L\ge64\) (Figure~\ref{fig:obcscaling}).
\item \textbf{A signature of criticality.}
  The normalized height of the maximum grows as \(\ln L\), with the slope
  derived in Section~\ref{subsec:app-log-proof}.
  This growth is consistent with the finite-size scaling
  associated with a logarithmic critical singularity in the
  thermodynamic limit; see, for example,
  \cite[Eqs.~(7)--(8)]{zhou2007fidelity}.
  The height-versus-\(\ln L\) plots are the right panels of
  Figures~\ref{fig:pbcscaling} and~\ref{fig:obcscaling}, and the
  lower-right panels of Figures~\ref{fig:C1} and~\ref{fig:C2}.
  The \(H_0\) lifetime, in contrast, has no feature at \(g_c\).
\item \textbf{Analytical detection of criticality.}
  For periodic boundaries with correlator dissimilarity, under
  Assumption~\ref{hyp:app-mst} and translation invariance, the proofs in
  Subsection~\ref{subsec:app-log-proof} establish that the peak
  location converges to \(g_c=1\) and its normalized height grows as
  \(A_{\mathrm C}\ln L\). Assumption~\ref{hyp:app-mst} is checked
  numerically for the parameter range studied; translation invariance
  follows from the symmetry of the periodic model. The numerical
  peak-location and peak-height trends are illustrated in the centre
  and lower-right panels, respectively, of Figure~\ref{fig:C1}.
\item \textbf{Robustness.}
  
  The conclusions hold for both distances \(D^{\mathrm{MI}}\) and
  \(D^C\), with exactly the same rule. For \(D^{\mathrm{MI}}\), this is
  illustrated in the right panels of Figures~\ref{fig:pbc6} and
  \ref{fig:obc12}, and the left and centre panels of
  Figures~\ref{fig:pbcscaling} and~\ref{fig:obcscaling}; the
  corresponding \(D^C\) results appear in the top-right and lower
  centre panels of Figures~\ref{fig:C1} and~\ref{fig:C2}.
\item \textbf{Ready for disorder.}
  Since \(\xi\) is defined on landscapes, the disorder averages of
  Section~\ref{subsec:disorder} apply directly, and
  Theorem~\ref{prop:mean-speed} computes the velocity of the
  
  empirical mean landscape from the individual death times only
  (Figures~\ref{fig:disorder-MI} and~\ref{fig:disorder-C}, panels
  (a)--(e)).
\end{enumerate}

\subsection{Analytical detection of the critical point by normalized landscape velocity}
\label{subsec:app-log-proof}

This section establishes the asymptotic behavior suggested by the
empirical data. In particular, Figure~\ref{fig:C1} (bottom row,
middle and right-hand panels) illustrates a trend in the normalized
landscape velocity for system sizes \(L \leq 512\); we show that
this trend persists as \(L \to \infty\).
Since the paper is primarily intended for a topological data analysis
audience, whereas the proofs require substantial background in physics,
we state only the results here. Details and proofs are
provided in Appendix~\ref{subsec:landscape-velocity-proofs}.

\begin{restatable}[Thermodynamic-limit landscape velocity]{corollary}{thermolimitcor}
\label{cor:app-log-infinite-velocity}
For the periodic chain in the even-parity sector with correlator
dissimilarity \(D^{\mathrm C}\), suppose Assumption~\ref{hyp:app-mst} holds
for all even \(L\) under consideration and all \(g\ne1\)
in a neighbourhood of \(1\).
Taking first the thermodynamic limit ($L\to\infty$) at fixed \(g\ne1\),
and then letting \(g\to1\), we obtain
\[
\lim_{L\to\infty}
\frac{\xi_L^{}(g)}{\sqrt{L-1}}
\sim
A_{\mathrm C}\ln\frac1{|g-1|}
\qquad(g\to1)
\]
where
\[
A_{\mathrm C}
=
\frac1\pi\sqrt{\frac12-\frac1\pi}
\approx0.135680.
\]
Here \(\sim\) denotes asymptotic equivalence: the ratio
of the two sides tends to \(1\) as \(g\to1\).
\end{restatable}

\begin{restatable}[Critical-field response]{proposition}{criticalfieldprop}
\label{prop:app-log-critical-response}
For the periodic chain in the even-parity sector, 
for correlator dissimilarity \(D^{\mathrm C}\), if
Assumption~\ref{hyp:app-mst} holds at \(g=1\) for all even \(L\)
under consideration, then as \(L\to\infty\)
through even integers,
\[
    \frac{\xi_L^{}(1)}{\sqrt{L-1}}
    =
    A_{\mathrm C}\ln L+O(1),
\]
where
\[
A_{\mathrm C}
=\frac1\pi\sqrt{\frac12-\frac1\pi}
\approx0.135680.
\]
\end{restatable}

\begin{restatable}[Normalized peak height]{proposition}{normalizedpeakprop}
\label{prop:app-log-peak}
For the periodic chain in the even-parity sector with correlator
dissimilarity \(D^{\mathrm C}\), let
\(I=[a,b]\subset(0,\infty)\), with \(a<1<b\).
Suppose Assumption~\ref{hyp:app-mst} holds for every
\(g\in I\) and every even \(L\) under consideration.
Then
\[
\max_{g\in I}
\frac{\xi_L^{}(g)}{\sqrt{L-1}}\sim A_{\mathrm C}\ln L,
\]
where
\[
A_{\mathrm C}
=\frac1\pi\sqrt{\frac12-\frac1\pi}
\approx0.135680.
\]
Every choice of maximizers
\(g_{\mathrm{max},L}\in I\) satisfies
\(g_{\mathrm{max},L}\to1\).
\end{restatable}

\subsection{Scope of the results and limitations}
\label{subsec:app-proved}

For the periodic chain, the proof in
Section~\ref{subsec:app-log-proof} uses two structural ingredients.
First, Assumption~\ref{hyp:app-mst} requires the nearest-neighbour path
to be a minimum spanning tree for the correlator dissimilarity; this
condition is checked numerically over the parameter range studied but
is not proved. Second, translation invariance makes all
nearest-neighbour dissimilarities equal. This is an exact consequence
of the periodic Hamiltonian and its translation-invariant ground-state
sector, not an additional hypothesis verified numerically.
For the open chain the convergence \(g_{\mathrm{max},L}\to1\) and the
\(\ln L\) growth are numerical results, for the following reasons.
\begin{enumerate}
\item Without translation invariance each bond has its own distance
  \(D_{j,j+1}(g)\), so one has to control \(L-1\) different functions
  instead of a single one.
\item The Bogoliubov modes of the open chain are not plane waves:
  their momenta are fixed by a transcendental equation that depends
  on \(g\) \cite{pfeuty}, there is an edge mode, and there is no
  explicit Riemann sum as in the periodic case.
\end{enumerate}
One route to a proof is to bound \(\xi(1)\) from below using only
the central bonds, with the exact solution of the open chain.

\paragraph{Data and code availability.}
The code used for the computations and figures is publicly
available on GitHub and archived on Zenodo
\cite{belchi-statuto-landscape-velocity-code}.






\section{Introducing randomness: landscape velocity under disorder averaging}
\label{subsec:disorder}

A \emph{disordered system} is a model in which some of
the parameters specifying the system are random,
according to a prescribed probability distribution \cite{binder-young}.
These parameters may describe, for example, spatial
variations in interaction strengths or local fields.
A \emph{disorder realization} is one particular choice
of all these random parameters. Once a realization
has been selected, it specifies a definite system.

For example, a transverse-field Ising chain with random
nearest-neighbour couplings can be described, in the same
spin convention and units as the application in Section \ref{sec:application}, by
\[
H^{(r)}(g)
=
-
\sum_{i}
J_i^{(r)}\,\sigma_i^x\sigma_{i+1}^x
-
h\sum_{i=1}^{L}\sigma_i^z.
\]
The vector
\[
\bigl(J_1^{(r)},\dots,J_{L_{\mathrm b}}^{(r)}\bigr)
\]
contains the couplings sampled for realization \(r\).
Setting \(J_i^{(r)}=1\) for every bond gives the clean
transverse-field Ising chain studied in the application. 

Each realization therefore gives its own
parameter-dependent state, correlation matrix,
dissimilarity matrix, and persistence landscape,
which we denote by \(\lambda^{(r)}(g)\).
Even at the same value of \(g\), these landscapes
can differ between realizations.
Disorder averaging combines the results from several
such realizations to describe behaviour across the
ensemble, rather than that of a single sampled system.
This averaging is distinct from taking a quantum
expectation value: the latter is computed within
the state of a fixed realization, whereas the former
averages quantities obtained from different choices
of the random model parameters.

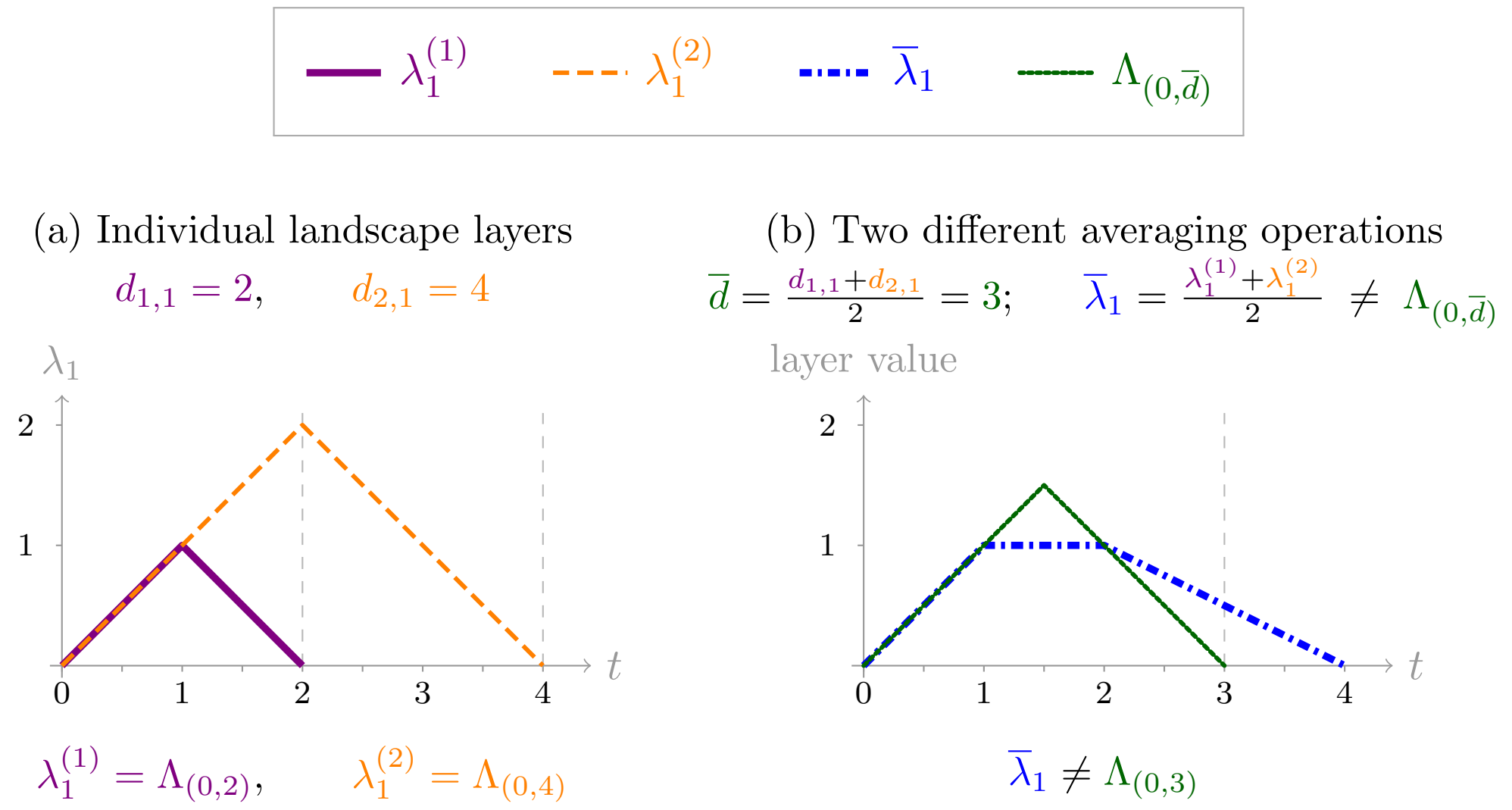
\begin{figure}[htbp]
\centering

\begingroup
\setlength{\fboxsep}{7pt}
\fcolorbox{gray!65}{white}{%
  \begin{tabular}{
    @{}l@{\hspace{7mm}}
       l@{\hspace{7mm}}
       l@{\hspace{7mm}}
       l@{}
  }

    \tikz[baseline=-0.6ex]{
      \draw[violet,line width=1.8pt,solid]
        (0,0) -- (0.65,0);
    }
    \({\color{violet}\lambda_1^{(1)}}\)

    &

    \tikz[baseline=-0.6ex]{
      \draw[
        orange,
        line width=1.1pt,
        dash pattern=on 4pt off 2pt
      ]
        (0,0) -- (0.65,0);
    }
    \({\color{orange}\lambda_1^{(2)}}\)

    &

    \tikz[baseline=-0.6ex]{
      \draw[
        blue,
        line width=1.8pt,
        dash pattern=on 3pt off 1.5pt on 1pt off 1.5pt
      ]
        (0,0) -- (0.65,0);
    }
    \({\color{blue}\overline\lambda_1}\)

    &

    \tikz[baseline=-0.6ex]{
      \draw[
        green!40!black,
        line width=1.1pt,
        densely dotted,
        line cap=round
      ]
        (0,0) -- (0.65,0);
    }
    \({\color{green!40!black}
       \Lambda_{(0,\overline d)}}\)

  \end{tabular}%
}
\endgroup

\par\vspace{5mm}

\makebox[\linewidth][c]{%
\begin{tikzpicture}[
  x=1.05cm,
  y=1.05cm,
  first landscape/.style={
    violet,
    line width=1.8pt,
    solid,
    line join=round
  },
  second landscape/.style={
    orange,
    line width=1.1pt,
    dash pattern=on 4pt off 2pt,
    line join=round
  },
  mean landscape/.style={
    blue,
    line width=1.8pt,
    dash pattern=on 3pt off 1.5pt on 1pt off 1.5pt,
    line join=round
  },
  averaged death tent/.style={
    green!40!black,
    line width=1.1pt,
    densely dotted,
    line cap=round,
    line join=round
  }
]

  %

  \begin{scope}

    \node[font=\small,align=center]
      at (2,3.60)
      {(a) Individual landscape layers};

    \node[font=\small]
      at (2,3.1)
      {$
        {\color{violet}d_{1,1}=2},
        \qquad
        {\color{orange}d_{2,1}=4}
      $};

    \draw[landscape axis]
      (-0.10,0) -- (4.40,0)
      node[right] {$t$};

    \draw[landscape axis]
      (0,-0.10) -- (0,2.25)
      node[above,font=\small] {$\lambda_1$};

    \draw[landscape guide]
      (2,0) -- (2,2.10);
    \draw[landscape guide]
      (4,0) -- (4,2.10);

    \foreach \x in {0,0.5,...,4}{
      \draw[gray!80]
        (\x,0) -- (\x,-0.05);
    }

    \foreach \x in {0,1,2,3,4}{
      \node[
        below=3pt,
        font=\scriptsize,
        fill=white,
        inner sep=1pt
      ] at (\x,0) {$\x$};
    }

    \foreach \y in {1,2}{
      \draw[gray!80]
        (0,\y) -- (-0.05,\y);
      \node[left=3pt,font=\scriptsize]
        at (0,\y) {$\y$};
    }

    \draw[first landscape]
      (0,0) -- (1,1) -- (2,0);

    \draw[second landscape]
      (0,0) -- (2,2) -- (4,0);

    \node[font=\small,align=center]
      at (2,-0.90)
      {$
        {\color{violet}
          \lambda_1^{(1)}=\Lambda_{(0,2)}},
        \qquad
        {\color{orange}
          \lambda_1^{(2)}=\Lambda_{(0,4)}}
      $};

  \end{scope}

  \begin{scope}[shift={(7cm,0)}]

    \node[font=\small,align=center]
      at (2,3.60)
      {(b) Two different averaging operations};

    \node[font=\small]
      at (2,3.1)
      {$
        {\color{green!40!black}\overline d}
        =
        \frac{
          {\color{violet}d_{1,1}}
          +
          {\color{orange}d_{2,1}}
        }{2}
        =
        {\color{green!40!black}3}
      $;\hspace{0.5cm}
      $
        {\color{blue}\overline\lambda_1}
        =
        \frac{
          {\color{violet}\lambda_1^{(1)}}
          +
          {\color{orange}\lambda_1^{(2)}}
        }{2}
        \;\ne\;
        {\color{green!40!black}
          \Lambda_{(0,\overline d)}}
      $
      };

    \draw[landscape axis]
      (-0.10,0) -- (4.40,0)
      node[right] {$t$};

    \draw[landscape axis]
      (0,-0.10) -- (0,2.25)
      node[above,font=\small] {layer value};

    \draw[landscape guide]
      (3,0) -- (3,2.10);

    \foreach \x in {0,0.5,...,4}{
      \draw[gray!80]
        (\x,0) -- (\x,-0.05);
    }

    \foreach \x in {0,1,2,3,4}{
      \node[
        below=3pt,
        font=\scriptsize,
        fill=white,
        inner sep=1pt
      ] at (\x,0) {$\x$};
    }

    \foreach \y in {1,2}{
      \draw[gray!80]
        (0,\y) -- (-0.05,\y);
      \node[left=3pt,font=\scriptsize]
        at (0,\y) {$\y$};
    }

    \draw[mean landscape]
      (0,0) -- (1,1) -- (2,1) -- (4,0);

    \draw[averaged death tent]
      (0,0) -- (1.5,1.5) -- (3,0);

    \node[font=\small,align=center]
      at (2,-0.90)
      {$
        {\color{blue}\overline\lambda_1}
        \ne
        {\color{green!40!black}
          \Lambda_{(0,3)}}
      $};

  \end{scope}

\end{tikzpicture}%
}
\caption{
Averaging landscapes is not the same as averaging
death times before constructing the landscape.
(a) Two realizations with a single finite bar each, with
\({\color{violet}d_{1,1}=2}\) and
\({\color{orange}d_{2,1}=4}\);
their landscape layers are solid \textcolor{violet}{violet} and dashed \textcolor{orange}{orange}.
(b) Their pointwise mean
\({\color{blue}\overline\lambda_1}\)
(dash-dotted \textcolor{blue}{blue}) has a flat segment on \([1,2]\).
By contrast, the averaged death
\({\color{green!40!black}\overline d=3}\)
produces the tent
\({\color{green!40!black}\Lambda_{(0,\overline d)}}\)
(dotted \textcolor{green!50!black}{dark green}), with a single peak at \(t=3/2\).
For example,
\({\color{blue}\overline\lambda_1(3/2)=1}\),
whereas
\({\color{green!40!black}
  \Lambda_{(0,\overline d)}(3/2)=3/2}\).
Thus the empirical mean need not itself be a
common-birth barcode landscape.
Both panels use the same scales, with equal
horizontal and vertical unit lengths.
}
\label{fig:mean-landscape-versus-mean-death}
\end{figure}

Fix \(R\geq1\) disorder realizations.
For realization \(r\), let \(\lambda^{(r)}(g)\) be
its landscape. 
Denote the empirical mean landscape
by
\[
\overline\lambda(g)
=
\frac1R\sum_{r=1}^R\lambda^{(r)}(g).
\]
The empirical mean landscape \(\overline\lambda(g)\)
belongs to \(\mathcal H\), but need not itself be the
persistence landscape of a common-birth barcode.
In particular, averaging the sorted death times and
substituting them into the single-barcode closed formula in Theorem \ref{prop:closed}
does not, in general, compute the velocity of the mean
landscape. However, in Theorem \ref{prop:mean-speed},
we prove a closed formula for the velocity of the
empirical mean landscape
that bypasses the need to compute it.
Figure~\ref{fig:mean-landscape-versus-mean-death}
illustrates that averaging the landscape layers and
constructing a landscape from the averaged death
times are different operations.

\subsection{Velocity of the empirical mean landscape versus mean of individual velocities}
\label{subsec:disorder-summaries}
Averaging the landscape before taking its velocity is
not the same as averaging the individual velocities.
We distinguish these two operations and compare them
with the root-mean-square individual velocity.

For the realization \(r\) landscape, \(\lambda^{(r)}(g)\), denote its  scalar landscape velocity as
\[
\xi_r(g)
=
\bigl\|\partial_g\lambda^{(r)}(g)\bigr\|_{\mathcal H}.
\]
When \(\lambda^{(r)}(g)\) is a landscape that comes from a common-birth barcode satisfying the hypotheses of
Theorem~\ref{prop:l2-speed}, this norm-based quantity
agrees with the landscape velocity in
Definition~\ref{def:velocity}.

The three scalar quantities to distinguish are
\begin{equation}
\label{eq:mean-speed}
\xi_{\overline\lambda}
:=
\left\|\partial_g\overline\lambda\right\|_{\mathcal H}
\end{equation}
\begin{equation}
\label{eq:mean-individual-speed}
\overline\xi
:=\frac1R\sum_{r=1}^R\xi_r
\end{equation}
\begin{equation}
\label{eq:rms-speed}    
\xi_{\mathrm{rms}}
:=\left(\frac1R\sum_{r=1}^R\xi_r^2\right)^{1/2}
\end{equation}

These are, respectively, the velocity of the
empirical mean landscape,
the mean of the individual velocities, and the root-mean-square
of the individual velocities.

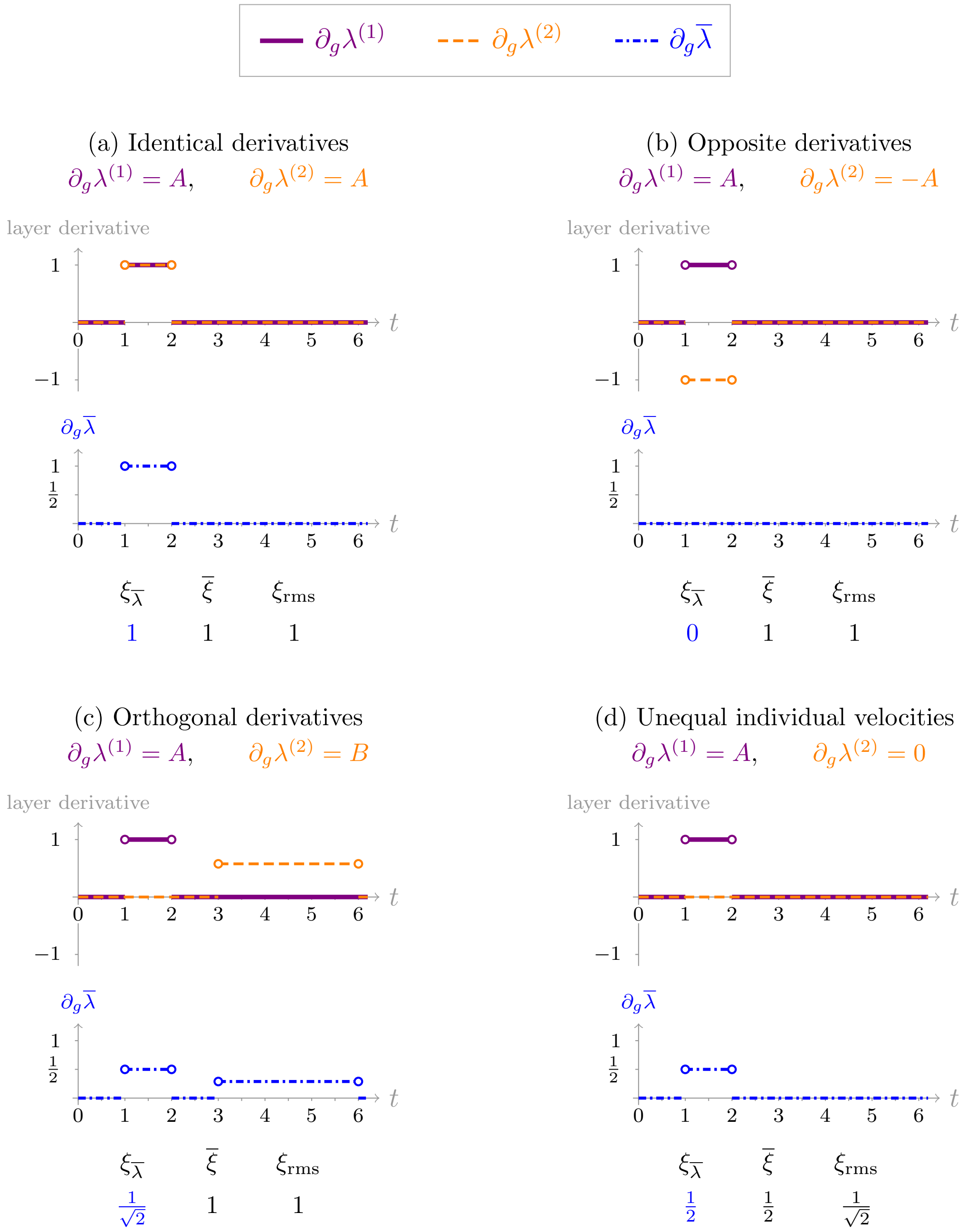
\begin{figure}[htbp]
\centering

\begingroup
\setlength{\fboxsep}{7pt}
\fcolorbox{gray!65}{white}{%
  \begin{tabular}{@{}l@{\hspace{7mm}}l@{\hspace{7mm}}l@{}}

    \tikz[baseline=-0.6ex]{
      \draw[violet,line width=1.8pt,solid]
        (0,0) -- (0.6,0);
    }
    \({\color{violet}\partial_g\lambda^{(1)}}\)

    &

    \tikz[baseline=-0.6ex]{
      \draw[
        orange,
        line width=1.1pt,
        dash pattern=on 4pt off 2pt
      ]
        (0,0) -- (0.6,0);
    }
    \({\color{orange}\partial_g\lambda^{(2)}}\)

    &

    \tikz[baseline=-0.6ex]{
      \draw[
        blue,
        line width=1.2pt,
        dash pattern=on 3pt off 1.5pt on 1pt off 1.5pt
      ]
        (0,0) -- (0.6,0);
    }
    \({\color{blue}\partial_g\overline\lambda}\)

  \end{tabular}%
}
\endgroup

\par\vspace{6mm}

\makebox[\linewidth][c]{%
\begin{tikzpicture}[
  x=0.65cm,
  y=0.80cm,
  first derivative/.style={
    violet,
    line width=1.8pt,
    solid,
    line join=round
  },
  second derivative/.style={
    orange,
    line width=1.1pt,
    dash pattern=on 4pt off 2pt,
    line join=round
  },
  mean derivative/.style={
    blue,
    line width=1.2pt,
    dash pattern=on 3pt off 1.5pt on 1pt off 1.5pt,
    line join=round
  }
]

  %
  %
  %

\foreach \xshift/\yshift/\mode/\heading in {
  0/0/1/{(a) Identical derivatives},
  7.8/0/2/{(b) Opposite derivatives},
  0/-8/3/{(c) Orthogonal derivatives},
  7.8/-8/4/{(d) Unequal individual velocities}
}{
    \begin{scope}[shift={(\xshift cm,\yshift cm)}]

      \ifcase\mode
      \or
        \def\secondexpression{A}
        \def\meanspeed{1}
        \def\averagespeed{1}
        \def\rmsspeed{1}
        \pgfmathsetmacro{\secondstart}{1}
        \pgfmathsetmacro{\secondend}{2}
        \pgfmathsetmacro{\secondheight}{1}
      \or
        \def\secondexpression{-A}
        \def\meanspeed{0}
        \def\averagespeed{1}
        \def\rmsspeed{1}
        \pgfmathsetmacro{\secondstart}{1}
        \pgfmathsetmacro{\secondend}{2}
        \pgfmathsetmacro{\secondheight}{-1}
      \or
        \def\secondexpression{B}
        \def\meanspeed{\frac1{\sqrt2}}
        \def\averagespeed{1}
        \def\rmsspeed{1}
        \pgfmathsetmacro{\secondstart}{3}
        \pgfmathsetmacro{\secondend}{6}
        \pgfmathsetmacro{\secondheight}{1/sqrt(3)}
      \or
        \def\secondexpression{0}
        \def\meanspeed{\frac12}
        \def\averagespeed{\frac12}
        \def\rmsspeed{\frac1{\sqrt2}}
        \pgfmathsetmacro{\secondstart}{1}
        \pgfmathsetmacro{\secondend}{2}
        \pgfmathsetmacro{\secondheight}{0}
      \fi

    \node[font=\small,align=center]
      at (3,3.10) {\heading};
    
    \node[font=\small]
      at (3,2.50)
      {$
        {\color{violet}\partial_g\lambda^{(1)}=A},
        \qquad
        {\color{orange}
          \partial_g\lambda^{(2)}=\secondexpression}
      $};

      \draw[landscape axis]
        (-0.12,0) -- (6.45,0)
        node[right] {$t$};

      \draw[landscape axis]
        (0,-1.20) -- (0,1.30)
        node[above,font=\scriptsize]
        {layer derivative};

      \foreach \x in {0,0.5,...,6}{
        \draw[gray!80]
          (\x,0) -- (\x,-0.05);
      }
      \foreach \x in {0,1,...,6}{
        \node[
          below=3pt,
          font=\scriptsize,
          fill=white,
          inner sep=1pt
        ] at (\x,0) {$\x$};
      }

      \foreach \y in {-1,1}{
        \draw[gray!80] (0,\y) -- (-0.07,\y);
        \node[left=3pt,font=\scriptsize]
          at (0,\y) {$\y$};
      }

      \draw[first derivative]
        (0,0) -- (1,0);
      \draw[first derivative]
        (1,1) -- (2,1);
      \draw[first derivative]
        (2,0) -- (6.20,0);

      \ifnum\mode=4

        \draw[second derivative]
          (0,0) -- (6.20,0);

      \else

        \draw[second derivative]
          (0,0) -- (\secondstart,0);
        \draw[second derivative]
          (\secondstart,\secondheight)
          -- (\secondend,\secondheight);
        \draw[second derivative]
          (\secondend,0) -- (6.20,0);

      \fi

      \foreach \x in {1,2}{
        \draw[draw=violet,fill=white,thick]
          (\x,1) circle[radius=1.5pt];
      }

      \ifnum\mode=4
      \else
        \draw[draw=orange,fill=white,thick]
          (\secondstart,\secondheight)
          circle[radius=1.5pt];
        \draw[draw=orange,fill=white,thick]
          (\secondend,\secondheight)
          circle[radius=1.5pt];
      \fi

      \begin{scope}[shift={(0,-2.8cm)}]

        \draw[landscape axis]
          (-0.12,0) -- (6.45,0)
          node[right] {$t$};

        \draw[landscape axis]
          (0,-0.12) -- (0,1.30)
          node[above,font=\scriptsize,text=blue]
          {$\partial_g\overline\lambda$};

        \foreach \x in {0,0.5,...,6}{
          \draw[gray!80]
            (\x,0) -- (\x,-0.05);
        }
        \foreach \x in {0,1,...,6}{
          \node[
            below=3pt,
            font=\scriptsize,
            fill=white,
            inner sep=1pt
          ] at (\x,0) {$\x$};
        }

        \foreach \y/\lab in {
          0.5/{\frac12},
          1/1
        }{
          \draw[gray!80]
            (0,\y) -- (-0.07,\y);
          \node[left=3pt,font=\scriptsize]
            at (0,\y) {$\lab$};
        }

        \ifnum\mode=2

          \draw[mean derivative]
            (0,0) -- (6.20,0);

        \else

          \ifnum\mode=1
            \pgfmathsetmacro{\meanheight}{1}
          \else
            \pgfmathsetmacro{\meanheight}{0.5}
          \fi

          \draw[mean derivative]
            (0,0) -- (1,0);
          \draw[mean derivative]
            (1,\meanheight) -- (2,\meanheight);

          \foreach \x in {1,2}{
            \draw[draw=blue,fill=white,thick]
              (\x,\meanheight)
              circle[radius=1.5pt];
          }

          \ifnum\mode=3

            \draw[mean derivative]
              (2,0) -- (3,0);
            \draw[mean derivative]
              (3,{1/(2*sqrt(3))})
              -- (6,{1/(2*sqrt(3))});
            \draw[mean derivative]
              (6,0) -- (6.20,0);

            \foreach \x in {3,6}{
              \draw[draw=blue,fill=white,thick]
                (\x,{1/(2*sqrt(3))})
                circle[radius=1.5pt];
            }

          \else

            \draw[mean derivative]
              (2,0) -- (6.20,0);

          \fi
        \fi

      \end{scope}

      \node[font=\small,align=center]
        at (3,-5.05)
        {%
          \begin{tabular}{@{}c@{\hspace{8mm}}c@{\hspace{8mm}}c@{}}
            \(\xi_{\overline\lambda}\)
            &
            \(\overline\xi\)
            &
            \(\xi_{\mathrm{rms}}\)
            \\[5pt]
            \({\color{blue}\meanspeed}\)
            &
            \(\averagespeed\)
            &
            \(\rmsspeed\)
          \end{tabular}%
        };

    \end{scope}
  }

\end{tikzpicture}%
}
\caption{
Comparison of the three scalar velocity summaries for
two realizations.
Only the first landscape layer is nonzero.
Write \(A(t)=\mathbf1_{(1,2)}(t)\) and
\(B(t)=3^{-1/2}\mathbf1_{(3,6)}(t)\), both of
\(L^2\) norm \(1\).
Each panel shows the individual derivatives above
(solid \textcolor{violet}{violet} for
\({\color{violet}\partial_g\lambda^{(1)}}\),
dashed \textcolor{orange}{orange} for
\({\color{orange}\partial_g\lambda^{(2)}}\))
and their mean below (dash-dotted \textcolor{blue}{blue}).
The three velocity summaries are listed underneath.
(a) Identical derivatives reinforce;
the \textcolor{orange}{orange} curve lies on the thicker \textcolor{violet}{violet} curve.
(b) Opposite derivatives cancel in the mean.
(c) Disjoint supports give orthogonality:
\(\langle A,B\rangle_{L^2}
=\int_{\mathbb R}A(t)B(t)\,dt=0\).
(d) The
\textcolor{orange}{second landscape is constant
with respect to \(g\) near the parameter considered},
so \({\color{orange}\partial_g\lambda^{(2)}=0}\).
Open circles omit endpoint values, which do not
affect the norms.
Corresponding graphs across panels use the same scales.
}
\label{fig:disorder-speed-comparison}
\end{figure}

Figure~\ref{fig:disorder-speed-comparison} compares
the three scalar summaries in four configurations
of two landscape derivatives, illustrating the effects
of reinforcement, cancellation, orthogonality, and
unequal individual velocities.

At a parameter value where all individual landscape
curves are differentiable in \(\mathcal H\), finite
summation commutes with differentiation:
\[
\partial_g\overline\lambda
=
\frac1R\sum_{r=1}^R\partial_g\lambda^{(r)}.
\]
Taking norms does not, in general, commute with averaging.
Indeed, the triangle inequality and the
Cauchy--Schwarz inequality give
\[
\underbrace{
\left\|
\frac1R\sum_{r=1}^R\partial_g\lambda^{(r)}
\right\|_{\mathcal H}
}_{\displaystyle \xi_{\overline\lambda}}
\;\leq\;
\underbrace{
\frac1{R}
\sum_{r=1}^R
\|\partial_g\lambda^{(r)}\|_{\mathcal H}
}_{\displaystyle \overline\xi}
\;\leq\;
\underbrace{
\left(
\frac1{R}\sum_{r=1}^R
\|\partial_g\lambda^{(r)}\|_{\mathcal H}^2
\right)^{1/2}
}_{\displaystyle \xi_{\mathrm{rms}}}.
\]

The difference between the squared root-mean-square
velocity and the squared velocity of the
empirical mean landscape
can be computed exactly:
\[
\begin{aligned}
\frac1R\sum_{r=1}^R
\|\partial_g\lambda^{(r)}-\partial_g\overline\lambda\|_{\mathcal H}^2
&=
\frac1R\sum_{r=1}^R
\left(
\|\partial_g\lambda^{(r)}\|_{\mathcal H}^2
-2\langle\partial_g\lambda^{(r)},
          \partial_g\overline\lambda\rangle_{\mathcal H}
+\|\partial_g\overline\lambda\|_{\mathcal H}^2
\right)
\\
&=
\frac1R\sum_{r=1}^R
\|\partial_g\lambda^{(r)}\|_{\mathcal H}^2
-2\left\langle
\frac1R\sum_{r=1}^R\partial_g\lambda^{(r)},
\partial_g\overline\lambda
\right\rangle_{\mathcal H}
+\|\partial_g\overline\lambda\|_{\mathcal H}^2
\\
&=
\frac1R\sum_{r=1}^R
\|\partial_g\lambda^{(r)}\|_{\mathcal H}^2
-2\left\langle
\partial_g\overline\lambda,
\partial_g\overline\lambda
\right\rangle_{\mathcal H}
+\|\partial_g\overline\lambda\|_{\mathcal H}^2
\\
&=
\frac1R\sum_{r=1}^R
\|\partial_g\lambda^{(r)}\|_{\mathcal H}^2
-2\|\partial_g\overline\lambda\|_{\mathcal H}^2
+\|\partial_g\overline\lambda\|_{\mathcal H}^2
\\
&=
\frac1R\sum_{r=1}^R
\|\partial_g\lambda^{(r)}\|_{\mathcal H}^2
-\|\partial_g\overline\lambda\|_{\mathcal H}^2.
\end{aligned}
\]
Hence
\begin{equation}\label{eq:velocity-variance}
\xi_{\mathrm{rms}}^2-\xi_{\overline\lambda}^2
=
\frac1{R}\sum_{r=1}^R
\|\partial_g\lambda^{(r)}-\partial_g\overline\lambda\|_{\mathcal H}^2.
\end{equation}
In particular,
\[
\xi_{\mathrm{rms}}=\xi_{\overline\lambda}
\quad\Longleftrightarrow\quad
\partial_g\lambda^{(r)}=\partial_g\overline\lambda
\text{ in }\mathcal H
\quad\text{for every }r=1,\dots,R.
\]

\paragraph{Statistical interpretation.}
The Hilbert-space framework also makes the derivative fields
amenable to standard statistical analysis.
At a fixed parameter value $g$, let
$V_i=\partial_g\lambda_i(g)$ be independent and identically
distributed measurable $H$-valued random variables.
If $\mathbb E\|V_i\|_H<\infty$, the strong law of large numbers
implies that the velocity of the empirical mean landscape
converges almost surely to $\|\mathbb E V_i\|_H$.
When differentiation and expectation may be interchanged,
this limit is the velocity of the population mean landscape.
If $\mathbb E\|V_i\|_H^2<\infty$, a Hilbert-space central limit
theorem also applies to the empirical mean derivative.


\subsection{\texorpdfstring{A closed formula for the velocity of the empirical mean landscape}{A closed formula for the velocity of the empirical mean landscape}}
\label{subsec:mean-velocity}

We now develop a closed formula, analogous to that in Theorem \ref{prop:closed}, for the velocity of the
empirical mean landscape.
Indeed, for \(H_0\) Vietoris--Rips, \v{C}ech, alpha, and other filtrations in this setting, the formula in Theorem~\ref{prop:mean-speed} computes the velocity \(\xi_{\overline\lambda}\) of the empirical mean landscape from the individual barcode death times, without constructing the individual landscapes.
Thus neither the individual landscapes nor their empirical mean need to be computed. Given $R$ independent samples with $M$ deaths each, the direct evaluation of Equation~\eqref{eq:mean-overlap} costs $O(R^2M)$ in additions and multiplications, whereas naive alternatives such as first constructing the mean landscape (requiring representation-dependent costs) may be more expensive. The empirical mean landscape is unique, but in general is not the landscape of any barcode, so Theorem~\ref{prop:closed} cannot be applied to it directly.

\begin{theorem}[Closed formula for the velocity of
the empirical mean landscape]
\label{prop:mean-speed}
Let \(I\subseteq\mathbb R\) be an open interval.
Suppose that, for each realization \(r=1,\dots,R\)
and every \(g\in I\), the finite part of the barcode is
\[
\mathcal B_g^{(r),\mathrm{fin}}
=
\bigl\{[0,d_{r,k}(g))\bigr\}_{k=1}^{M},
\qquad
d_{r,1}(g)\geq\cdots\geq d_{r,M}(g)>0,
\]
where \(R\geq1\) and \(M\geq1\) are fixed.
Let \(\lambda^{(r)}(g)\in\mathcal H\) be the landscape
constructed from these finite bars, with any infinite
bars omitted, and define
\[
\overline\lambda(g)
=
\frac1R\sum_{r=1}^R\lambda^{(r)}(g).
\]

Assume that every sorted death function \(d_{r,k}\)
is differentiable at \(g_0\in I\).
Then the individual landscapes and their empirical mean
are differentiable at \(g_0\) in \(\mathcal H\).

In the formulas below, all deaths, their derivatives,
and landscape derivatives are evaluated at \(g_0\),
and \([x]_+=\max\{x,0\}\).
The squared velocity of the
empirical mean landscape is
\begin{equation}\label{eq:mean-overlap}
\boxed{
\xi_{\overline\lambda}^2
=\frac1{R^2}\sum_{r=1}^R\sum_{s=1}^R\sum_{k=1}^M
\bigl(\partial_g d_{r,k}\bigr)
\bigl(\partial_g d_{s,k}\bigr)
\left[\min\{d_{r,k},d_{s,k}\}
-\frac12\max\{d_{r,k},d_{s,k}\}\right]_+.
}
\end{equation}
where
\[
\xi_{\overline\lambda}
=
\|\partial_g\overline\lambda\|_{\mathcal H}
\]
In contrast, the squared root-mean-square individual velocity is
\begin{equation}\label{eq:rms-closed}
\boxed{
\xi_{\mathrm{rms}}^2
=
\frac1{2R}
\sum_{r=1}^R\sum_{k=1}^M
d_{r,k}\bigl(\partial_g d_{r,k}\bigr)^2.
}
\end{equation}
\end{theorem}

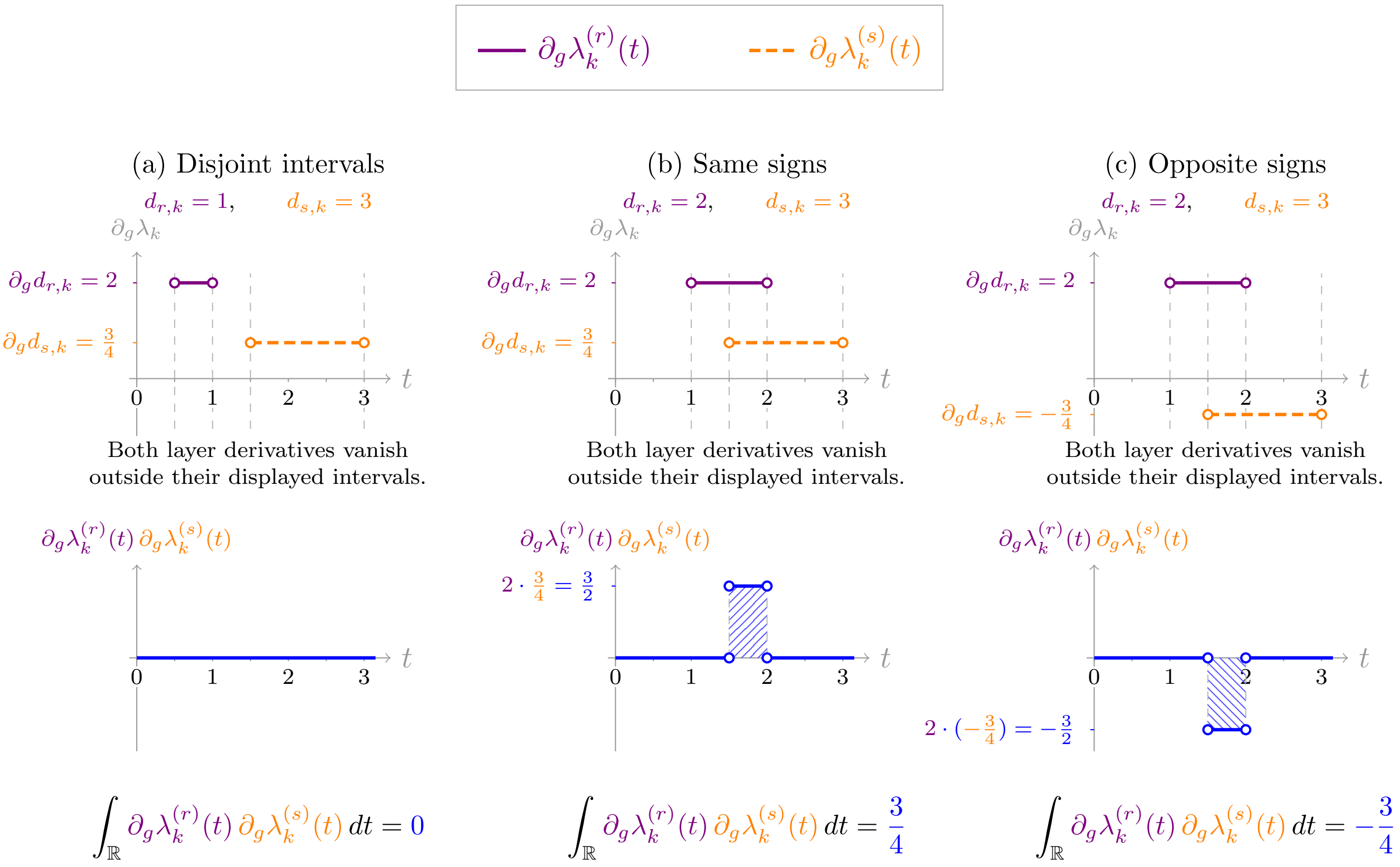
\begin{figure}[htbp]
\centering

\begingroup
\setlength{\fboxsep}{7pt}
\fcolorbox{gray!65}{white}{%
  \begin{tabular}{@{}l@{\hspace{12mm}}l@{}}
    \tikz[baseline=-0.6ex]{
      \draw[landscape curve,violet,solid]
        (0,0) -- (0.6,0);
    }
    \({\color{violet}\partial_g\lambda_k^{(r)}(t)}\)
    &
    \tikz[baseline=-0.6ex]{
      \draw[
        landscape curve,
        orange,
        dash pattern=on 4pt off 2pt
      ]
        (0,0) -- (0.6,0);
    }
    \({\color{orange}\partial_g\lambda_k^{(s)}(t)}\)
  \end{tabular}%
}
\endgroup

\par\vspace{6mm}

\makebox[\linewidth][c]{%
\begin{tikzpicture}[x=0.95cm,y=0.60cm]

  \foreach \xshift/\firstdeath/\seconddeath/\sgn/\heading/\result in {
    0/1/3/1/{(a) Disjoint intervals}/{0},
    6/2/3/1/{(b) Same signs}/{\frac34},
    12/2/3/-1/{(c) Opposite signs}/{-\frac34}
  }{

    \begin{scope}[shift={(\xshift cm,0)}]

      \pgfmathsetmacro{\secondheight}{3*\sgn/4}
      \pgfmathsetmacro{\productheight}{3*\sgn/2}
      \pgfmathsetmacro{\leftend}
        {max(\firstdeath/2,\seconddeath/2)}
      \pgfmathsetmacro{\rightend}
        {min(\firstdeath,\seconddeath)}


      \node[font=\small,align=center]
        at (1.60,4.45) {\heading};

      \node[font=\scriptsize]
        at (1.60,3.65)
        {$
          {\color{violet}d_{r,k}=\firstdeath},
          \qquad
          {\color{orange}d_{s,k}=\seconddeath}
        $};

      \draw[landscape axis]
        (-0.10,0) -- (3.35,0)
        node[right] {$t$};

      \draw[landscape axis]
        (0,-1.20) -- (0,2.65)
        node[above,font=\scriptsize]
        {$\partial_g\lambda_k$};

      \draw[landscape guide]
        ({\firstdeath/2},-1.05)
        -- ({\firstdeath/2},2.20);
      \draw[landscape guide]
        (\firstdeath,-1.05)
        -- (\firstdeath,2.20);
      \draw[landscape guide]
        ({\seconddeath/2},-1.05)
        -- ({\seconddeath/2},2.20);
      \draw[landscape guide]
        (\seconddeath,-1.05)
        -- (\seconddeath,2.20);

      \foreach \x in {0,0.5,1,1.5,2,2.5,3}{
        \draw[gray!80]
          (\x,0) -- (\x,-0.07);
      }

      \foreach \x in {0,1,2,3}{
        \node[
          below=3pt,
          font=\scriptsize,
          fill=white,
          inner sep=1pt
        ] at (\x,0) {$\x$};
      }

      \draw[violet] (0,2) -- (-0.05,2);
      \node[
        left=3pt,
        font=\scriptsize,
        text=violet
      ] at (0,2)
        {$\partial_g d_{r,k}=2$};

      \draw[orange]
        (0,\secondheight)
        -- (-0.05,\secondheight);

      \ifnum\sgn>0
        \node[
          left=3pt,
          font=\scriptsize,
          text=orange
        ] at (0,\secondheight)
          {$\partial_g d_{s,k}=\frac34$};
      \else
        \node[
          left=3pt,
          font=\scriptsize,
          text=orange
        ] at (0,\secondheight)
          {$\partial_g d_{s,k}=-\frac34$};
      \fi

      \draw[landscape curve,violet,solid]
        ({\firstdeath/2},2)
        -- (\firstdeath,2);

      \draw[
        landscape curve,
        orange,
        dash pattern=on 4pt off 2pt
      ]
        ({\seconddeath/2},\secondheight)
        -- (\seconddeath,\secondheight);

      \draw[draw=violet,fill=white,thick]
        ({\firstdeath/2},2)
        circle[radius=1.6pt];
      \draw[draw=violet,fill=white,thick]
        (\firstdeath,2)
        circle[radius=1.6pt];

      \draw[draw=orange,fill=white,thick]
        ({\seconddeath/2},\secondheight)
        circle[radius=1.6pt];
      \draw[draw=orange,fill=white,thick]
        (\seconddeath,\secondheight)
        circle[radius=1.6pt];

      \node[font=\scriptsize,align=center]
        at (1.60,-1.80)
        {Both layer derivatives vanish\\
         outside their displayed intervals.};

      \begin{scope}[shift={(0,-3.5cm)}]

        \draw[landscape axis]
          (-0.10,0) -- (3.35,0)
          node[right] {$t$};

        \draw[landscape axis]
          (0,-1.95) -- (0,1.95)
          node[above,font=\scriptsize]
          {$
            {\color{violet}\partial_g\lambda_k^{(r)}(t)}
            \,
            {\color{orange}\partial_g\lambda_k^{(s)}(t)}
          $};

        \foreach \x in {0,0.5,1,1.5,2,2.5,3}{
          \draw[gray!80]
            (\x,0) -- (\x,-0.07);
        }

        \foreach \x in {0,1,2,3}{
          \node[
            below=3pt,
            font=\scriptsize,
            fill=white,
            inner sep=1pt
          ] at (\x,0) {$\x$};
        }

        \ifdim\leftend pt<\rightend pt

          \ifnum\sgn>0
            \fill[
              pattern=north east lines,
              pattern color=blue!65
            ]
              (\leftend,0)
              rectangle
              (\rightend,\productheight);
          \else
            \fill[
              pattern=north west lines,
              pattern color=blue!65
            ]
              (\leftend,0)
              rectangle
              (\rightend,\productheight);
          \fi

          \draw[landscape guide]
            (\leftend,0)
            -- (\leftend,\productheight);
          \draw[landscape guide]
            (\rightend,0)
            -- (\rightend,\productheight);

          \draw[landscape curve,blue]
            (0,0) -- (\leftend,0);
          \draw[landscape curve,blue]
            (\leftend,\productheight)
            -- (\rightend,\productheight);
          \draw[landscape curve,blue]
            (\rightend,0) -- (3.15,0);

          \foreach \x in {\leftend,\rightend}{
            \draw[draw=blue,fill=white,thick]
              (\x,0) circle[radius=1.6pt];
            \draw[draw=blue,fill=white,thick]
              (\x,\productheight)
              circle[radius=1.6pt];
          }

          \draw[blue]
            (0,\productheight)
            -- (-0.05,\productheight);

          \ifnum\sgn>0
            \node[
              left=3pt,
              font=\scriptsize,
              text=blue
            ] at (0,\productheight)
              {$
                {\color{violet}2}
                \cdot
                {\color{orange}\frac34}
                =\frac32
              $};
          \else
            \node[
              left=3pt,
              font=\scriptsize,
              text=blue
            ] at (0,\productheight)
              {$
                {\color{violet}2}
                \cdot
                ({\color{orange}-\frac34})
                =-\frac32
              $};
          \fi

        \else

          \draw[landscape curve,blue]
            (0,0) -- (3.15,0);

        \fi

        \node[font=\small]
          at (1.60,-3.55)
          {$\displaystyle
            \int_{\mathbb R}
            {\color{violet}\partial_g\lambda_k^{(r)}(t)}
            \,
            {\color{orange}\partial_g\lambda_k^{(s)}(t)}
            \,dt
            ={\color{blue}\result}
          $};

      \end{scope}

    \end{scope}
  }

\end{tikzpicture}%
}
\caption{
Pairwise contributions to the derivative inner product
in \eqref{eq:velocity-overlap}, at a fixed parameter
value and for a fixed layer \(k\).
The upper row shows the two landscape-layer derivatives:
solid \textcolor{violet}{violet} for realization \textcolor{violet}{\(r\)},
dashed \textcolor{orange}{orange} for realization \textcolor{orange}{\(s\)},
as identified in the boxed legend.
Their nonzero heights are the death-time derivatives
\({\color{violet}\partial_g d_{r,k}}\) and
\({\color{orange}\partial_g d_{s,k}}\),
as indicated on the vertical axes.
Each derivative is nonzero on its falling-side interval,
meaning the interval of filtration values
\({\color{violet}(d_{r,k}/2,d_{r,k})}\) or
\({\color{orange}(d_{s,k}/2,d_{s,k})}\),
and vanishes elsewhere, apart from the omitted
breakpoint values.
In (a), the deaths are
\({\color{violet}d_{r,k}=1}\) and
\({\color{orange}d_{s,k}=3}\);
in (b)--(c), they are
\({\color{violet}d_{r,k}=2}\) and
\({\color{orange}d_{s,k}=3}\).
In every column,
\({\color{violet}\partial_g d_{r,k}=2}\);
the second death-time derivative is
\({\color{orange}\partial_g d_{s,k}=3/4}\)
in (a)--(b) and
\({\color{orange}\partial_g d_{s,k}=-3/4}\)
in (c).
The lower row shows the product of the layer derivatives
in \textcolor{blue}{blue}.
In (a), the falling-side intervals
\({\color{violet}(1/2,1)}\) and
\({\color{orange}(3/2,3)}\) are disjoint,
so the contribution is \({\color{blue}0}\).
In (b)--(c), the intervals
\({\color{violet}(1,2)}\) and
\({\color{orange}(3/2,3)}\)
intersect in \((3/2,2)\), of length \(1/2\).
Multiplying this length by the product heights
\({\color{violet}2}\cdot
 {\color{orange}(3/4)}={\color{blue}3/2}\)
and
\({\color{violet}2}\cdot
 {\color{orange}(-3/4)}={\color{blue}-3/2}\)
gives the contributions
\({\color{blue}3/4}\) and
\({\color{blue}-3/4}\), respectively.
Positive product regions are filled with diagonal lines
sloping upwards from left to right; negative product
regions are filled with diagonal lines sloping downwards
from left to right.
Open circles omit breakpoint values, which do not
affect the integrals.
All columns use the same coordinate scales.
}
\label{fig:derivative-overlaps}
\end{figure}

\begin{proof}
The common-birth property gives
\[
\lambda_k^{(r)}(g;t)
=
\max\{0,\min\{t,d_{r,k}(g)-t\}\},
\qquad 1\leq k\leq M,
\]
and all higher layers vanish.
The argument of Theorem~\ref{prop:l2-speed},
applied to these \(M\) finite common-birth bars,
shows that each landscape is differentiable at \(g_0\)
in \(\mathcal H\).
Its derivative is represented almost everywhere by
\[
\partial_g\lambda_k^{(r)}(t)
=
\bigl(\partial_g d_{r,k}\bigr)
\mathbf1_{\left(\frac{d_{r,k}}2,d_{r,k}\right)}(t),
\qquad 1\leq k\leq M.
\]
Since the empirical mean is a finite sum,
it is also differentiable, with
\[
\partial_g\overline\lambda
=
\frac1R\sum_{r=1}^R\partial_g\lambda^{(r)}.
\]

We first compute the inner product of two individual
landscape derivatives.
By the definition of the inner product in
\eqref{eq:hilbert},
\[
\begin{aligned}
\left\langle
\partial_g\lambda^{(r)},\partial_g\lambda^{(s)}
\right\rangle_{\mathcal H}
&=
\sum_{k=1}^{M}
\int_{\mathbb R}
\partial_g\lambda_k^{(r)}(t)\,
\partial_g\lambda_k^{(s)}(t)\,dt
\\
&=
\sum_{k=1}^{M}
\bigl(\partial_g d_{r,k}\bigr)
\bigl(\partial_g d_{s,k}\bigr)
\int_{\mathbb R}
\mathbf1_{\left(\frac{d_{r,k}}2,d_{r,k}\right)}(t)\,
\mathbf1_{\left(\frac{d_{s,k}}2,d_{s,k}\right)}(t)\,dt.
\end{aligned}
\]
The death-time derivatives have been taken outside
the integral because they do not depend on \(t\).

The product of the indicator functions is the indicator
of the intersection
\[
\left(\frac{d_{r,k}}2,d_{r,k}\right)
\cap
\left(\frac{d_{s,k}}2,d_{s,k}\right).
\]
If nonempty, this intersection has left endpoint
\[
\max\left\{
\frac{d_{r,k}}2,\frac{d_{s,k}}2
\right\}
=
\frac12\max\{d_{r,k},d_{s,k}\}
\]
and right endpoint
\[
\min\{d_{r,k},d_{s,k}\}.
\]
Its length is therefore
\[
\left[
\min\{d_{r,k},d_{s,k}\}
-\frac12\max\{d_{r,k},d_{s,k}\}
\right]_+,
\]
where the use of \([x]_+=\max\{x,0\}\) also accounts for an empty
intersection or intervals meeting only at an endpoint.
Thus, the integral is
the length of the intersection:
\[
\int_{\mathbb R}
\mathbf1_{\left(\frac{d_{r,k}}2,d_{r,k}\right)}(t)\,
\mathbf1_{\left(\frac{d_{s,k}}2,d_{s,k}\right)}(t)\,dt
=
\left[
\min\{d_{r,k},d_{s,k}\}
-\frac12\max\{d_{r,k},d_{s,k}\}
\right]_+.
\]
Consequently,
\begin{equation}\label{eq:velocity-overlap}
\left\langle
\partial_g\lambda^{(r)},\partial_g\lambda^{(s)}
\right\rangle_{\mathcal H}
=
\sum_{k=1}^{M}
\bigl(\partial_g d_{r,k}\bigr)
\bigl(\partial_g d_{s,k}\bigr)
\left[
\min\{d_{r,k},d_{s,k}\}
-\frac12\max\{d_{r,k},d_{s,k}\}
\right]_+.
\end{equation}
Figure~\ref{fig:derivative-overlaps} illustrates how
the overlap of the falling-side intervals
\((d_{r,k}/2,d_{r,k})\) and
\((d_{s,k}/2,d_{s,k})\), together with the product
of the death-time derivatives, determines each
layer's contribution to
\eqref{eq:velocity-overlap}.
Only corresponding layers are paired, as required
by the inner product in \(\mathcal H\).

We now expand the squared velocity of the
empirical mean landscape:
\[
\begin{aligned}
\xi_{\overline\lambda}^{\,2}
&=
\|\partial_g\overline\lambda\|_{\mathcal H}^2
\\
&=
\left\langle
\frac1R\sum_{r=1}^R\partial_g\lambda^{(r)},
\frac1R\sum_{s=1}^R\partial_g\lambda^{(s)}
\right\rangle_{\mathcal H}
\\
&=
\frac1{R^2}
\sum_{r=1}^R\sum_{s=1}^R
\left\langle
\partial_g\lambda^{(r)},\partial_g\lambda^{(s)}
\right\rangle_{\mathcal H}.
\end{aligned}
\]
Substituting \eqref{eq:velocity-overlap} yields
\eqref{eq:mean-overlap}.

Finally, taking \(s=r\) in \eqref{eq:velocity-overlap}
gives
\[
\begin{aligned}
\|\partial_g\lambda^{(r)}\|_{\mathcal H}^2
&=
\sum_{k=1}^M
\bigl(\partial_g d_{r,k}\bigr)^2
\left[d_{r,k}-\frac12d_{r,k}\right]_+
\\
&=
\frac12\sum_{k=1}^M
d_{r,k}\bigl(\partial_g d_{r,k}\bigr)^2,
\end{aligned}
\]
since \(d_{r,k}>0\).
Using the definition of the root-mean-square velocity,
we obtain
\[
\begin{aligned}
\xi_{\mathrm{rms}}^2
&=
\frac1R\sum_{r=1}^R\xi_r^2
=
\frac1{R}\sum_{r=1}^R
\|\partial_g\lambda^{(r)}\|_{\mathcal H}^2
\\
&=
\frac1{2R}
\sum_{r=1}^R\sum_{k=1}^M
d_{r,k}\bigl(\partial_g d_{r,k}\bigr)^2,
\end{aligned}
\]
which proves \eqref{eq:rms-closed}.
\end{proof}

\paragraph{Recovery of the single-barcode formula.}
When \(R=1\), there is only one realization, so
\(\overline\lambda=\lambda^{(1)}\).
For every \(k\), the overlap factor in
\eqref{eq:mean-overlap} becomes
\[
\left[
\min\{d_{1,k},d_{1,k}\}
-\frac12\max\{d_{1,k},d_{1,k}\}
\right]_+
=
\frac{d_{1,k}}2,
\]
because \(d_{1,k}>0\).
Consequently,
\[
\xi_{\overline\lambda}^{\,2}
=
\frac1{2}\sum_{k=1}^M
d_{1,k}\bigl(\partial_g d_{1,k}\bigr)^2.
\]
This is the single-barcode closed formula; in the setting of
Theorem~\ref{prop:closed}, one takes
\(M=L-1\).

\paragraph{Diagonal and cross contributions.}
For a general number \(R\) of realizations, separating
the terms with \(r=s\) from those with \(r\ne s\) gives
\[
\begin{aligned}
\xi_{\overline\lambda}^{\,2}
&=
\frac1{R^2}
\sum_{r=1}^R
\|\partial_g\lambda^{(r)}\|_{\mathcal H}^2
+
\frac1{R^2}
\sum_{\substack{r,s=1\\r\ne s}}^R
\left\langle
\partial_g\lambda^{(r)},\partial_g\lambda^{(s)}
\right\rangle_{\mathcal H}
\\
&=
\frac{\xi_{\mathrm{rms}}^2}{R}
+
\frac1{R^2}
\sum_{\substack{r,s=1\\r\ne s}}^R
\left\langle
\partial_g\lambda^{(r)},\partial_g\lambda^{(s)}
\right\rangle_{\mathcal H}.
\end{aligned}
\]
Thus \(\xi_{\mathrm{rms}}^2/R\) is exactly the
diagonal contribution to the squared velocity of
the
empirical mean landscape.
The difference between the two quantities is
determined entirely by the cross terms:
\[
\xi_{\overline\lambda}^{\,2}
-\frac{\xi_{\mathrm{rms}}^2}{R}
=
\frac1{R^2}
\sum_{\substack{r,s=1\\r\ne s}}^R
\left\langle
\partial_g\lambda^{(r)},\partial_g\lambda^{(s)}
\right\rangle_{\mathcal H}.
\]

In particular, suppose that, for some small
\(\varepsilon>0\),
\[
\left|
\sum_{\substack{r,s=1\\r\ne s}}^R
\left\langle
\partial_g\lambda^{(r)},\partial_g\lambda^{(s)}
\right\rangle_{\mathcal H}
\right|
\leq
\varepsilon
\sum_{r=1}^R
\|\partial_g\lambda^{(r)}\|_{\mathcal H}^2.
\]
Then
\[
\left|
\xi_{\overline\lambda}^{\,2}
-\frac{\xi_{\mathrm{rms}}^2}{R}
\right|
\leq
\varepsilon\frac{\xi_{\mathrm{rms}}^2}{R},
\]
which makes precise the approximation
\[
\xi_{\overline\lambda}^{\,2}
\approx
\frac{\xi_{\mathrm{rms}}^2}{R}.
\]
This condition concerns the total signed sum of
cross terms, not merely the size of each individual
inner product.
There are \(R(R-1)\) ordered cross terms, so many
small contributions can accumulate when \(R\)
is large.

A positive cross term
\[
\left\langle
\partial_g\lambda^{(r)},\partial_g\lambda^{(s)}
\right\rangle_{\mathcal H}>0
\]
increases the squared velocity of the
empirical mean landscape
relative to the diagonal contribution, whereas
a negative cross term decreases it.
Thus derivative fields can reinforce or partially
cancel one another.
If every pair of distinct derivative fields is
orthogonal, all cross terms vanish and
\[
\xi_{\overline\lambda}^{\,2}
=
\frac{\xi_{\mathrm{rms}}^2}{R}
\]
holds exactly.

\paragraph{Independence does not imply orthogonality.}
The preceding identities concern a fixed collection
of realizations.
If the realizations are instead regarded as random,
their independence does not imply that their
derivative fields are orthogonal, either for a
particular sample or in expectation.

More precisely, suppose that, at the parameter value
under consideration, the derivative fields
\(\partial_g\lambda^{(1)},\dots,\partial_g\lambda^{(R)}\)
are independent and identically distributed as
\(\mathcal H\)-valued random elements \cite{ledoux-talagrand}, with
\[
\mathbb E
\|\partial_g\lambda^{(1)}\|_{\mathcal H}^2<\infty.
\]
For \(r\ne s\), independence gives
\[
\begin{aligned}
\mathbb E\left[
\left\langle
\partial_g\lambda^{(r)},\partial_g\lambda^{(s)}
\right\rangle_{\mathcal H}
\right]
&=
\left\langle
\mathbb E[\partial_g\lambda^{(r)}],
\mathbb E[\partial_g\lambda^{(s)}]
\right\rangle_{\mathcal H}
\\
&=
\left\|
\mathbb E[\partial_g\lambda^{(1)}]
\right\|_{\mathcal H}^2.
\end{aligned}
\]
Hence these expected cross terms vanish only if
the derivative field has zero mean:
\[
\mathbb E[\partial_g\lambda^{(1)}]=0
\quad\text{in }\mathcal H.
\]
Indeed, taking expectations in the diagonal--cross
decomposition yields
\[
\mathbb E[\xi_{\overline\lambda}^{\,2}]
=
\frac1{R}
\mathbb E\|\partial_g\lambda^{(1)}\|_{\mathcal H}^2
+
\frac{R-1}{R}
\left\|
\mathbb E[\partial_g\lambda^{(1)}]
\right\|_{\mathcal H}^2.
\]
Under the additional zero-mean assumption, this reduces to
\[
\mathbb E[\xi_{\overline\lambda}^{\,2}]
=
\frac1R\,\mathbb E[\xi_{\mathrm{rms}}^2].
\]
This is an identity in expectation; it does not
assert that the empirical cross terms vanish
for each collection of realizations.

\subsection{Numerical illustration with random couplings}
\label{subsec:disorder-numerical}

We illustrate the three disorder summaries
of Section~\ref{subsec:disorder-summaries} on the random-bond open
transverse-field Ising chain, a disordered version of the model in
Section~\ref{subsec:app-model}:
\[
H^{(r)}(g)
=
-\sum_{i=1}^{L-1}J_i^{(r)}\sigma_i^x\sigma_{i+1}^x
-g\sum_{i=1}^{L}\sigma_i^z,
\qquad
J_i^{(r)}=e^{u_i^{(r)}},
\qquad
u_i^{(r)}\overset{\mathrm{iid}}{\sim}\mathrm{Unif}[-1,1].
\]
Fisher's study of the random transverse-field Ising chain
describes it as the simplest nontrivial interacting quantum model
with quenched randomness \cite{fisher-prl}. Here we use a
bond-disordered version with independent logarithms
\(u_i^{(r)}=\ln J_i^{(r)}\) uniformly distributed on \([-1,1]\),
while keeping the transverse field uniform across sites.
Here \(J_i^{(r)}\) is the coupling on bond \((i,i+1)\), meaning
the interaction strength multiplying
\(\sigma_i^x\sigma_{i+1}^x\). The exponential parametrization makes
each coupling positive; equivalently, its logarithm \(u_i^{(r)}\)
is uniform on \([-1,1]\), so the couplings are log-uniform on
\([e^{-1},e]\).
The field \(g\) is uniform across the sites. Since
\(\mathbb E[\ln J_i]=0\), Fisher's criterion
\(\mathbb E[\ln J_i]=\ln g\) places the critical point of the
infinite chain at \(g_c=1\)
\cite{fisher-prl,fisher-prb}.

For each \(L\in\{16,32,64,128\}\), we independently sample
\(R=500\) coupling realizations. A coupling realization is one
sampled vector \((J_1^{(r)},\ldots,J_{L-1}^{(r)})\), specifying the
couplings on all bonds of a chain. At each realization and each of
the \(201\) equally spaced fields in \([0.5,1.5]\), with
\(\Delta g=5\times10^{-3}\), we compute the even-parity ground state
using the free-fermion method described in
Section~\ref{subsec:app-model}; the method applies unchanged to
site-dependent couplings, meaning that the coupling \(J_i^{(r)}\)
may vary from bond to bond rather than being constant along the chain.
From this state we form the mutual-information
and correlator dissimilarities in
\eqref{eq:app-dmi} and \eqref{eq:app-dc}.

For each realization \(r\), the \(L-1\) finite \(H_0\) death
times, ordered as \(d_{r,1}(g)\geq\cdots\geq d_{r,L-1}(g)\), are the
edge weights of a minimum spanning tree of the complete graph on the
sites; we compute them with Prim's algorithm. Unlike in the clean
chain, this tree is not always the nearest-neighbour path. The fraction
of realization--field pairs \((r,g)\) for which it is that path
decreases with \(L\): for \(D^{\mathrm{MI}}\) it is \(0.60\) at
\(L=16\) and \(0.29\) at \(L=128\), while for \(D^C\) it is \(0.47\)
and \(0.16\), respectively. These computations use the general
closed formulas in Theorems~\ref{prop:closed} and~\ref{prop:mean-speed}, which do not require the tree to be a
nearest-neighbour path.

As in Section~\ref{subsec:app-numerics}, we approximate
\(\partial_g d_{r,k}\) by second-order centred finite differences of
the sorted death times, separately for each realization. These
derivatives give the individual velocity \(\xi_r\) through
Theorem~\ref{prop:closed}, and the three ensemble summaries
\(\xi_{\overline\lambda}\), \(\overline\xi\), and \(\xi_{\mathrm{rms}}\)
through \eqref{eq:mean-speed}--\eqref{eq:rms-speed} (equivalently,
through the death-time formulas
\eqref{eq:mean-overlap} and \eqref{eq:rms-closed}). We divide each
velocity by \(\sqrt{L-1}\), as in the clean-chain application.
The ratio \(\xi_{\overline\lambda}/\xi_{\mathrm{rms}}\) is already
dimensionless and, by \eqref{eq:velocity-variance}, lies in
\([0,1]\); it equals \(1\) precisely when all realization derivative
fields coincide.

To show uncertainty in this ratio, we use a nonparametric
bootstrap over the \(R\) realizations: we draw \(200\) resamples of
size \(R\) with replacement, recompute both velocities and their ratio
for each resample, and take the \(2.5\)th and \(97.5\)th percentiles
as a \(95\%\) band. Despite the disorder in coupling strengths, the same qualitative tendencies visible in the clean-chain case (Section~\ref{sec:application}) persist: the position of the maximum landscape velocity narrows around $g_c=1$ as system size increases, and its height grows with $L$, indicating that the landscape velocity seems to continue to probe quantum criticality even under random perturbations.

\begin{figure}[htbp]
\centering
\includegraphics[width=0.6\textwidth]{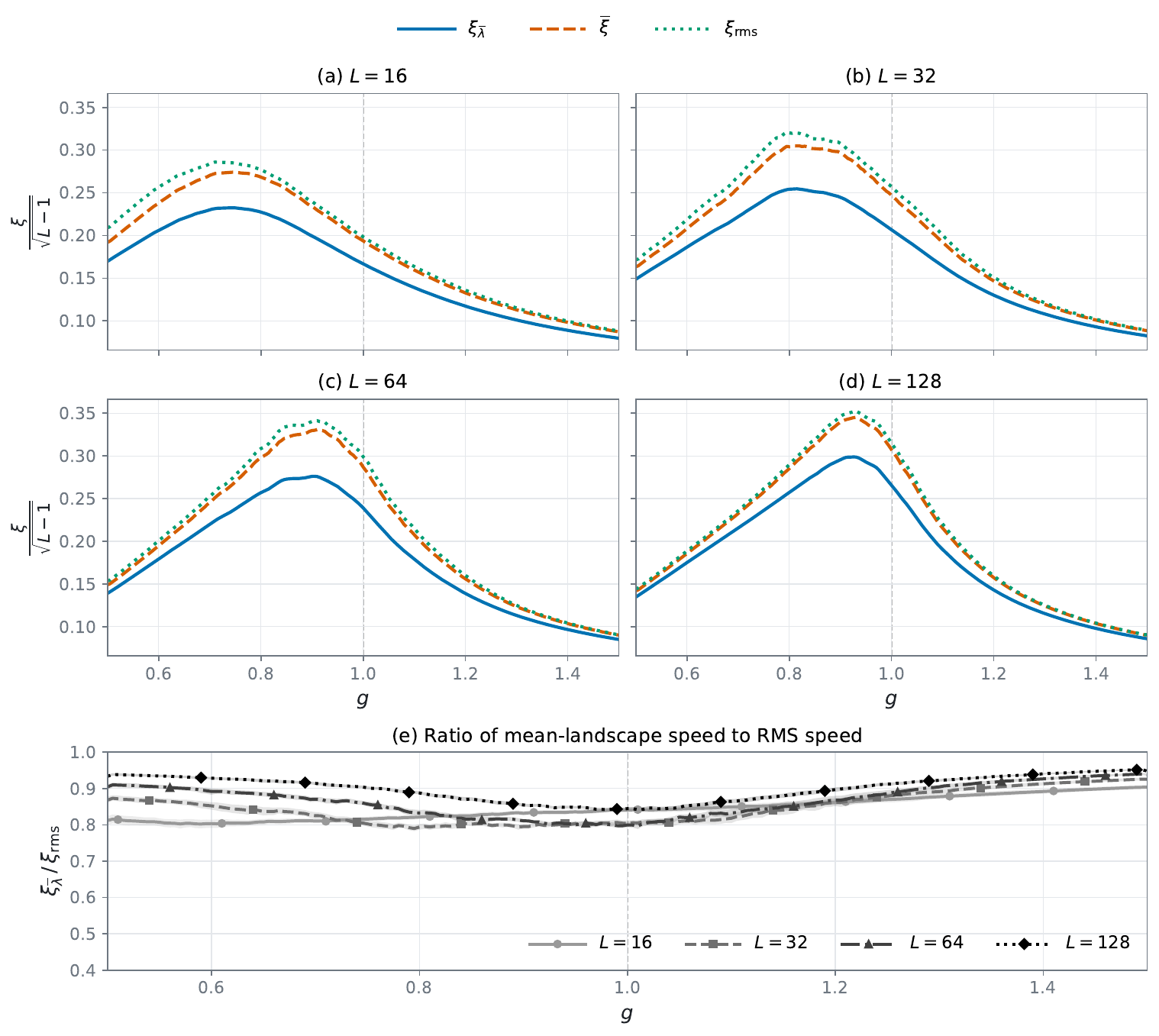}
\caption{\textbf{Empirical-mean and individual landscape velocities under disorder for the
random-bond open chain with \(D^{\mathrm{MI}}\).}
Panels (a)--(d) show \(L=16,32,64,128\), respectively. The solid
\textcolor{blue}{blue}, dashed \textcolor{orange}{orange}, and dotted
\textcolor{green!50!black}{green} curves are the velocity of the
empirical mean landscape, the mean individual velocity, and the root-mean-square
individual velocity, each divided by \(\sqrt{L-1}\). The vertical
dashed line marks \(g_c=1\). Panel (e) shows
\(\xi_{\overline\lambda}/\xi_{\mathrm{rms}}\) for all four sizes; the
shaded bands are the \(95\%\) percentile bootstrap intervals from
\(200\) resamples. Each size uses \(R=500\) independent disorder
realizations. Across panels (a)--(d), the maxima of the normalized
velocities appear to move towards \(g=1\) as \(L\) increases,
suggesting convergence to the critical field in the large-size limit.}
\label{fig:disorder-MI}
\end{figure}

\begin{figure}[htbp]
\centering
\includegraphics[width=0.6\textwidth]{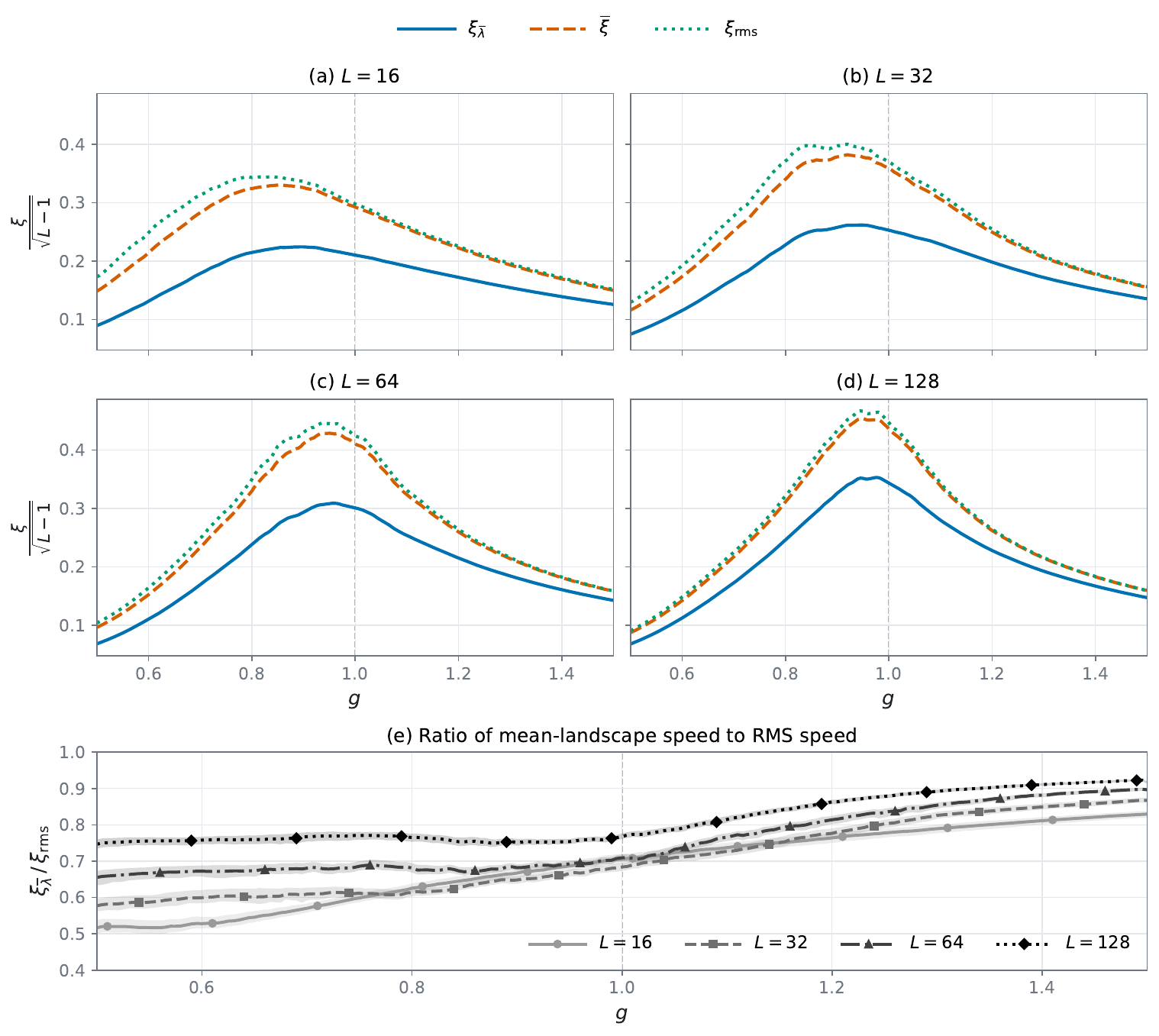}
\caption{\textbf{Empirical-mean and individual landscape velocities under disorder for the
random-bond open chain with \(D^C\).}
Panels (a)--(d) show \(L=16,32,64,128\), respectively. The solid
\textcolor{blue}{blue}, dashed \textcolor{orange}{orange}, and dotted
\textcolor{green!50!black}{green} curves are the velocity of the
empirical mean landscape, the mean individual velocity, and the root-mean-square
individual velocity, each divided by \(\sqrt{L-1}\). The vertical
dashed line marks \(g_c=1\). Panel (e) shows
\(\xi_{\overline\lambda}/\xi_{\mathrm{rms}}\) for all four sizes; the
shaded bands are the \(95\%\) percentile bootstrap intervals from
\(200\) resamples. Each size uses \(R=500\) independent disorder
realizations. Across panels (a)--(d), the maxima of the normalized
velocities appear to move towards \(g=1\) as \(L\) increases,
suggesting convergence to the critical field in the large-size limit.}
\label{fig:disorder-C}
\end{figure}

Figures~\ref{fig:disorder-MI} and
\ref{fig:disorder-C} show the results for \(D^{\mathrm{MI}}\) and
\(D^C\), respectively. In each figure, panels (a)--(d) show the three
normalized velocities for one system size per panel, on a common
vertical scale; panel (e) shows the ratio across all four sizes.

\clearpage
\section{Conclusion}
\label{sec:conclusion}
We have introduced the landscape velocity, a new topological summary that measures how rapidly persistence landscapes reorganize as a parameter varies.
Unlike summaries that operate on barcodes (which lack vector-space structure), we define velocity on landscapes. However, we do not require explicit landscape computation, the landscape velocity admits closed formulas that depend only on death times and their derivatives.
For the common-birth setting—the most common scenario in practice—the landscape velocity can be computed without ever constructing a landscape.

\textbf{Main theoretical contributions.} Theorem~\ref{prop:closed} establishes that for parameter-dependent families of finite barcodes born at filtration value zero, the squared landscape velocity equals half the sum of the squared derivatives of the death times, weighted by the death times themselves.
For Vietoris--Rips, \v{C}ech, and alpha filtrations, section~\ref{subsec:rips-mst} shows how the closed formula simplifies to depend only on minimum-spanning-tree edge weights—eliminating the need to compute the barcode itself. Proposition~\ref{prop:regularity} shows that sorted death times are differentiable almost everywhere when input dissimilarities are locally Lipschitz.
Theorem~\ref{prop:l2-speed} proves that this quantity is precisely the squared norm of the landscape derivative in the Hilbert space of landscapes, justifying the interpretation as metric speed.
Theorem~\ref{prop:mean-speed} extends the closed formula to empirical mean landscapes over multiple realizations, computing the velocity of the mean without ever constructing the mean landscape itself.

\textbf{Physical validation.} As an application to quantum criticality, we reanalyzed Olsthoorn's persistent homology study of the transverse-field Ising chain.
Our single rule—using only $H_0$, measuring landscape velocity, and locating the critical point at its maximum—succeeds for all four tested configurations (periodic and open boundary conditions, mutual-information and correlator dissimilarities), whereas Olsthoorn's methods require different summaries depending on boundary conditions.
The position of the maximum of the landscape velocity converges to the true critical point $g_c=1$ faster than Olsthoorn's estimators for all tested configurations. This convergence is mathematically proved for periodic boundaries with correlator dissimilarity (Section~\ref{subsec:app-log-proof}) under the nearest-neighbor minimum-spanning-tree assumption; for other configurations, comparable convergence is observed numerically. The peak height grows logarithmically with system size (proved for periodic boundaries with correlator dissimilarity; observed numerically for other configurations), indicating that the landscape velocity captures the singular topological response expected at quantum criticality.

\textbf{Implications for statistics and learning.} By defining the velocity through persistence landscapes rather than barcodes, we retain access to the rich statistical framework available in Hilbert spaces: laws of large numbers, central limit theorems, confidence bands, hypothesis tests, and machine-learning pipelines.
These require appropriate regularity of the landscape derivative and integrability of moments, as stated in the body.
The closed formulas ensure that this statistical machinery can be deployed without the computational overhead of explicitly constructing landscapes or their averages.
This is essential for applications to disordered systems, where one must average over many random realizations, and for noisy quantum devices, where correlators or mutual information can be estimated directly from measurements.

\textbf{Future directions.} Several extensions merit investigation.
\emph{Beyond common-birth.} Most naturally, one should develop closed formulas for barcodes where birth times also depend on the parameter, or for multiparameter families.
\emph{Stability.} While persistence barcodes are stable to perturbations of the input, the stability of the landscape velocity with respect to noise in the death times, and its finite-sample statistical properties, should be analyzed systematically.
\emph{Higher homology.} The application focused on $H_0$ because it is computationally cheaper and because the transverse-field Ising chain has no $H_1$ bar with open boundaries. Extending to higher homology degrees would broaden applicability.
\emph{Quantum applications.} The planned companion paper will explore applications to random quantum spin systems and to quantum criticality in other models, leveraging the disorder-averaging formulas of Section~\ref{subsec:disorder}.
\emph{Topological data analysis beyond physics.} The landscape velocity may capture dynamical phase transitions in time-dependent point clouds, identify abrupt changes in datasets parameterized by control variables, or improve machine-learning pipelines by replacing fixed summaries with adaptive ones sensitive to rate of change.

The central message is that parameter-dependent topological information, encoded in landscapes and measured through velocity, yields efficient closed formulas and statistical power without sacrificing computational efficiency.
This opens new avenues for using topological data analysis in quantum and classical systems where detecting phase transitions, critical points, or sudden reorganization is the scientific goal.

\begin{appendices}

\section{Proofs for analytical detection of the critical point
by normalized landscape velocity}
\label{subsec:landscape-velocity-proofs}
We now proof the results in Section~\ref{subsec:app-log-proof}.
In this section, to be explicit about the velocity $\xi$ being computed on a system with $L$ spins, we will write it as
$\xi=\xi_L$. We will then prove that, for the periodic chain with correlator
dissimilarity \(D^{\mathrm C}\), the normalized landscape velocity
\(\xi_L^{}/\sqrt{L-1}\) detects the critical point.
Under Assumption~\ref{hyp:app-mst}, the location of its maximum
converges to \(g_c=1\), while its maximum height grows
logarithmically with \(L\). Pfeuty's exact solution identifies
\(g_c=1\) as the quantum critical point \cite{pfeuty}. Thus the peak
location convergence establishes detection, and the logarithmic
growth shows that the peak also records the singular response.

We also establish the logarithmic divergence obtained by
first taking \(L\to\infty\) at fixed \(g\ne1\) and then
letting \(g\to1\).
We derive the logarithmic amplitude \(A_{\mathrm C}\),
that is, the coefficient multiplying the leading logarithmic
term, and show that the same coefficient governs both this
thermodynamic-limit divergence and the finite-size growth.
Consequently, \(A_{\mathrm C}\) is the asymptotic slope of
the normalized peak height plotted against \(\ln L\), as in
the bottom-right panel of Figure~\ref{fig:C1}.

\textbf{Throughout this subsection, \(L\) is even, the boundary
conditions are periodic, and we consider the ground state
in the even-parity sector,} i.e.\ the eigenspace of
\(P=\prod_{j=1}^{L}\sigma_j^z\) with eigenvalue \(+1\).
Under the Jordan--Wigner transformation, this corresponds
to even fermion-number parity and antiperiodic fermionic
boundary conditions, giving the positive momenta
\(k_n=(2n-1)\pi/L\), \(n=1,\ldots,L/2\).
We work with the Hamiltonian~\eqref{eq:app-hamiltonian} of
Section~\ref{sec:application},
\[
H_L(g)
=
-\sum_{j=1}^{L}\sigma_j^x\sigma_{j+1}^x
-g\sum_{j=1}^{L}\sigma_j^z,
\qquad
\sigma_{L+1}^{\alpha}=\sigma_1^{\alpha}.
\]

Let \(E_{0,L}(g)\) denote the lowest eigenvalue of \(H_L(g)\)
restricted to the even-parity sector, and write
\[
e_L(g)=\frac{E_{0,L}(g)}{L},
\qquad
C_L(g)=\langle\sigma_j^x\sigma_{j+1}^x\rangle,
\]
where \(e_L(g)\) is the ground-state energy per site and
the expectation defining \(C_L(g)\) is taken in that ground
state.
Translation invariance makes \(C_L(g)\) independent of \(j\).

\begin{assumptions}[Periodic-chain MST condition]
\label{hyp:app-mst}

For each even chain length \(L\) under consideration and each
field value \(g\) in the parameter range specified in the
relevant result, assume that the nearest-neighbour path
\(1-2-\cdots-L\) is a minimum spanning tree of the complete graph
with edge weights \(D_{ij}^{\mathrm C}(g)\). The minimum spanning
tree need not be unique.
\end{assumptions}

We numerically verified Assumption~\ref{hyp:app-mst} with periodic boundary conditions and the correlator dissimilarity \(D^{\mathrm{C}}\)
for \(g \in [0.5, 1.5]\) and \(L \in \{10, 12, 14, \ldots, 512\}\).

Translation invariance means that the model is unchanged by shifting all sites by one position with periodic boundary conditions. In the clean (undisordered) periodic Ising chain, this symmetry holds automatically and does not require numerical verification: it is a consequence of the translation symmetry of the Hamiltonian~\eqref{eq:app-hamiltonian} on the ring.
By translation invariance, all nearest-neighbour edges have
the same weight. Under Assumption~\ref{hyp:app-mst}, the finite
\(H_0\) death times therefore satisfy
\[
d_{1,L}^{}(g)
=
\cdots
=
d_{L-1,L}^{}(g)
=
D_{j,j+1}^{\mathrm C}(g)
=
1-C_L(g),
\qquad j=1,\dots,L-1.
\]
Here \(D_{j,j+1}^{\mathrm C}(g)\) is the dissimilarity between
sites \(j\) and \(j+1\). In \(d_{k,L}^{}(g)\), the first subscript
\(k\) indexes the sorted death times, while the second
subscript \(L\) records the chain length. 

\paragraph{From death times to the correlation response.}
Under 
Assumption~\ref{hyp:app-mst}
the finite death times coincide.
The closed formula \eqref{eq:app-closed} therefore gives
\begin{equation}
\label{eq:app-log-normalized-velocity}
\frac{\xi_L^{}(g)}{\sqrt{L-1}}
=
\left|\partial_g d_{1,L}^{}(g)\right|
\sqrt{\frac{d_{1,L}^{}(g)}{2}}.
\end{equation}
Since \(d_{1,L}^{}(g)=1-C_L(g)\), the normalized
landscape velocity is the magnitude of the field derivative
of the nearest-neighbour correlation, multiplied by
\(\sqrt{d_{1,L}^{}(g)/2}\).

The landscape-velocity results below are conditional on

Assumption~\ref{hyp:app-mst}.
We verified its condition numerically for the periodic chain with
correlator dissimilarity.
The corresponding statements about the energy and correlation
responses do not require Assumption~\ref{hyp:app-mst}.

The energy identities and logarithmic energy-response
asymptotics recalled below are classical consequences of
the exact solution of the transverse-field Ising chain
\cite{pfeuty}. We include their derivations for completeness,
to fix the normalization and the coefficient of the
logarithmic term, and to make the argument accessible
without assuming familiarity with the model.
Lemma~\ref{lem:app-energy-response} is a direct application
of the Hellmann--Feynman identity, while
Proposition~\ref{prop:app-log-infinite-energy} follows
by estimating the singular integral obtained from the
exact ground-state energy.
These classical results provide the physical input for our
new conclusions concerning the normalized landscape velocity:
its thermodynamic-limit logarithmic divergence
(Corollary~\ref{cor:app-log-infinite-velocity}),
its finite-size growth at the critical field
(Proposition~\ref{prop:app-log-critical-response}),
and the asymptotic growth of its maximum together with
the convergence of its maximizers to \(g_c=1\)
(Proposition~\ref{prop:app-log-peak}).

\begin{lemma}[Energy-response identity]
\label{lem:app-energy-response}
For the periodic chain in the even-parity ground state, for even values of $L$ and for \(g>0\),
$$
\frac{\xi_L^{}(g)}{\sqrt{L-1}}
=
g\sqrt{\frac{d_{1,L}^{}(g)}2}
\bigl(-\partial_g^2e_L(g)\bigr).
$$
\end{lemma}

\begin{proof}
Set
\[
Z_L(g)=\langle\sigma_j^z\rangle.
\]
Taking the expectation of the Hamiltonian and dividing by
\(L\) gives
\[
e_L(g)=-C_L(g)-gZ_L(g).
\]
In this finite ferromagnetic chain with an \textbf{even number of sites} $L$, the selected ground state in the fixed even-parity sector is nondegenerate and real-analytic in $g>0$, \textbf{justifying eigenvalue and eigenvector differentiation.}
The Hellmann--Feynman identity, applied within the fixed
even-parity sector, gives
\[
\partial_g e_L(g)
=
-\frac1L\sum_{j=1}^{L}\langle\sigma_j^z\rangle
=
-Z_L(g).
\]
Substitution yields
\[
C_L(g)=-e_L(g)+g\,\partial_g e_L(g).
\]
Differentiation gives
\[
\partial_g C_L(g)
=
-\partial_g e_L(g)
+\partial_g e_L(g)
+g\,\partial_g^2e_L(g)
=
g\,\partial_g^2e_L(g).
\]
Using 
\(d_{1,L}^{}=1-C_L\) leads to
\[
\partial_g d_{1,L}^{}(g)
=-g\,\partial_g^2e_L(g).
\]
Substituting in 
\eqref{eq:app-log-normalized-velocity}, for \(g>0\), we obtain the result.
\end{proof}

\paragraph{Mode representation.}
The Jordan--Wigner transformation and Bogoliubov
diagonalization give the fermionic dispersion
\begin{equation}
\label{eq:fermionic_dispersion}
    \omega_k(g)
    =
    \sqrt{1+g^2-2g\cos k}
\end{equation}
and, in the even-parity sector of the periodic chain,
\[
e_L(g)
=
-\frac{2}{L}\sum_{n=1}^{L/2}\omega_{k_n}(g),
\qquad
k_n=\frac{(2n-1)\pi}{L}.
\]
In the thermodynamic limit,
\[
e_\infty(g)
=
-\frac1\pi\int_0^\pi\omega_k(g)\,dk
\]
\cite{pfeuty}.
The sum contains only the positive momenta; the factor \(2\)
accounts for the corresponding negative momenta.

Direct differentiation gives
\[
\partial_g\omega_k(g)
=
\frac{g-\cos k}{\omega_k(g)},
\qquad
\partial_g^2\omega_k(g)
=
\frac{\sin^2k}{\omega_k(g)^3}.
\]
It follows that
\begin{equation}
\label{eq:app-log-finite-response}
-\partial_g^2e_L(g)
=
\frac2L\sum_{n=1}^{L/2}
\frac{\sin^2k_n}{\omega_{k_n}(g)^3}
>0.
\end{equation}
For fixed \(g\ne1\), differentiation under the integral
is justified because the dispersion is bounded away from
zero. Hence
\begin{equation}
\label{eq:app-log-infinite-response}
-\partial_g^2e_\infty(g)
=
\frac1\pi\int_0^\pi
\frac{\sin^2k}{\omega_k(g)^3}\,dk.
\end{equation}

\begin{proposition}[Logarithmic energy response]
\label{prop:app-log-infinite-energy}

For the periodic chain in the even-parity sector, as \(g\to1\)
with \(g\ne1\),
\[
-\partial_g^2e_\infty(g)
=
\frac1\pi\ln\frac1{|g-1|}
+O(1).
\]
\end{proposition}

\begin{proof}
Write \(\delta=g-1\). Then \eqref{eq:fermionic_dispersion} gives
\[
\omega_k(g)^2
=
\delta^2+2g(1-\cos k).
\]
Fix a sufficiently small \(k_0>0\), independent of
\(\delta\), and split the integral in
\eqref{eq:app-log-infinite-response} at \(k_0\).

For \(k\in[k_0,\pi]\) and \(g\) in a fixed sufficiently
small neighbourhood of \(1\), the denominator is bounded
away from zero. This part of the integral is therefore
uniformly bounded.

For \(0\le k\le k_0\), Taylor expansion gives
\[
\sin^2k=k^2+O(k^4),
\qquad
2g(1-\cos k)=gk^2+O(k^4),
\]
with uniform constants for \(g\) near \(1\).
Moreover, \(\delta^2+2g(1-\cos k)\) and
\(\delta^2+gk^2\) are uniformly comparable.
The mean value theorem applied to \(x^{-3/2}\) therefore
shows that
\[
\left|
\frac{\sin^2k}
     {(\delta^2+2g(1-\cos k))^{3/2}}
-
\frac{k^2}{(\delta^2+gk^2)^{3/2}}
\right|
\le
C\left(
\frac{k^4}{(\delta^2+k^2)^{3/2}}
+
\frac{k^6}{(\delta^2+k^2)^{5/2}}
\right)
\le C'k.
\]
Thus replacing the original integrand by the comparison
kernel changes the integral by a uniformly bounded amount.

The comparison integral is explicit:
\[
\begin{aligned}
\int_0^{k_0}
\frac{k^2}{(\delta^2+gk^2)^{3/2}}\,dk
&=
g^{-3/2}
\left[
\operatorname{arcsinh}
\left(\frac{\sqrt g\,k_0}{|\delta|}\right)
-
\frac{\sqrt g\,k_0}
     {\sqrt{\delta^2+gk_0^2}}
\right]
\\
&=
g^{-3/2}\ln\frac1{|\delta|}+O(1).
\end{aligned}
\]
Since \(g^{-3/2}=1+O(\delta)\) and
\(\delta\ln(1/|\delta|)\to0\), this equals
\[
\ln\frac1{|\delta|}+O(1).
\]
Multiplication by \(1/\pi\) proves the result.
\end{proof}

The source of the logarithm is the small-momentum region:
at criticality the integrand behaves as \(1/k\).
For \(g\ne1\), the scale \(|g-1|\) cuts off this singular
behaviour.

We now restate Corollary~\ref{cor:app-log-infinite-velocity}, concerning
the divergence of the maximum height in the thermodynamic limit.

\thermolimitcor*

\begin{proof}
For fixed \(g\ne1\), the finite sums for the energy and its
first two derivatives converge to their integral
representations. Consequently, Lemma~\ref{lem:app-energy-response} implies
\[
\lim_{L\to\infty}
\frac{\xi_L^{}(g)}{\sqrt{L-1}}
 =
g\sqrt{\frac{d_\infty^{}(g)}2}
\bigl(-\partial_g^2e_\infty(g)\bigr),
\]
where
\[
d_\infty^{}(g)
=
1+e_\infty(g)-g\,\partial_g e_\infty(g).
\]

The energy-response identity also gives
\[
C_\infty(g)
=
\frac1\pi\int_0^\pi
\frac{1-g\cos k}{\omega_k(g)}\,dk.
\]
Since
\[
\omega_k(g)^2-(1-g\cos k)^2
=
g^2\sin^2k,
\]
the absolute value of this integrand is bounded by \(1\).
Dominated convergence therefore applies at \(g=1\).
At that point,
\[
\omega_k(1)=2\sin(k/2),
\qquad
\frac{1-\cos k}{\omega_k(1)}
=
\sin(k/2),
\]
and hence
\[
C_\infty(1)
=
\frac1\pi\int_0^\pi\sin(k/2)\,dk
=
\frac2\pi.
\]
Thus
\[
d_\infty^{}(1)=1-\frac2\pi.
\]
The prefactor
\(g\sqrt{d_\infty^{}(g)/2}\) therefore tends to
\(\sqrt{1/2-1/\pi}\).
The conclusion follows from
Proposition~\ref{prop:app-log-infinite-energy}.
\end{proof}

We now restate and prove Proposition~\ref{prop:app-log-critical-response},
which describes finite-size growth at the critical field.

\criticalfieldprop*

\begin{proof}
We want to express $\frac{\xi_L^{}(1)}{\sqrt{L-1}}$
using Lemma~\ref{lem:app-energy-response}, so we will first compute 
$-\partial_g^2e_L(1)$ and $d_{1,L}^{}(1)$.

In the even-parity sector, we can use \eqref{eq:app-log-finite-response}.
At \(g=1\), the summand becomes
\[
F(k)
=
\frac{\sin^2k}{(2\sin(k/2))^3}
=
\frac{\cos^2(k/2)}{2\sin(k/2)}.
\]
As $k\to0^+$,
\begin{equation}
\label{eq:app-critical-integrand-asymptotic}
F(k)=\frac{1}{k}+O(k).
\end{equation}
Thus \(F(k)-1/k\) extends to a bounded function on
\([0,\pi]\), and
\[
\begin{aligned}
-\partial_g^2e_L(1)
&=
\frac2L\sum_{n=1}^{L/2}\frac1{k_n}+O(1)
\\
&=
\frac2\pi\sum_{n=1}^{L/2}\frac1{2n-1}+O(1)
\\
&=
\frac1\pi\ln L+O(1).
\end{aligned}
\]

The correlation at criticality is
\[
C_L(1)
=
\frac2L\sum_{n=1}^{L/2}\sin(k_n/2)
=
\frac1{L\sin(\pi/(2L))}.
\]
Consequently,
\[
d_{1,L}^{}(1)
=
1-C_L(1)
=
1-\frac1{L\sin(\pi/(2L))}
=
1-\frac2\pi+O(L^{-2}).
\]
Substitution in
Lemma~\ref{lem:app-energy-response} proves the
landscape-velocity statement.
\end{proof}

This calculation concerns the value at the fixed field
\(g=1\). Since the maximum is at least this value, it gives
a lower bound for the peak height. A matching upper bound
is needed to establish the asymptotic of the maximum itself.

We now restate Proposition~\ref{prop:app-log-peak}, concerning
the finite-size growth of the maximum.

\normalizedpeakprop*

\begin{proof}
All limits in \(L\) are taken through the even integers
under consideration.
Define
\[
H_L
=
\max_{g\in I}
\frac{\xi_L(g)}{\sqrt{L-1}},
\qquad
R_L(g)=-\partial_g^2e_L(g),
\qquad
B_L(g)=g\sqrt{\frac{d_{1,L}(g)}2}.
\]
By Lemma~\ref{lem:app-energy-response},
\[
\frac{\xi_L(g)}{\sqrt{L-1}}
=
B_L(g)R_L(g).
\]
Both factors are nonnegative.

We first establish a uniform upper bound for the response.
Fix \(k_0\in(0,\pi/2)\).
Since
\[
\omega_k(g)^2
=
(g-\cos k)^2+\sin^2k
\ge \sin^2k,
\]
for \(0<k\le k_0\) we have
\[
\frac{\sin^2k}{\omega_k(g)^3}
\le \frac1{\sin k}
=
\frac1k+O(k).
\]
On the compact set \(I\times[k_0,\pi]\), the dispersion
is bounded away from zero: its only zero for \(g>0\)
and \(k\in[0,\pi]\) occurs at \((g,k)=(1,0)\).
Using \eqref{eq:app-log-finite-response} and the same
odd harmonic sum as in
Proposition~\ref{prop:app-log-critical-response}, we obtain
\[
0\le R_L(g)\le \frac1\pi\ln L+K,
\qquad g\in I,
\]
where \(K\) is independent of \(L\) and \(g\).

We next control the prefactor near \(g=1\).
The finite-size correlation is
\[
C_L(g)
=
\frac2L\sum_{n=1}^{L/2}
\frac{1-g\cos k_n}{\omega_{k_n}(g)},
\qquad
k_n=\frac{(2n-1)\pi}{L}.
\]
Its summand is bounded in absolute value by \(1\), since
\[
\omega_k(g)^2-(1-g\cos k)^2
=
g^2\sin^2k\ge0.
\]
In particular, \(0\le d_{1,L}(g)=1-C_L(g)\le2\), so
\(B_L(g)\le b\) on \(I\).

Let \(g_L\in I\) satisfy \(g_L\to1\).
For a fixed \(\eta\in(0,\pi)\), split the correlation sum
into \(k_n<\eta\) and \(k_n\ge\eta\).
The absolute value of the first contribution is at most
\[
\frac2L\#\{n:k_n<\eta\}
=
\frac{\eta}{\pi}+O(L^{-1}).
\]
On \([\eta,\pi]\), the summands converge uniformly to
their critical-field values,
\[
\frac{1-\cos k}{\omega_k(1)}
=
\sin(k/2).
\]
The corresponding Riemann sums therefore converge to
\[
\frac1\pi\int_\eta^\pi\sin(k/2)\,dk.
\]
Letting \(\eta\downarrow0\) gives
\[
C_L(g_L)\longrightarrow\frac2\pi,
\qquad
B_L(g_L)\longrightarrow
B_\star:=\sqrt{\frac12-\frac1\pi}.
\]
The same argument applies along every subsequence of
even lengths.

We now justify the required shrinking-neighbourhood
bound. Set
\[
M_L(\varepsilon)
=
\sup_{\substack{g\in I\\|g-1|\le\varepsilon}}
B_L(g).
\]
We claim that
\[
\lim_{\varepsilon\downarrow0}
\limsup_{L\to\infty}M_L(\varepsilon)
=
B_\star.
\]
The lower bound follows from
\(M_L(\varepsilon)\ge B_L(1)\to B_\star\).
For the upper bound, suppose the displayed limit were
strictly larger than \(B_\star\).
Then there would exist \(\delta>0\), a sequence
\(\varepsilon_j\downarrow0\), and increasing even
lengths \(L_j\to\infty\) such that
\[
M_{L_j}(\varepsilon_j)>B_\star+2\delta.
\]
By the definition of the supremum, we can choose
\(g_j\in I\) such that
\[
|g_j-1|\le\varepsilon_j,
\qquad
B_{L_j}(g_j)>
M_{L_j}(\varepsilon_j)-\delta
>B_\star+\delta.
\]
But \(g_j\to1\), contradicting the sequential convergence
proved above. This establishes the claim.

To bound the peak height, fix
\[
0<\varepsilon<\min\{1-a,b-1\}
\]
and split the maximization into
\[
H_L
=
\max\{H_L^{\mathrm{near}}(\varepsilon),
       H_L^{\mathrm{far}}(\varepsilon)\},
\]
where
\[
H_L^{\mathrm{near}}(\varepsilon)
=
\sup_{\substack{g\in I\\|g-1|\le\varepsilon}}
B_L(g)R_L(g),
\qquad
H_L^{\mathrm{far}}(\varepsilon)
=
\sup_{\substack{g\in I\\|g-1|>\varepsilon}}
B_L(g)R_L(g).
\]
Since
\[
\omega_k(g)^2
=
(g-1)^2+2g(1-\cos k)
\ge(g-1)^2,
\]
the response formula gives
\[
R_L(g)\le\varepsilon^{-3}
\qquad\text{when }|g-1|>\varepsilon.
\]
Together with \(B_L(g)\le b\), this implies
\[
H_L^{\mathrm{far}}(\varepsilon)
\le b\varepsilon^{-3}.
\]
On the other hand, the uniform logarithmic estimate gives
\[
H_L^{\mathrm{near}}(\varepsilon)
\le
M_L(\varepsilon)
\left(\frac1\pi\ln L+K\right).
\]
Hence, for each fixed \(\varepsilon\),
\[
\limsup_{L\to\infty}\frac{H_L}{\ln L}
\le
\frac1\pi
\limsup_{L\to\infty}M_L(\varepsilon).
\]
Letting \(\varepsilon\downarrow0\), we obtain
\[
\limsup_{L\to\infty}\frac{H_L}{\ln L}
\le
\frac{B_\star}{\pi}
=
A_{\mathrm C}.
\]
On the other hand,
\[
H_L\ge\frac{\xi_L(1)}{\sqrt{L-1}},
\]
so Proposition~\ref{prop:app-log-critical-response} yields
\[
\liminf_{L\to\infty}\frac{H_L}{\ln L}
\ge A_{\mathrm C}.
\]
Therefore
\[
A_{\mathrm C}
\le
\liminf_{L\to\infty}\frac{H_L}{\ln L}
\le
\limsup_{L\to\infty}\frac{H_L}{\ln L}
\le
A_{\mathrm C},
\]
so all inequalities become equalities and 
\[
H_L\sim A_{\mathrm C}\ln L.
\]

Finally, for every fixed \(\varepsilon>0\), the normalized
velocity is uniformly bounded on
\(\{g\in I:|g-1|\ge\varepsilon\}\), whereas its value at
\(g=1\) diverges.
Thus every maximizer eventually satisfies
\(|g_{\mathrm{max},L}-1|<\varepsilon\).
Since \(\varepsilon\) is arbitrary, all choices of
maximizers satisfy \(g_{\mathrm{max},L}\to1\).
\end{proof}

\paragraph{Interpretation.}
Both logarithmic laws reflect the same small-momentum
singularity: at the critical field, the energy-response
integrand behaves as \(1/k\), as expressed in
\eqref{eq:app-critical-integrand-asymptotic}.
In the infinite system, moving away from criticality
regularizes this singularity at the momentum scale
\(\lvert g-1\rvert\).
At finite size and at criticality, the smallest allowed
momentum is \(k_1=\pi/L\), and the corresponding discrete
sum grows as \(\ln L\).
The preceding propositions establish that the same
amplitude \(A_{\mathrm C}\) governs the thermodynamic-limit
divergence and the normalized finite-size peak growth.

We do not determine the convergence rate of the peak location
or the additive constant in the peak-height asymptotics.

\paragraph{Scope.}
The proof above establishes the logarithmic coefficient only for
periodic boundaries and the correlator dissimilarity \(D^C\).
The analysis also relies on translation invariance and
the periodic-chain momentum grid.
It therefore does not establish the open-chain scaling.
The open-boundary results in the figures remain numerical
observations unless supported by a separate argument.

\section{Notation and Dependencies}
\label{app:notation}

This appendix provides a comprehensive notation table and a conceptual dependency graph for the main definitions, hypotheses, and results.

\subsection{Notation Table}
\begin{longtable}{@{}p{0.29\textwidth}p{0.68\textwidth}@{}}
\toprule
\textbf{Symbol} & \textbf{Meaning and First Reference} \\
\midrule
\endhead

\(g\) & External parameter indexing families of data and barcodes (e.g., transverse field in Section \ref{sec:application}) \\
\(t\) or \(\epsilon\) & Filtration threshold in persistent homology \\
\(I \subseteq \mathbb{R}\) & Open interval of parameter values in \(g\) \\

\(\mathcal{K} = \{K(t)\}_{t \in T}\) & Filtration of simplicial complexes (Definition \ref{def:filtration}) \\
\(H_p(-;\mathbb{F})\) & Simplicial homology in degree $p$ with coefficients in field $\mathbb{F}$ \\
\(\mathcal{B}_g^{\mathrm{fin}} = \{[b_m(g), d_m(g))\}_{m=1}^{M(g)}\) & Finite bars of the persistence barcode $\mathcal{B}_g$ (Definition \ref{def:filtration}) \\
\(b_m(g)\) & Birth time of bar $m$ as a function of parameter $g$ \\
\(d_m(g)\) & Death time of bar $m$ as a function of parameter $g$ \\
\(M(g)\) & Number of finite bars at parameter value $g$ \\

\(\Lambda_{(b,d)}(t)\) & Tent function for bar $[b,d)$, which is $\max\{0, \min\{t-b, d-t\}\}$ (Definition \ref{def:landscape}) \\
\(\lambda_k(g;t)\) & $k$-th largest tent value at filtration threshold $t$ for parameter $g$ \\
\(\lambda(g) = (\lambda_k(g;\cdot))_{k \geq 1}\) & Persistence landscape as sequence of layers \\

\(\mathcal{H} = \ell^2(\mathbb{N}; L^2(\mathbb{R};\mathbb{R}))\) & Hilbert space of landscapes (Equation \ref{eq:hilbert}) \\
\(\langle \cdot, \cdot \rangle_{\mathcal{H}}\) & Inner product in $\mathcal{H}$ (Equation \ref{eq:hilbert}) \\
\(\|\cdot\|_{\mathcal{H}}\) & Norm in $\mathcal{H}$ (Equation \ref{eq:hilbert}) \\

\(\partial_g \lambda_k(g;t)\) & Pointwise derivative of landscape layer $\lambda_k$ with respect to $g$ \\
\(\partial_g \lambda(g)\) & Vector-valued derivative of landscape in $\mathcal{H}$ \\
\(\xi(g)\) & Landscape velocity, defined pointwise (Definition \ref{def:velocity}) and equal to Hilbert-space speed where differentiable (Theorem \ref{prop:l2-speed}) \\
\(\xi_L(g)\) & Same quantity as \(\xi(g)\), with the
dependence on the point-cloud size \(L\) made explicit.
Notation used only in
Section~\ref{subsec:app-log-proof} \\
\(\xi(g)/\sqrt{L-1}\) & Normalized (per-bar) landscape velocity, where \(L-1\) is the number of finite bars in each barcode \\
\(\xi_r(g)\) & Landscape velocity for realization $r$ in disorder setting. Notation used only in Section~\ref{subsec:disorder} \\
\(\xi_{\overline{\lambda}}(g)\) & Velocity of empirical mean landscape (Equation \ref{eq:mean-speed}) \\
\(\overline{\xi}(g)\) & Mean of individual landscape velocities (Equation \ref{eq:mean-individual-speed}) \\
\(\xi_{\mathrm{rms}}(g)\) & Root-mean-square of individual velocities (Equation \ref{eq:rms-speed}) \\

\(L\) & Chain length (number of qubits) in quantum application \\
\(H(g)\) & Transverse-field Ising Hamiltonian (Equation \ref{eq:app-hamiltonian}) \\
\(g = h/J\) & Dimensionless field ratio \\
\(g_c = 1\) & Critical field value for transverse-field Ising chain \\

\(D_{ij}(g)\) & Dissimilarity between vertices $i$ and $j$ \\
\(D^{\mathrm{MI}}_{ij}(g) = 1 - \frac{M_{ij}(g)}{2\ln 2}\) & Mutual-information dissimilarity (Equation \ref{eq:app-dmi}) \\
\(D^{C}_{ij}(g) = 1 - \langle \sigma_i^x \sigma_j^x \rangle(g)\) & Correlator dissimilarity (Equation \ref{eq:app-dc}) \\
\(\mathrm{VR}_t(X, D_g)\) & Vietoris--Rips filtration at parameter $g$ \\

\(\rho_i\) & One-site reduced density matrix \\
\(\rho_{ij}\) & Two-site reduced density matrix \\
\(M_{ij}\) & Von Neumann mutual information (Section \ref{subsec:app-model}) \\
\(C_L(g) = \langle \sigma_j^x \sigma_{j+1}^x \rangle\) & Nearest-neighbor correlator at chain length $L$ \\

\(e_L(g) = \frac{E_{0,L}(g)}{L}\) & Ground-state energy per site \\
\(A_{\mathrm{C}}\) & Logarithmic amplitude of normalized landscape velocity (Section \ref{subsec:app-log-proof}), 
which is the slope in the bottom-right panel of Figure~\ref{fig:C1} \\
\(g_{\mathrm{max},L}\) & Position of maximum of landscape velocity for chain length $L$ \\

\bottomrule
\end{longtable}

\subsection{Conceptual Dependency Graph}

The following diagram illustrates how definitions, hypotheses, and main results depend on one another. Arrows indicate logical dependencies: an arrow from $A$ to $B$ means that $B$ uses or builds on $A$.

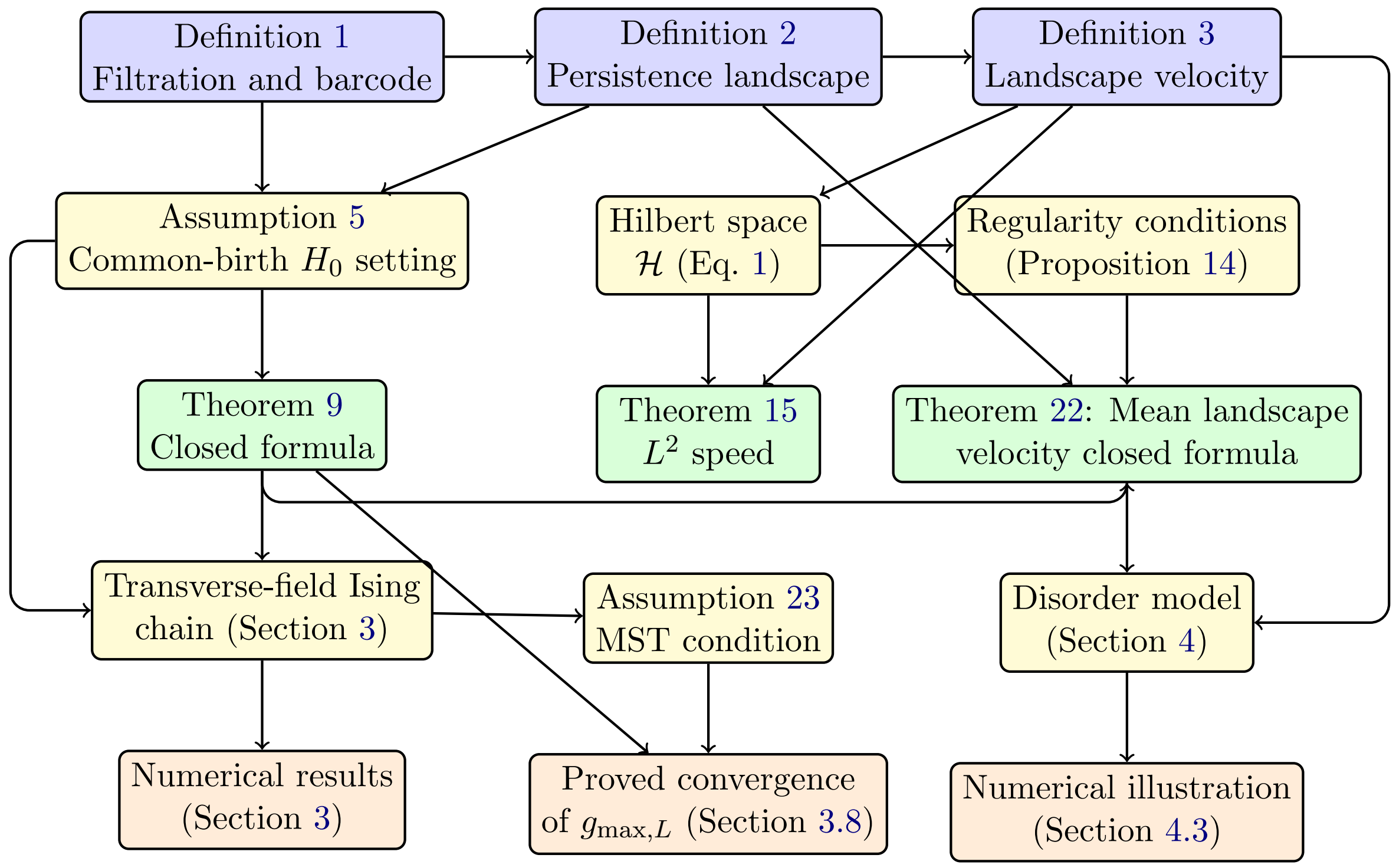
\begin{figure}[h]
\centering
\begin{tikzpicture}[
  node distance=1.cm and 1.cm,
  thick,
  every node/.style={
    draw,
    align=center,
    rounded corners=3pt,
    minimum width=2.5cm
  },
  definition/.style={fill=blue!15},
  assumption/.style={fill=yellow!20},
  proposition/.style={fill=green!15},
  application/.style={fill=orange!15}
]

\node[definition] (filtration) {Definition~\ref{def:filtration}\\Filtration and barcode};
\node[definition, right=of filtration] (landscape) {Definition~\ref{def:landscape}\\Persistence landscape};
\node[definition, right=of landscape] (velocity) {Definition~\ref{def:velocity}\\Landscape velocity};

\node[assumption, below=of filtration] (common_birth) {Assumption~\ref{ass:common-birth}\\Common-birth $H_0$ setting};
\node[assumption, below=of landscape] (hilbert) {Hilbert space\\$\mathcal{H}$ (Eq.~\ref{eq:hilbert})};
\node[assumption, below=of velocity] (regularity) {Regularity conditions\\(Proposition~\ref{prop:regularity})};

\node[proposition, below=of common_birth] (closed) {Theorem~\ref{prop:closed}\\Closed formula};
\node[proposition, below=of hilbert] (l2speed) {Theorem~\ref{prop:l2-speed}\\$L^2$ speed};
\node[proposition, below=of regularity] (mean) {Theorem~\ref{prop:mean-speed}: Mean landscape\\ velocity closed formula};

\node[assumption, below=of closed] (ising) {Transverse-field Ising\\chain (Section~\ref{sec:application})};
\node[assumption, below=of l2speed] (mst) {Assumption~\ref{hyp:app-mst}\\MST condition};
\node[assumption, below=of mean] (disorder) {Disorder model\\(Section~\ref{subsec:disorder})};

\node[application, below=of ising] (detection) {Numerical results\\(Section~\ref{sec:application})};
\node[application, below=of mst] (proof) {Proved convergence\\ of $g_{\max,L}$ (Section~\ref{subsec:app-log-proof}) };
\node[application, below=of disorder] (scaling) {Numerical illustration\\(Section~\ref{subsec:disorder-numerical})};

\draw[->] (filtration) to (landscape);
\draw[->] (landscape) to (velocity);
\draw[->] (velocity) to (hilbert);
\draw[->] (filtration) to (common_birth);
\draw[->] (landscape) to (common_birth);
\draw[->] (hilbert) to (regularity);
\draw[->] (common_birth) to (closed);
\draw[->] (velocity) to (l2speed);
\draw[->] (hilbert) to (l2speed);
\draw[->, rounded corners=6pt]
  (closed.south)
  -- ++(0,-0.35cm)
  -| (mean.south);
\draw[->] (landscape) to (mean);
\draw[->] (regularity) to (mean);
\draw[->] (closed) to (ising);
\draw[->, rounded corners=6pt]
  (common_birth.west)
  -- ++(-0.5cm,0)
  |- (ising.west);
\draw[->] (ising) to (mst);
\draw[->] (mean) to (disorder);
\draw[->, rounded corners=6pt]
  (velocity.east)
  -- ++(1.2cm,0)
  |- (disorder.east);
\draw[->] (ising) to (detection);
\draw[->] (mst) to (proof);
\draw[->] (closed) to (proof);
\draw[->] (disorder) to (scaling);

\end{tikzpicture}
\caption{Conceptual dependency structure. Definitions (blue) and assumptions (yellow) feed into propositions (green), which enable applications (orange).}
\label{fig:dependency}
\end{figure}

\textbf{Reading guide.} To understand the theoretical core, one should read:
\emph{Foundations:} Definitions \ref{def:filtration}--\ref{def:velocity} and Assumption \ref{ass:common-birth}.
\emph{Main results:} Theorems \ref{prop:closed}, \ref{prop:l2-speed}, and \ref{prop:mean-speed}.
\emph{Regularity:} Proposition \ref{prop:regularity} justifies the differentiability assumption in Theorem \ref{prop:closed}.
\emph{Application:} Sections \ref{sec:application} and \ref{subsec:app-log-proof} apply these results to quantum criticality.
\emph{Extensions:} Section \ref{subsec:disorder} extends the closed formula to disorder-averaged landscapes.

\end{appendices}

\end{document}